\documentclass[10pt]{article}
\usepackage{mathtools}
\mathtoolsset{showonlyrefs}
\usepackage{subfigure}
\usepackage{amssymb,amsmath,amsfonts,bbm,pifont,upgreek,bbold,accents,ulem} 
\usepackage[colorlinks=true]{hyperref}
\hypersetup{urlcolor=blue, citecolor=red}
\DeclareFontEncoding{FMS}{}{}
\DeclareFontSubstitution{FMS}{futm}{m}{n}
\DeclareFontEncoding{FMX}{}{}
\DeclareFontSubstitution{FMX}{futm}{m}{n}
\DeclareSymbolFont{fouriersymbols}{FMS}{futm}{m}{n}
\DeclareSymbolFont{fourierlargesymbols}{FMX}{futm}{m}{n}
\DeclareMathDelimiter{\VERT}{\mathord}{fouriersymbols}{152}{fourierlargesymbols}{147}
\newtheorem{theorem}{Theorem}[section]
\newtheorem{definition}{Definition}[section]
\newtheorem{proposition}{Proposition}[section]
\newtheorem{lemma}{Lemma}[section]
\newtheorem{remark}{Remark}[section]
\newtheorem{corollary}{Corollary}[section]
\renewcommand\b{{\beta}}
\renewcommand\d{{\delta}}
\expandafter\chardef\csname pre amssym.def
at\endcsname=\the\catcode`\@
\catcode`\@=11
\def\undefine#1{\let#1\undefined}
\def\newsymbol#1#2#3#4#5{\let\next@\relax
 \ifnum#2=\@ne\let\next@\msafam@\else
 \ifnum#2=\tw@\let\next@\msbfam@\fi\fi
 \mathchardef#1="#3\next@#4#5}
\def\mathhexbox@#1#2#3{\relax
 \ifmmode\mathpalette{}{\m@th\mathchar"#1#2#3}%
 \else\leavevmode\hbox{$\m@th\mathchar"#1#2#3$}\fi}
\def\hexnumber@#1{\ifcase#1 0\or 1\or 2\or 3\or 4\or 5\or 6\or 7\or
8\or
 9\or A\or B\or C\or D\or E\or F\fi}
\ifcase\@ptsize
 \font\tenmsb=msbm10
 \font\sevenmsb=msbm7
 \font\fivemsb=msbm5
\or
 \font\tenmsb=msbm10 scaled \magstephalf
 \font\sevenmsb=msbm7 scaled \magstephalf
 \font\fivemsb=msbm5  scaled \magstephalf
\or
 \font\tenmsb=msbm10 scaled \magstep1
 \font\sevenmsb=msbm7 scaled \magstep1
 \font\fivemsb=msbm5 scaled \magstep1
\fi
\newfam\msbfam
\textfont\msbfam=\tenmsb
\scriptfont\msbfam=\sevenmsb
\scriptscriptfont\msbfam=\fivemsb
\edef\msbfam@{\hexnumber@\msbfam}
\def\Bbb#1{\fam\msbfam\relax#1}
\def\widehat#1{\setboxz@h{$\m@th#1$}%
 \ifdim\wdz@>\tw@ em\mathaccent"0\msbfam@5B{#1}%
 \else\mathaccent"0362{#1}\fi}
\def\widetilde#1{\setboxz@h{$\m@th#1$}%
 \ifdim\wdz@>\tw@ em\mathaccent"0\msbfam@5D{#1}%
 \else\mathaccent"0365{#1}\fi}

\def\RIfM@{\relax\ifmmode}
\def\nonmatherr@#1{\errmessage{\string#1\space allowed only in math mode}}
\def\Bbb{\RIfM@\expandafter\Bbb@\else
 \expandafter\nonmatherr@\expandafter\Bbb\fi}
\def\Bbb@#1{{\Bbb@@{#1}}}
\def\Bbb@@#1{\fam\msbfam\relax#1}
\def\setboxz@h{\setbox\z@\hbox}
\def\wdz@{\wd\z@}
\catcode`\@=\csname pre amssym.def at\endcsname

\begin{document}

\title{\bf{Exponential stability of} fast driven systems, with an application to celestial mechanics\thanks{The authors are indebted to the anonymous Reviewers for their helpful suggestions. This research is supported by the ERC project 677793 Stable and Chaotic Motions in the Planetary Problem. Figures~\ref{figure1} and~\ref{figure2}  have been produced with {\sc mathematica}. {\bf MSC2000 numbers:}
primary:
34C20, 70F10,  37J10, 37J15, 37J40;
secondary: 
34D10,  70F07, 70F15, 37J25, 37J35. {\bf Keywords:} normal form theory; three--body problem; renormalizable integrability.
}}

\author{  
Qinbo Chen\footnote{Department of Mathematics, University of Padua, via Trieste 63, 3121 Padua, Italy. {\tt qinbochen1990@gmail.com}}\ \ and
Gabriella Pinzari\footnote{(corresponding author) Department of Mathematics, University of Padua, via Trieste 63, 3121 Padua, Italy. {\tt gabriella.pinzari@unipd.it}}}

\date{{February 2, 2021}}
\maketitle

\begin{abstract}\footnotesize{
We construct a normal form suited to {\it fast driven systems}. We call so systems including  actions ${\rm I}$,  angles {$\psi$}, and one  fast coordinate $y$, moving under the action of a vector--field $N$ depending only on ${\rm I}$ and $y$ and with vanishing ${\rm I}$--components. {In absence of the  coordinate $y$, such systems have been extensively investigated and it is known that, after a small perturbing term is switched on,  the normalised actions ${\rm I}$ turn to have exponentially small variations  compared to the size of the perturbation. We obtain the same result of the classical situation, with the additional benefit that  } no trapping argument is needed, as no small denominator arises. {We use the result to prove that, in the three--body problem, the level sets of a certain function called {\it Euler integral} have exponentially small variations in a short time,  closely to collisions.}}
\end{abstract}

\maketitle
{\footnotesize\tableofcontents}

\renewcommand{\theequation}{\arabic{equation}}
\setcounter{equation}{0}

\section{Description of the results}
We consider a $(n+1+m)$--dimensional vector--field $N$ which, expressed in local coordinates $({\rm I}, y, {{\psi}})\in{\mathbb P}={\mathbb I}\times {\mathbb Y}\times {\mathbb T}^m$ (where ${\mathbb I}\subset{\mathbb R}^n$, ${\mathbb Y}\subset{\mathbb R}$ are open and connected; ${\mathbb T}={\mathbb R}/(2\pi{\mathbb Z})$ is the standard torus), has the form
\begin{align}\label{X0}N({\rm I}, y) =v({\rm I}, y)\partial_y+\omega({\rm I}, y) \partial_{{\psi}}\, . \end{align}
The motion equations of $N$
\begin{align}\label{unperturbed vectorfield}\left\{\begin{array}{l}
\displaystyle \dot {\rm I}=0 \\\\
\displaystyle \dot y= v({\rm I}, y)\\\\
\displaystyle \dot{{\psi}}=\omega({\rm I}, y)
\end{array}\right.\end{align}
can be integrated in cascade:
\begin{align}\label{motions1}\left\{\begin{array}{l}
\displaystyle {\rm I}(t)={\rm I}_0 \\\\
\displaystyle  y(t)=\eta({\rm I}_0, t)\\\\
\displaystyle {{\psi}}(t)={{\psi}}_0+\int_{t_0}^t \omega({\rm I}_0, \eta({\rm I}_0, t'))d  t'
\end{array}\right.\end{align}
with $\eta({\rm I}_0, \cdot)$ being the general solution of the one--dimensional equation $\dot y(t)= v({\rm I}_0, y)$. 
This formula shows that along the solutions of $N$
the coordinates ${\rm I}$ (``actions'') remain constant, while the motion of the coordinates ${{\psi}}$  (``angles'') is coupled with the motion of the ``driving'' coordinate $y$. We assume that $v$ is suitably far from vanishing (for the problem considered in the paper $|v|$ has a positive lower bound). It is to be noted that, {without further assumptions on the function $v$ (like, for example, of being ``small'', or having a stationary point)}  {nothing prevents to the $y$ coordinate to move {\it fast}.} For this reason -- {with slight abuse due to the fact that fastness may nowise occur} -- we refer to  the solutions in~\eqref{motions1} as {\it fast driven system}. {The main risk of such kind of system is that}
the solution $q(t)=({\rm I}(t), y(t), {{\psi}}(t))$ of $N$ in~\eqref{motions1} 
 leaves the domain ${\mathbb P}$ {at a finite time}.
It is then convenient to define the {\it exit time from ${\mathbb P}$ under $N$}, or, more in general, the {\it exit time from a given $W\subseteq{\mathbb P}$ under $N$
 }, and denote it as $t^{N, W}_{\rm ex}$,  the {(possibly infinite)} first time  that $q(t)$
leaves  $W$.\\ Let us now replace the vector--field $N({\rm I}, y)$  with a new vector--field of the form
\begin{align}\label{perturbed}X({\rm I}, y, {{\psi}})=N({\rm I}, y)+P({\rm I}, y, {{\psi}}) \end{align}
where the ``perturbation''
\begin{align}P=P_{ 1}({\rm I}, y, {{\psi}})d{\rm I}+P_{ 2}({\rm I}, y, {{\psi}})dy+P_{ 3}({\rm I}, y, {{\psi}})d{{\psi}}\end{align} 
is, in some sense, ``small'' (see the next section for precise statements). Let $t^{X, W}_{\rm ex}$ be the exit time from $W$ under $X$, and let $\epsilon$ be   {a uniform upper bound} for {the absolute value of $P_{ 1}$}  on $W$. Then, one has a linear--in--time {\it a--priori} bound for the variations of ${\rm I}$, as follows
\begin{align}\label{aprioribound}|{\rm I}(t)-{\rm I}(0)|\le \epsilon t\quad \forall\ t:\ |t|< t^{X, W}_{\rm ex}\qquad W\subseteq{\mathbb P}\,.\end{align}
We are interested in improving the bound~\eqref{aprioribound}. To the readers who are familiar with  Kolmogorov--Arnold--Moser ({{\sc kam}}) or  Nekhorossev  theories,  this kind of problems is well known: {see}~\cite{arnold63, nehorosev77, poschel93, guzzoCB16}, or 
{\cite{cellettiC98, giorgilliLS09, locatelliG07, guzzoEP2020, volpiLS18} for applications to  realistic models.}
{Those are theories originally formulated for Hamiltonian vector--fields (next extended to more general ODEs), hence, in particular, with $n=m$ and the coordinate $y$ absent. In those cases the unperturbed motions of the  coordinates $({\rm I}, {{\psi}})$ are \begin{align}\label{unperturbedKAM}{\rm I}(t)={\rm I}_0\,,\quad {{\psi}}(t)={{\psi}}_0+\omega({\rm I}_0)t\end{align} and the properties of the motions after the perturbing term is switched on depend on the arithmetic properties of the frequency vector $\omega({\rm I}_0)$. Under suitable non--commensurability assumptions of $\omega({\rm I}_0)$ (referred to as ``Diophantine conditions''), {\sc kam} theory ensures the possibility of continuing the unperturbed motions \eqref{unperturbedKAM} for all times. Conversely, if $\omega({\rm I})$ satisfies, on an open set, an analytic property known as ``steepness'' (which is satisfied, e.g., if $\omega$ does not vanish and moreover if it is the gradient of a convex function), Nekhorossev theory allows to infer -- for {\it all} orbits -- a bound as in~\eqref{aprioribound}, with $e^{-C/\epsilon^{a}}$ replacing $\epsilon $  and $t_{\rm ex}^{X, W}=e^{C/\epsilon^{b}}$, with suitable $a$, $b$, $C>0$. It is to be remarked that in the Nekhorossev regime the exponential scale of $t_{\rm ex}^{X, W}$ is an intrinsic consequence of steepness, responsible of a process known as ``capture in resonance''.}
In  {the case considered in the paper such phenomenon does not seem to exist and hence the exit time} $t_{\rm ex}^{X, W}$ {has no reason to be long}. Nevertheless, motivated by an application to celestial mechanics  described below, we are interested with replacing $\epsilon$ in~\eqref{aprioribound} with a smaller number. 
We shall prove the following  result {(note that steepness conditions are not needed here).}  

\vskip.2in
\noindent
{\bf Theorem A} {\it Let $X=N+P$ be real--analytic, where $N$ is as in~\eqref{X0}, with $v\not\equiv 0$. Under suitable ``smallness'' assumptions involving $\omega$, $\partial\omega$, $\partial v$ and $P$, the bound in~\eqref{aprioribound} holds with $e^{-C/\epsilon^{a}}$ replacing $\epsilon$, with a suitable $a$, $C>0$.}

\vskip.1in
\noindent
A quantitative  statement of Theorem A  is given in {Theorem}~\ref{Normal Form Lemma} below. In addition, in view  of our application, we also discuss a version to the case when analyticity in ${{\psi}}$ fails; this is Theorem~\ref{Normal Form LemmaNEW}.

\noindent
To describe how we shall use {Theorem A (more precisely, Theorem~\ref{Normal Form LemmaNEW})}, we make a digression on the three--body problem and the {\it renormalizable integrability} of the simply averaged Newtonian potential~\cite{pinzari19}. The Hamiltonian governing the motions of  a three--body problem in the plane where the masses are $1$, $\mu$ and ${\kappa} appa$, is (see, e.g.,~\cite{fejoz04})
\begin{align} {\rm H}_{\rm 3b}&=\left(1+\frac{1}{{\kappa}}\right)\frac{\|{\mathbf y}\|^2}{2}+\left(1+\frac{1}{\mu}\right)\frac{\|{\mathbf y}'\|^2}{2}-\frac{{\kappa}}{\|{\mathbf x}\|}-\frac{\mu}{\|{\mathbf x}'\|}-\frac{{\kappa}\mu}{\|{\mathbf x}-{\mathbf x}'\|}+{{\mathbf y}\cdot{\mathbf y}'}\end{align}
where ${\mathbf y}$, ${\mathbf y}'\in {\mathbb R}^2$; ${\mathbf x}$, ${\mathbf x}'\in {\mathbb R}^2$, with ${\mathbf x}\ne 0\ne {\mathbf x}'$ and ${\mathbf x}\ne {\mathbf x}'$, are impulse--position coordinates; $\|\cdot\|$ denotes the Euclidean norm and the gravity constant  has been chosen equal to $1$, by a proper choice of the units system.
We rescale
\begin{align} 
({\mathbf y}', {\mathbf y})\to \frac{{{\kappa}}^2}{1+{{\kappa}}}({\mathbf y}', {\mathbf y})\, ,\quad ({\mathbf x}', {\mathbf x})\to \frac{1+{{\kappa}}}{{{\kappa}}^2}({\mathbf x}', {\mathbf x})\end{align}
 multiply the Hamiltonian by $\frac{1+{{\kappa}}}{{{\kappa}}^3}$ and obtain
\begin{align}\label{HH}{\rm H}_{\rm 3b}({\mathbf y}', {\mathbf y}, {\mathbf x}', {\mathbf x})=\frac{\|{\mathbf y}\|^2}{2}-\frac{1}{\|{\mathbf x}\|}
 +\delta\left(\frac{\|{\mathbf y}'\|^2}{2}-\frac{\alpha}{\|{\mathbf x}-{\mathbf x}'\|}-\frac{\beta}{\|{\mathbf x}'\|}\right)
 +\gamma {\mathbf y}\cdot{\mathbf y}'\end{align}
 with \begin{align}
\alpha:=\frac{\mu^2(1+{\kappa})}{{\kappa}(1+\mu)}\, ,\quad
\beta:=\frac{\mu^2(1+{\kappa})}{{\kappa}^2(1+\mu)}\, ,\quad \gamma :=\frac{{\kappa}}{1+{\kappa}}\, ,\quad \delta:=\frac{{\kappa}(1+\mu)}{\mu(1+{\kappa})}\,.
 \end{align}
 In order to simplify the analysis a little bit, we introduce a main assumption.
The Hamiltonian
 ${\rm H}_{\rm 3b}$  in~\eqref{HH} includes the Keplerian term  
\begin{align}\label{Kep}\frac{\|\mathbf y\|^2}{2}-\frac{1}{\|\mathbf x\|}=-\frac{1}{2\Lambda^2}\, . \end{align}
 We  assume that  this term is  ``leading'' in the  Hamiltonian. By averaging theory, this assumption allows us to replace
(at the cost of a small error) ${\rm H}_{\rm 3b}$ by its   $\ell$--average
\begin{align}\label{ovlH}\overline{\rm H}=-\frac{1}{2\Lambda^2}+\delta{\rm H}\ \end{align}
 where $\ell$ is the mean anomaly associated to~\eqref{Kep}, and\footnote{ Remark  that 
  ${\mathbf y}(\ell)$ has  vanishing $\ell$--average so that the last term in~\eqref{HH} does not survive.}
 \begin{align}\label{secular}
{\rm H}:=\frac{\|{\mathbf y}'\|^2}{2 }-\alpha{\rm U}-\frac{\beta}{\|{\mathbf x}'\|}
\end{align}
 with
\begin{align}\label{Usb}{\rm U}:=\frac{1}{2\pi}\int_0^{2\pi}\frac{d\ell}{\|{\mathbf x}'-{\mathbf x}(\ell)\|}\end{align}
being the  ``simply\footnote{Here, ``simply'' is used as opposed to the more familiar ``doubly'' averaged Newtonian potential, most often encountered in the literature; e.g.~\cite{laskarR95, fejoz04, pinzari-th09, chierchiaPi11b, chierchiaPi11c}.} averaged Newtonian potential''. {We recall that the mean anomaly $\ell$ is defined as the area spanned by ${\mathbf x}$ on the Keplerian ellipse generated by \eqref{Kep} relatively to the perihelion ${\mathbf P}$ of the ellipse, in $2\pi$ units.}
From now on we focus on the motions of the averaged Hamiltonian~\eqref{secular}, bypassing any quantitative statement concerning the averaging procedure, {as this would lead much beyond the purposes of the paper\footnote{{As we consider a region in phase space close where ${\mathbf x}'$ is very close to the instantaneous Keplerian orbit of ${\mathbf x}$, quantifying the values of the mass parameters and the distance which allow for the averaging procedure is a delicate (even though crucial) question, which, by its nature, demands careful use  of regularisations. Due to the non--trivial underlying analysis, we choose to limit ourselves to 
point out that
 the renormalizable integrability of the Newtonian potential has a nontrivial dynamical impact  on the simply averaged three--body problem, which explain the existence of the motions herewith discussed, which would not be justified otherwise. }}}. 
Neglecting the first term in~\eqref{ovlH}, which is an inessential additive constant for $\overline{\rm H}$ and reabsorbing the constant $\delta$ with a  time change, we are led  to look at the Hamiltonian ${\rm H}$ in~\eqref{secular}.  We denote as $\mathbb E$ the Keplerian ellipse generated by  Hamiltonian~\eqref{Kep}, for  negative values of the energy. Without loss
 of generality,  
assume $\mathbb E$ is not a circle and\footnote{We can do this as  the Hamiltonian ${\rm H}_{\rm 3b}$ rescale by a factor $\beta^{-2}$ as $({\mathbf y}', {\mathbf y})\to\beta^{-1}({\mathbf y}', {\mathbf y})$ and $({\mathbf x}', {\mathbf x})\to\beta^2({\mathbf x}', {\mathbf x})$.}
$\Lambda=1$.
 Remark that, as  the mean anomaly  $\ell$ is averaged out, we loose any information concerning the position of ${\mathbf x}$ on $\mathbb E$, so we shall only need two couples of coordinates for determining the shape of 
 $\mathbb E$ and the vectors ${\mathbf y}'$, ${\mathbf x}'$. These are:
 
\begin{itemize}
\item[\tiny\textbullet] the  ``Delaunay couple'' $({\rm G}, {\rm g})$, where ${\rm G}$ is the Euclidean length of ${\mathbf x}\times {\mathbf y}$ and ${\rm g}$ detects the perihelion. We remark that ${\rm g}$ is measured with respect to ${\mathbf x}'$ (instead of with respect to a fixed direction), as  the  SO(2) reduction we use 
 a rotating frame which moves with ${\mathbf x}'$ (compare the formulae in~\eqref{coord} below);

\item[\tiny\textbullet] the ``radial--polar couple''$({\rm R}, {\rm r})$, where ${\rm r}:=\|{\mathbf x}'\|$ and ${\rm R}:=\frac{{\mathbf y}'\cdot{\mathbf x}'}{\|{\mathbf x}'\|}$.
\end{itemize}

\noindent
Using the coordinates above, the Hamiltonian in~\eqref{secular} becomes
\begin{align}\label{3bpav}{\rm H}({\rm R}, {\rm G}, {\rm r}, {\rm g})=\frac{{\rm R}^2}{2}+\frac{({\rm C}-{\rm G})^2}{2{\rm r}^2}-\alpha{\rm U}({\rm r}, {\rm G}, {\rm g})
 -\frac{\beta}{{\rm r}}\end{align}
 where ${\rm C}=\|{\mathbf x}\times{\mathbf y}+{\mathbf x}'\times{\mathbf y}'\|$ is the total angular momentum of the system, and we have assumed  ${\mathbf x}\times{\mathbf y}\parallel{\mathbf x}'\times{\mathbf y}'$, so that $\|{\mathbf x}'\times{\mathbf y}'\|={\rm C}-\|{\mathbf x}\times{\mathbf y}\|={\rm C}-{\rm G}$.
 
\begin{figure}[htbp]
\subfigure[\label{a}]{\includegraphics[width=0.30\textwidth]{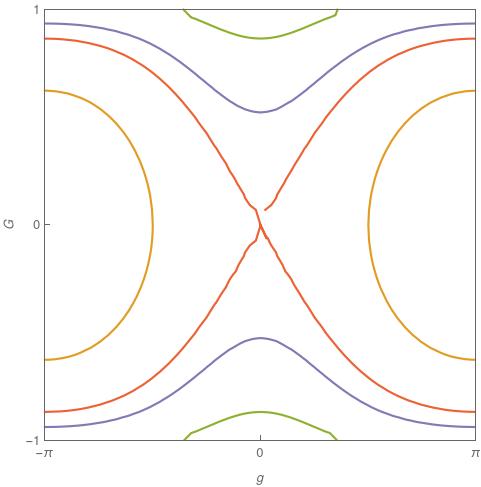}}
\hfill
\subfigure[\label{b}]{\includegraphics[width=0.30\textwidth]{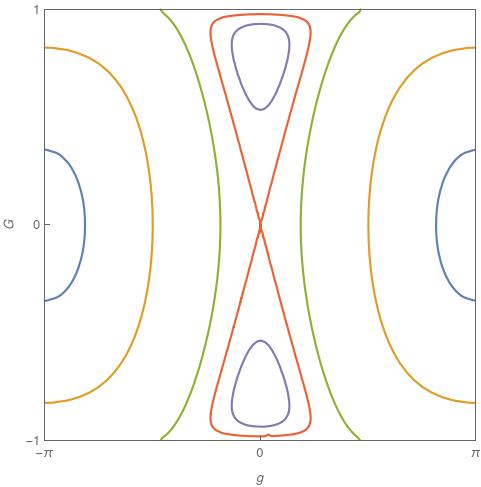}}
\hfill
\subfigure[\label{c}]{\includegraphics[width=0.30\textwidth]{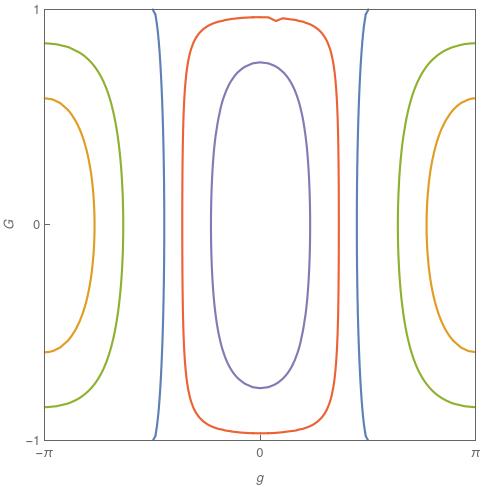}}
\caption{Sections, at ${\rm r}$ fixed,  of the level surfaces of ${\rm E}$.  (a): $0<{\rm r}<1$; (b): $1<{\rm r}< 2$; (c): ${\rm r}>2$.}  
\label{figure1}
\subfigure[\label{d}]{\includegraphics[width=0.30\textwidth]{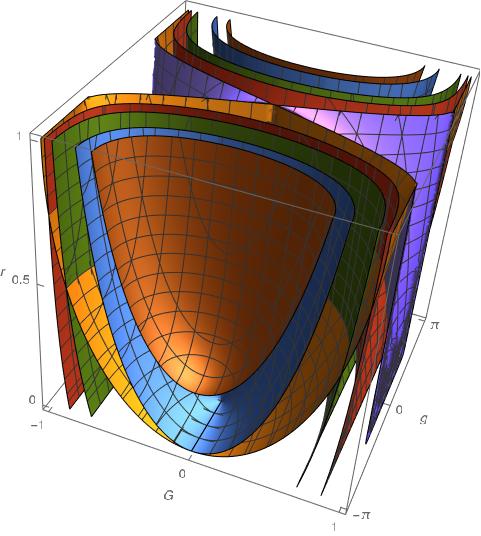}}\hfill
\subfigure[\label{e}]{\includegraphics[width=0.30\textwidth]{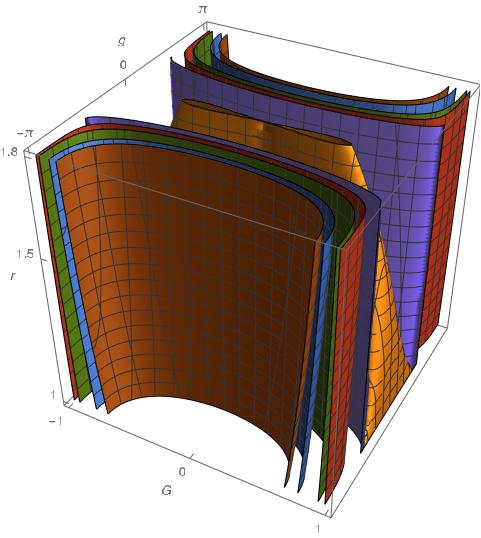}}\hfill
\subfigure[\label{f}]{\includegraphics[width=0.30\textwidth]{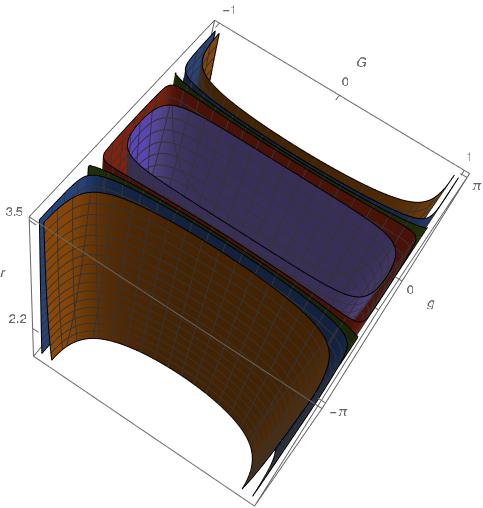}}
\caption{Logs of the level surfaces of ${\rm E}$ in the space $({\rm g}, {\rm G}, {\rm r})$. (a): $0<{\rm r}<1$; (b): $1<{\rm r}< 2$; (c): ${\rm r}>2$.}  
\label{figure2}
\end{figure}

\noindent
The Hamiltonian
~\eqref{3bpav} is now wearing 2 degrees--of--freedom. As the energy is conserved, its motions evolve on the $3$--dimensional manifolds ${\cal M}_c=\{{\rm H}=c\}$. On each of such manifolds  the evolution  is associated to a $3$--dimensional vector--field $X_c$, given by the velocity field of some triple of coordinates on
 ${\cal M}_c$. As an example, one can take  the triple $({\rm r}, {\rm G}, {\rm g})$, even though  a more convenient choice will be done below.   
To describe the motions we are looking for, we need to  recall a remarkable property of the function ${\rm U}$, pointed out in~\cite{pinzari19}. First of all, one has to note that ${\rm U}$  is integrable, as it is a function of $({\rm r}, {\rm G}, {\rm g})$ only. But the main point is that there exists a function ${\rm F}$ of two arguments such that
\begin{align}\label{relation***}{\rm U}({\rm r}, {\rm G}, {\rm g})={\rm F}({\rm E}({\rm r}, {\rm G}, {\rm g}), {\rm r})\end{align}
where
\begin{align}\label{E}{\rm E}({\rm r}, {\rm G}, {\rm g})={\rm G}^2+{\rm r}\sqrt{1-{\rm G}^2}\cos{\rm g}\,.\end{align}
The function ${\rm E}$ is referred to as {\it Euler integral}, and we express~\eqref{relation***} by saying  that ${\rm U}$ is {\it renormalizable integrability via the Euler integral}.  Such cirumstance implies that  the level sets of ${\rm E}$, namely the curves \begin{align}\label{level curves}{\rm G}^2+{\rm r}\sqrt{1-{\rm G}^2}\cos{\rm g}={\cal E}\end{align} are also level sets of ${\rm U}$. On the other hand, the phase portrait of~\eqref{level curves} keeping ${\rm r}$ fixed is completely explicit and has been  studied in~\cite{pinzari20b}. We recall it now. Let us fix (by periodicity of ${\rm g}$) the strip $[-\pi, \pi]\times [-1, 1]$.
For $0<{\rm r}<1$ or $1<{\rm r}<2$ it includes two minima $(\pm\pi, 0)$  on the ${\rm g}$--axis; two symmetric maxima on the ${\rm G}$--axis and one saddle point at $(0, 0)$. 
When ${\rm r}>2$ the saddle point disappears and $(0, 0)$ turns to be a maximum. 
The phase portrait includes two separatrices  when $0<{\rm r}<1$ or $1<{\rm r}<2$; one separatrix if ${\rm r}>2$. 
These are the level sets $$\left\{\begin{array}{l}{\cal S}_0({\rm r})=\{{\cal E}={\rm r}\}\,,\quad 0<{\rm r}<1\,,\ 1<{\rm r}<2\\\\
{\cal S}_1({\rm r})=\{{\cal E}=1\}\,,\quad 0<{\rm r}<1\,,\ 1<{\rm r}<2\,,\ {\rm r}>2
\end{array}\right.$$ 
with ${\cal S}_0({\rm r})$ being the separatrix through the saddle; ${\cal S}_1({\rm r})$ the level set through circular orbits.
 Rotational motions in between ${\cal S}_0({\rm r})$ and ${\cal S}_1({\rm r})$, do exist only for $0<{\rm r}<1$. The minima and the maxima are surrounded by librational motions and different motions (librations about different equilibria or rotations) are separated by ${\cal S}_0({\rm r})$ and ${\cal S}_1({\rm r})$. All of this is represented in  Figure~\ref{figure1}.
In Figure~\ref{figure2} the same level sets are drawn in the 3--dimensional space $({\rm r}, {\rm G}, {\rm g})$. The spatial visualisation turns out to be useful for the purposes of the paper, as the coordinate ${\rm r}$, which stays fixed under ${\rm E}$, is instead  moving under ${\rm H}$, due to its dependence on ${\rm R}$; see~\eqref{3bpav}. We denote as ${\cal S}_0$ the union of all the ${\cal S}_0({\rm r})$ with $0\le {\rm r}\le 2$.
It is to be noted  that, while ${\rm E}$ is perfectly defined along ${\cal S}_0$,  ${\rm U}$ is not so. Indeed, as
\begin{align}\label{S0}{\cal S}_0({\rm r})=\Big\{({\rm G}, {\rm g}):\quad {\rm G}^2+{\rm r}\sqrt{1-{\rm G}^2}\cos{\rm g}={\rm r}\, ,\ -1\le {\rm G}\le 1\, ,\quad {\rm g}\in {\mathbb T}\Big\}\qquad 0\le{\rm r}<2\end{align}
we have\footnote{Rewriting~\eqref{S0} as
\begin{align}\label{relation1}
{\rm r}=\frac{{\rm G}^2}{1-\sqrt{1-{\rm G}^2}\cos{\rm g}}
\end{align}
 tells us that $({\rm G}, {\rm g})\in {\cal S}_0({\rm r})$ if and only if
${\mathbf x}'$  occupies in the ellipse ${\mathbb E}$ the position with true anomaly $\nu=\pi-{\rm g}$.} ${\rm U}({\rm r}, {\rm G}, {\rm g})=\infty$  for $({\rm G}, {\rm g})\in{\cal S}_0({\rm r})$, for all $0\le{\rm r}\le2$.

\noindent
The natural question now raises  whether 
any of the ${\cal E}$--levels in Figure~\ref{figure2} is an ``approximate'' invariant manifold for  the Hamiltonian ${\rm H}$ in~\eqref{3bpav}.  In~\cite{pinzari20a} and~\cite{diruzzaDP20} a positive answer has been given for case ${\rm r}>2$, corresponding to panels (c). 
In this paper, we  want to focus on motions close to 
${\cal S}_0$ with ${\rm r}$ in a left neighbourhood of $2$ (panels (b)). Such portion  of phase space is denoted as ${\cal C}$. By the discussion above, motions in ${\cal C}$ are to be understood as  ``quasi--collisional''.

 \noindent
To state our result, we denote as ${\rm r}_{\rm s}(A)$ the value of ${\rm r}$ such that the area encircled by ${\cal S}_0({\rm r}_{\rm s}(A))$ is $A$. Then the set $\{\exists\ A:\ {\rm r}={\rm r}_{\rm s}(A)\}$ corresponds to ${\cal S}_0$. We prove:

\vskip.1in
\noindent
{\bf Theorem B} {\it Inside the region ${\cal C}$  there exists an open set $W$ such that along any motion with initial datum in $W$, for all $t$ with $|t|\le t_{\rm ex}^{X, W}$, the ratio between the absolute variations of the  Euler integral ${\rm E}$ from time $0$ to time $t$, for all $|t|\le t_{\rm ex}^{X, W}$, and the {\it a--priori} bound $\epsilon t$ (where $\epsilon	:=|P_1|_\infty${, with $P_1$ being the action component of the vector--field})
 does not exceed $C e^{-L^3/C}$, provided that the initial value of ${\rm r}$ is  $e^{-L}$ away from ${\rm r}_{\rm s}(A)$, with $L>0$ sufficiently large.
}
\vskip.1in
\noindent
The proof of Theorem B, fully given in the next section, relies
on a careful choice of coordinates $(A, y, \psi)$ on ${\cal M}_c$, where $y$ is diffeomorphic to ${\rm r}$, while $(A, \psi)$ are the action--angle coordinates of ${\rm E}({\rm r}, \cdot, \cdot)$, such that the associated vector--field has the form in~\eqref{perturbed} with $n=m=1$. The diffeomorphism ${\rm r}\to y$ allows  $X_c$  to keep its regularity  upon ${\cal S}_0$.

\vskip.1in
\noindent
Before switching to proofs, we recall how the theme of collisions in  $ N$--body problems (with ${N}\ge 3$) has been treated so far. 
As the literature in the field  in countless,  by no means we claim completeness.
In the late 1890s H. Poincar\'e~\cite{Poincare:1892} conjectured  the existence of special solutions in a model of the three--body problem usually referred to as planar, circular, restricted three--body problem ({\sc pcrtbp}).
According to Poincar\'e's conjecture, when
one of the primaries has a small mass $\mu$, the orbit of an infinitesimal body approaching a close encounter with the small primary consists of 
two Keplerian arcs glueing so as to form a cusp.  These solutions were named by him {\it second species solutions}, and their  existence has been next proved in~\cite{bolotin2005, bolotin2006a, bolotin2006b, bolotin2006c, bolotinN2013, marcoN1995, henrard1980}.
 In the early 1900s, J.~Chazy classified all the possible  final motions of the three--body problem, including the possibility of collisions
~\cite{chazy1922}. The study was reconsidered in~\cite{alekseev1971, alekseev1981}. 
 After the advent of {{\sc kam}} theory, the existence of almost--collisional quasi--periodic orbits was proven~\cite{chenciner1988, fejoz02, lei15}.
 The papers~\cite{saari1971, saari1973, fleischerK2019a, fleischerK2019b, moeckel1980, moeckel1981, moeckel1989, moeckel2007}
 deal with rare occurrence of collisions or the existence of chaos in the proximity of collisions.
In
\cite{guardiaKZ19} it is proved that for {\sc pcrtbp} there exists an open set in phase space of fixed measure, where the set of initial points which lead to collision is 
 ${\rm O}(\mu^\alpha)$ dense with some $0<\alpha<1$. In \cite{levicivita06} it is proved that, after collision regularisation, {\sc pcrtbp} is integrable in a neighbourhood of collisions.
In~\cite{cardinG18a, cardinG18} the  result has been recently extended to the spatial version, often denoted {\sc scrtbp}. 
 
\section{A Normal Form Theorem for fast driven systems}
In the next Sections~\ref{Weighted norms}--\ref{Proof of NFL}  we state and prove a Normal Form Theorem ({\sc nft}) for  real--analytic systems. For the purpose of the paper, in Section~\ref{A generalisation when the dependence} we generalise the result,   allowing the dependence on the angular coordinate $\psi$ to be just $C^{\ell_*}$ ($\ell_*\in{\mathbb N}$), rather than  holomorphic. In all cases,  we limit to the case $n=m=1$. Generalisations to $n$, $m\ge 1$ are straightforward.

\subsection{Weighted norms}\label{Weighted norms}
Let us consider a 3--dimensional vector--field
 \begin{align}\label{vectorfield}  ({\rm I}, y, \psi)\in{\mathbb P}_{r, \sigma, s}:={\mathbb I}_r \times{\mathbb Y}_\sigma \times {{\mathbb T}}_s\to X=(X_1, X_2, X_3)\in{\mathbb C}^3\end{align} 
 where ${\mathbb I}\subset {\mathbb R}$, ${\mathbb Y}\subset {\mathbb R}$ are open and connected; ${\mathbb T}= {\mathbb R}/(2\pi {\mathbb Z})$,
 which has the form~\eqref{perturbed}.
 As usual, if $A\subset {\mathbb R}$ and $r${,$s>0$},  the symbols $A_r${, ${\mathbb T}_s$} denote the complex $r$, {$s$}--neighbourhoods of $A${,${\mathbb T}$:}
\begin{align}A_r:=\bigcup_{x\in A}B_r(x)\,,\qquad {{\mathbb T}_s:=\big\{\psi=\psi_1+{\rm i}\psi_2:\ \psi_1\in {\mathbb T}\,,\ \psi_2\in {\mathbb R}\,,\ |\psi_2|<s\big\}\,,}\end{align} with $B_r(x)$ being the complex ball centred at $x$ with radius $r$.
We assume each $X_i$ to be holomorphic in
${\mathbb P}_{r, \sigma, s}$, meaning the it has a finite weighted norm defined below. If this holds, we simply write $X\in {\cal O}^3_{r, \sigma, s}$.

\noindent
For functions $f:\ ({\rm I}, y, \psi)\in{\mathbb I}_r\times {\mathbb Y}_\sigma\times {\mathbb T}_s\to {\mathbb C}$, we write $f\in{\cal O}_{r, \sigma, s}$ if $f$ is holomorphic in ${\mathbb P}_{r, \sigma, s}$. We let
\begin{align}\label{normfOLD}\|f\|_{u}:=\sum_{k\in {\mathbb Z}}\,\sup_{{\mathbb I}_r\times{\mathbb Y}_\sigma}
|f_{\kappa}({\rm I}, y)|
\,e^{|k|s}\qquad u=(r, \sigma, s)
\end{align}
where
\begin{align}f=\sum_{k\in {\mathbb Z}} f_k({\rm I}, y)e^{{\rm i} k\psi}\end{align}
is the Fourier series associated to $f$ relatively to the $\psi$--coordinate.  For $\psi$--independent 
  functions or vector--fields we simply write $\|\,\cdot\,\|_{r, \sigma}$.
\\
For vector--fields $X:\ ({\rm I}, y, \psi)\in{\mathbb I}_r\times {\mathbb Y}_\sigma\times {\mathbb T}_s\to X=(X_1, X_2, X_3)\in{\mathbb C}^3$, we write $X\in{\cal O}^3_{r, \sigma, s}$ if $X_i\in{\cal O}_{r, \sigma, s}$ for $i=1$, $2$, $3$. We 
define the {\it weighted norms} 
\begin{align}\label{normXOLD}\VERT X \VERT_{u}^{w}:=\sum_i w^{-1}_i\|X_i\|_{u}\end{align}
where
$w=(w_1$, $w_2$, $w_3)\in {\mathbb R}_+^3$ are the {\it weights}.
The wighted norm affords the following properties.

\begin{itemize}
\item[{\tiny\textbullet}] {Monotonicity:} 
\begin{align}\label{monotonicity}
\VERT X \VERT_{u}^{w}\le \VERT X \VERT_{u'}^{w}\, ,\quad \VERT X \VERT_{u}^{w'}\le \VERT X \VERT_{u}^{w}\quad \forall\ u\le u'\, ,\ w\le w'
\end{align}
where $u\le u'$ means $u_i\le u_i'$ for $i=1$, $2$, $3$.

\item[{\tiny\textbullet}] {Homogeneity:} 
\begin{align}\label{homogeneity}\VERT X \VERT_{u}^{\alpha w}=\alpha^{-1}\VERT X \VERT_{u}^{w}\qquad \forall\ \alpha>0\, . \end{align}
\end{itemize}
  
\subsection{The Normal Form {Theorem}}
We now state the main result of this section.   Observe that the nature of the system does not give rise to any non--resonance condition or ultraviolet cut--off. 
We name Normal Form {Theorem} the following

\begin{theorem}[{\sc nft}]\label{Normal Form Lemma}
Let $u=(r, \sigma, s)$; $X=N+P\in {\cal O}^3_{u}$ and let $w=(\rho$, $\tau$, $t)\in {\mathbb R}_+^3$.
Put
\begin{align}\label{delta1}
 Q:=3\,{\rm diam}({\mathbb Y}_\sigma)\left\|\frac{1}{v}\right\|_{r, \sigma}\end{align}
 and\footnote{${\rm diam}({\cal A})$ denotes diameter of the set ${\cal A}$.} assume that for some ${p}\in {\mathbb N}$,  $s_2\in {\mathbb R}_+$, the following inequalities are satisfied:
 \begin{align}\label{NEWu+positive}
0<\rho<\frac{r}{8}\, ,\quad 0<\tau< e^{-s_2}\frac{\sigma}{8}\, ,\quad 0<t<\frac{s}{10}
\end{align}
and
 \begin{align}\label{NEWnewcond2}\chi&:= \frac{{\rm diam}({\mathbb Y}_\sigma)}{s_2}\left\|\frac{\partial_y v}{v}\right\|_{r, \sigma}
\le 1	\\\label{theta1}\theta_1&:= 2\,e^{s_2}{\rm diam}({\mathbb Y}_\sigma)\left\|\frac{\partial_y\omega}{v}\right\|_{r, \sigma}\frac{\tau}{t}\le 1\\
\label{theta2}\theta_2&:= 4\,{\rm diam}({\mathbb Y}_\sigma)\left\|\frac{\partial_{\rm I} v}{v}\right\|_{r, \sigma}\frac{\rho}{\tau}\le 1\\
\label{theta3}\theta_3&:= 8\,{\rm diam}({\mathbb Y}_\sigma)
\left\|\frac{\partial_{\rm I}\omega}{v}\right\|_{r, \sigma}\frac{\rho}{t}\le 1\\
\label{NEWnewsmallness}\eta^2&:= \max\left\{\frac{{\rm diam}({\mathbb Y}_\sigma)}{t}\left\|\frac{\omega}{v}\right\|_{r, \sigma}\, ,\ 2^7\,e^{2 s_2}Q^2  (\VERT P\VERT_{u}^{w})^2\right\}<\frac{1}{{p}}\, . \end{align}
Then, with
\begin{align}u_*=(r_\star, \sigma_\star, s_\star)\, ,\quad r_\star:=r-8\rho\, , \quad\sigma_\star=\sigma-8 e^{s_2}\tau\, , \quad s_\star=s-10 t\end{align}
 there exists  a real--analytic change of coordinates $\Phi_\star$
such that $X_\star:=\Phi_\star X\in {\cal O}^3_{u_\star}$ and
 $X_\star=N+P_\star$,
with
\begin{align}\label{P*}\VERT P_\star\VERT^w_{u_\star}<2^{-({p}+1)}\VERT P\VERT^w_{u}\, . \end{align}
\end{theorem}

\begin{remark}[Proof of Theorem A]\rm
{Theorem}~\ref{Normal Form Lemma} immediately implies Theorem~A, with $C=\min\{2^{-7}Q^{-2}e^{-2 s_2} \varrho^2 \log 2\,,\ t/{\rm diam}{\mathbb Y_\sigma}\}$, $a=2$, provided that $\varrho:=\frac{\epsilon^2}{(\VERT P\VERT_u^w)^2
}$ is of ``order one'' with respect to $\epsilon$. The mentioned ``smallness assumptions'' correspond to
conditions~\eqref{NEWu+positive}--\eqref{theta3} and $\left\|\frac{\omega}{v}\right\|_{r, \sigma}\ll (\VERT P\VERT_u^w)^2$.
\end{remark}

\subsection{The Step Lemma}

\noindent
We denote as  
\begin{align}\label{Lie}e^{{\cal L}_Y}=\sum_{{k}\ge 0}\frac{{\cal L}^{{k}}_Y}{{k}!}\end{align}
 the formal Lie series associated to $Y$, where \begin{align}[Y, X]=J_X Y-J_Y X\, ,\quad (J_Z)_{ij}:=\partial_j Z_i\end{align} denotes Lie brackets of two vector--fields, with
 \begin{align}{\cal L}_Y:=[Y, \cdot]\end{align}
being the Lie operator.

\begin{lemma}\label{iteration lemma}
Let $X=N
+P\in {\cal O}^3_{u}$, with $u=(r, \sigma, s)$, $N$ as in~\eqref{NZ}, $s_1$, $s_2>0$. Assume \begin{align}\label{existence}
\frac{{\rm diam}({\mathbb Y}_\sigma)}{s_1}\left\|\frac{\omega}{v}\right\|_{r, \sigma}\le 1\, ,\quad \frac{{\rm diam}({\mathbb Y}_\sigma)}{s_2}\left\|\frac{\partial_y v}{v}\right\|_{r, \sigma}
\le 1\end{align} 
and that 
$P$  is so small that

\begin{align}\label{smallcond}
Q \VERT P\VERT^{w} _{u}<1\qquad Q:=3{\rm diam}({\mathbb Y}_\sigma)\left\|\frac{1}{v}\right\|_{r, \sigma}\, ,\quad w=(\rho, \tau, t)
\end{align}
Let $\rho_*$, $\tau_*$, $t_*$ be defined via 
\begin{align}\label{bounds1}
\frac{1}{\rho_*}&=\frac{1}{\rho}-{\rm diam}({\mathbb Y}_\sigma)\left\|\frac{\partial_{\rm I} v}{v}\right\|_{r, \sigma}\left(\frac{1}{\tau}-e^{s_2}{\rm diam}({\mathbb Y}_\sigma)
\left\|\frac{\partial_y\omega}{v}\right\|_{r, \sigma}\frac{1}{t}
\right)\nonumber\\
&- {\rm diam}({\mathbb Y}_\sigma)
\left(\left\|\frac{\partial_{\rm I}\omega}{v}\right\|_{r, \sigma}+e^{s_2}{\rm diam}({\mathbb Y}_\sigma)\left\|\frac{\partial_{\rm I} v}{v}\right\|_{r, \sigma} \left\|\frac{\partial_y\omega}{v}\right\|_{r, \sigma}\right)\frac{1}{t}\nonumber\\
\frac{1}{\tau_*}&=\frac{e^{-s_2}}{\tau}-{\rm diam}({\mathbb Y}_\sigma)\left\|\frac{\partial_y\omega}{v}\right\|_{r, \sigma}\frac{1}{t}\nonumber\\
t_*&=t
\end{align}
and assume
\begin{align}\label{w*positive}w_*=(\rho_*, \tau_*, t_*)\in {\mathbb R}_+^3\, ,\qquad u_*=(r-2\rho_*, \sigma-2\tau_*, s-3s_1-2t_*)\in {\mathbb R}_+^3\, . \end{align} Then there exists $Y\in {\cal O}^3_{u_*+w_*}$  such that
$X_+:=e^{{\cal L}_Y}X\in {\cal O}^3_{u_*}$ and 
$X_+=N+
P_+$, with
\begin{align}\label{P+}
 \VERT P_+\VERT^{w_*}_{u_*}
\le \frac{2Q
\left(\VERT P\VERT^w_u\right)^2
}{1-Q
\left\VERT P\right\VERT^w_u}
\end{align}
\end{lemma}

\noindent
In the next section, we shall use {Lemma}~\ref{iteration lemma} in the following ``simplified'' form. 
\begin{lemma}[Step Lemma]\label{simplifiedsteplemma}
If~\eqref{existence},~\eqref{smallcond} and~\eqref{w*positive} are replaced with  
\begin{align}\label{newcond1}\begin{split}
&2\,e^{s_2}{\rm diam}({\mathbb Y}_\sigma)\left\|\frac{\partial_y\omega}{v}\right\|_{r, \sigma}\frac{\tau}{t}\le 1\\
&4\,{\rm diam}({\mathbb Y}_\sigma)\left\|\frac{\partial_{\rm I} v}{v}\right\|_{r, \sigma}\frac{\rho}{\tau}\le 1\\
&8\,{\rm diam}({\mathbb Y}_\sigma)
\left\|\frac{\partial_{\rm I}\omega}{v}\right\|_{r, \sigma}\frac{\rho}{t}\le 1
\end{split}
\end{align}
\begin{align}\label{newcond2}
\frac{{\rm diam}({\mathbb Y}_\sigma)}
{t}\left\|\frac{\omega}{v}\right\|_{r, \sigma}\le 1\, ,\quad \frac{{\rm diam}({\mathbb Y}_\sigma)}{s_2}\left\|\frac{\partial_y v}{v}\right\|_{r, \sigma}
\le 1
\end{align}
\begin{align}
\label{u+positive}
0<\rho<\frac{r}{4}\, ,\quad 0<\tau<\frac{\sigma}{4}e^{-s_2}\, ,\quad 0<t<\frac{s}{5}
\end{align}
\begin{align}
\label{newsmallness} 
2Q\VERT P\VERT ^w_u<1
\end{align}
then $X_+=N+P_+\in {\cal O}^3_{u_+}$ and \begin{align}\label{finalineq}\VERT P_+\VERT_{u_+}^{w}\le 8 e^{s_2}Q
(\VERT P \VERT_u^w)^2\, . \end{align}
with
\begin{align}u_+:=(r-4\rho, \sigma-4\tau e^{s_2}, s-5t)\, . \end{align}
\end{lemma}
{\bf Proof\ }The  inequality in~\eqref{newcond2} guarantees that one can take $s_1=t$, while the inequalities in~\eqref{newcond1} and~\eqref{u+positive} imply
\begin{align}\frac{1}{\rho_*}\ge \frac{1}{2\rho}\, ,\quad \frac{1}{\tau_*}\ge \frac{e^{-s_2}}{2\tau}\end{align}
whence, as $t_*=t$,
\begin{align}w_*<2 e^{s_2} w\, ,\qquad u_*\ge u_+>0\, . \end{align}
Then~\eqref{finalineq} is implied by~\eqref{P+}, monotonicity and homogeneity~\eqref{monotonicity}--\eqref{homogeneity}, and the inequality in~\eqref{newsmallness}. $\square$

\noindent
 {To prove {Lemma}~\ref{iteration lemma}, we look for a change of coordinates which conjugates the vector--field $X=N+P$ to a new vector--field $X_+=N_++P_+$, where $P_+$ depends in the coordinates ${\rm I}$ at higher orders. The procedure we follow is reminiscent of classical techniques of normal form theory, where one chooses the transformation so that $X_+=e^{{\cal L}_Y} X$, with the operator $e^{{\cal L}_Y}$ being defined as in \eqref{Lie}. As in the classical case, $Y$ will be chosen as the solution of a certain 
 ``homological equation'' which allows to eliminate the first order terms depending on $\psi$ of $P$. 
 However, as stated in Lemma \ref{iteration lemma}, differently from the classical situation, one can take $N=N_+$, which is another way of saying that
 it is  possible to choose $Y$ such in a way to solve
  \begin{align}\label{homeq111}
{\cal L}_N[Y]=P
\end{align}
regardless  $P$ has vanishing average or not -- or, in other words, that  {\it also the resonant terms} of the perturbing term will be killed. Note also that no ``ultraviolet cut--off'' is used.
   Equation \eqref{homeq111} is precisely what is discussed in  Lemma~\ref{estimates} and Proposition~\ref{homeq1} below.
}

\noindent
Fix $y_0\in {\mathbb Y}$; $v$, $\omega:{\mathbb I}\times {\mathbb Y}\to {\mathbb R}$, with $v\not\equiv 0$. We define, formally, the operators ${\cal F}_{v, \omega}$ and ${\cal G}_{v, \omega}$
as acting on  functions $g:{\mathbb I}\times {\mathbb Y}\times {\mathbb T}\to {\mathbb R}$ as
\begin{align}\label{FandG}
{\cal F}_{v, \omega}[g]({\rm I}, y, \psi)&:= \int_{y_0}^y\frac{g\left({\rm I}, \eta, \psi+\int_y^\eta \frac{\omega({\rm I}, \eta')}{v({\rm I}, \eta')}d\eta' \right)}{v({\rm I}, \eta)}d\eta\nonumber\\
{\cal G}_{v, \omega}[g]({\rm I}, y, \psi)&:= \int_{y_0}^y\frac{g\left({\rm I}, \eta, \psi+\int_y^\eta \frac{\omega({\rm I}, \eta')}{v({\rm I}, \eta')}d\eta' \right)e^{-\int_y^\eta \frac{\partial_y v({\rm I}, \eta')}{v({\rm I}, \eta')}d\eta'}}{v({\rm I}, \eta)}d\eta
\end{align}

\noindent
Observe that, when existing,
${\cal F}_{v, \omega}$, ${\cal G}_{v, \omega}$  send zero--average functions to zero--average functions.

\noindent
The existence ${\cal F}_{v, \omega}$, ${\cal G}_{v, \omega}$ is established by the following 

\begin{lemma}\label{estimates}
If inequalities~\eqref{existence} hold,
then
\begin{align}{\cal F}_{v, \omega}\, ,\ {\cal G}_{v, \omega}:\quad {\cal O}_{r, \sigma, s}\to {\cal O}_{r, \sigma, s-s_1}\end{align}
and
\begin{align}\|{\cal F}_{v, \omega}[g]\|_{r, \sigma, s-s_1}\le {\rm diam}({\mathbb Y}_\sigma)\left\|\frac{g}{v}\right\|_{r, \sigma, s}\, ,\quad \|{\cal G}_{v, \omega}[g]\|_{r, \sigma, s-s_1}\le e^{s_2}\,{\rm diam}({\mathbb Y}_\sigma) \left\|\frac{g}{v}\right\|_{r, \sigma, s}\end{align}
 \end{lemma}
 The proof of Lemma~\ref{estimates} is obvious from the definitions~\eqref{FandG}.

 \begin{proposition}\label{homeq1}
Let 
\begin{align}\label{NZ}N=(0, v({\rm I}, y), \omega({\rm I}, y))\, ,\qquad Z=(Z_1({\rm I}, y, \psi), Z_2({\rm I}, y, \psi), Z_3({\rm I}, y, \psi))\end{align}
belong to ${\cal O}^3_{r, \sigma, s}$ and assume~\eqref{existence}.
Then  the ``homological equation''
\begin{align}\label{homeq}{\cal L}_N[Y]=Z\end{align}
has a solution $Y\in {\cal O}_{r, \sigma, s-3s_1}$ verifying
\begin{align}\label{bounds}
\VERT Y\VERT_{r, \sigma, s-3 s_1}^{\rho_*, \tau_*, t_*}\le {\rm diam}({\mathbb Y}_\sigma)\left\|\frac{1}{v}\right\|_{r, \sigma}\VERT Z\VERT_{r, \sigma, s}^{\rho, \tau, t}
\end{align}
with $\rho_*$, $\tau_*$, $t_*$ as in~\eqref{bounds1}.
\end{proposition}
{\bf Proof\ }We expand $Y_j$ and $Z_j$ along the Fourier basis
\begin{align}Y_j({\rm I}, y, \psi)=\sum_{k\in {\mathbb Z}} Y_{j, k}({\rm I}, y)e^{{\rm i} k\psi}\, ,\quad Z_j({\rm I}, y, \psi)=\sum_{k\in {\mathbb Z}} Z_{j, k}({\rm I}, y)e^{{\rm i} k\psi}\, ,\quad  j=1,\ 2,\ 3\end{align}
Using \begin{align}{\cal L}_N[Y]=[N, Y]=J_Y N-J_N Y\end{align}
where
$(J_Z)_{ij}=\partial_j Z_i$
are the Jacobian matrices,
we rewrite~\eqref{homeq} as

\begin{align}\label{equations}
Z_{1, k}({\rm I}, y)&=v({\rm I}, y)\partial_y Y_{1, k} +{\rm i} k \omega({\rm I}, y) Y_{1, k} \nonumber\\
Z_{2, k}({\rm I}, y)&=v({\rm I}, y)\partial_y Y_{2, k} +({\rm i} k \omega({\rm I}, y)-\partial_y v({\rm I}, y)) Y_{2, k}-\partial_{\rm I} v({\rm I}, y) Y_{1, k}\nonumber\\
Z_{3, k}({\rm I}, y)&=v({\rm I}, y)\partial_y Y_{3, k} +{\rm i} k \omega({\rm I}, y) Y_{3, k}-\partial_{\rm I} \omega({\rm I}, y) Y_{1, k}-\partial_y \omega({\rm I}, y) Y_{2, k}\,.
\end{align}
Regarding~\eqref{equations} as equations for  $Y_{j, k}$, we find the solutions
\begin{align}
Y_{1, k}&=\int_{y_0}^y\frac{Z_{1, k}({\rm I}, \eta)}{v({\rm I}, \eta)}e^{{\rm i} k\int_y^\eta \frac{\omega({\rm I}, \eta')}{v({\rm I}, \eta')}d\eta'}d\eta\nonumber\\
Y_{2, k}&=\int_{y_0}^y\frac{Z_{2, k}({\rm I}, \eta)+\partial_{\rm I} v
Y_{1, k}
}{v({\rm I}, \eta)}e^{\int_y^\eta \frac{{\rm i} k \omega({\rm I}, \eta')-\partial_y v({\rm I}, \eta')}{v({\rm I}, \eta')}d\eta'}d\eta\nonumber\\
Y_{3, k}&=\int_{y_0}^y\frac{Z_{3, k}({\rm I}, \eta)+\partial_{\rm I} \omega({\rm I}, \eta) Y_{1, k}+\partial_y \omega({\rm I}, \eta) Y_{2, k}}{v({\rm I}, \eta)}e^{{\rm i} k\int_y^\eta \frac{\omega({\rm I}, \eta')}{v({\rm I}, \eta')}d\eta'}d\eta\nonumber\\
\end{align} 
multiplying by $e^{{\rm i} k\psi}$ and summing over $k\in {\mathbb Z}$ we find \begin{align}\label{Yi}
Y_1&={\cal F}_{v, \omega}[Z_1]\nonumber\\
 Y_2&={\cal G}_{v, \omega}[Z_2]+{\cal G}_{v, \omega}[\partial_{\rm I} v\,Y_1]\, ,\nonumber\\ 
 Y_3&={\cal F}_{v, \omega}[Z_3]+{\cal F}_{v, \omega}[\partial_{\rm I} \omega\,Y_1]+{\cal F}_{v, \omega}[\partial_y \omega\,Y_2]\, . \end{align} 
Then, by Lemma~\ref{estimates},
\begin{align}
&\|Y_1\|_{r, \sigma, s-s_1}\le  {\rm diam}({\mathbb Y}_\sigma)\left\|\frac{1}{v}\right\|_{r, \sigma}\left\|Z_1\right\|_{r, \sigma, s}\nonumber\\
&\|Y_2\|_{r, \sigma, s-2s_1}
\le e^{s_2}{\rm diam}({\mathbb Y}_\sigma)\left\|\frac{1}{v}\right\|_{r, \sigma}\left\|Z_2\right\|_{r, \sigma, s-s_1}+e^{s_2}{\rm diam}({\mathbb Y}_\sigma)^2\left\|\frac{1}{v}\right\|_{r, \sigma}\left\|\frac{\partial_{\rm I} v}{v}\right\|_{r, \sigma}\left\|Z_1\right\|_{r, \sigma, s}\nonumber\\
&\|Y_3\|_{r, \sigma, s-3s_1}
%
%
\le  {\rm diam}({\mathbb Y}_\sigma)\left\|\frac{1}{v}\right\|_{r, \sigma}\left\|Z_3\right\|_{r, \sigma, s-2s_1}
+e^{s_2}{\rm diam}({\mathbb Y}_\sigma)^2\left\|\frac{1}{v}\right\|_{r, \sigma}\left\|\frac{\partial_y\omega}{v}\right\|_{r, \sigma}
\left\|Z_2\right\|_{r, \sigma, s-s_1}
\nonumber\\
&\qquad+{\rm diam}({\mathbb Y}_\sigma)^2\left\|\frac{1}{v}\right\|_{r, \sigma}
\left(\left\|\frac{\partial_{\rm I}\omega}{v}\right\|_{r, \sigma}+e^{s_2}{\rm diam}({\mathbb Y}_\sigma)\left\|\frac{\partial_{\rm I} v}{v}\right\|_{r, \sigma} \left\|\frac{\partial_y\omega}{v}\right\|_{r, \sigma}\right)\left\|Z_1\right\|_{r, \sigma, s}
\end{align}
Multiplying the inequalities above by $\rho_*^{-1}$, $\tau_*^{-1}$, $t_*^{-1}$ respectively and taking the sum, we find~\eqref{bounds}, with 
\begin{align}
\frac{1}{\rho}&=\frac{1}{\rho_*}+e^{s_2}{\rm diam}({\mathbb Y}_\sigma)\left\|\frac{\partial_{\rm I} v}{v}\right\|_{r, \sigma}\frac{1}{\tau_*}+{\rm diam}({\mathbb Y}_\sigma)
\left(\left\|\frac{\partial_{\rm I}\omega}{v}\right\|_{r, \sigma}+e^{s_2}{\rm diam}({\mathbb Y}_\sigma)\left\|\frac{\partial_{\rm I} v}{v}\right\|_{r, \sigma} \left\|\frac{\partial_y\omega}{v}\right\|_{r, \sigma}\right)\frac{1}{t_*}\nonumber\\
\frac{1}{\tau}&=\frac{e^{s_2}}{\tau_*}+e^{s_2}{\rm diam}({\mathbb Y}_\sigma)\left\|\frac{\partial_y\omega}{v}\right\|_{r, \sigma}\frac{1}{t_*}\nonumber\\
\frac{1}{t}&=\frac{1}{t_*}\, . 
\end{align}
We recognise that, under conditions~\eqref{w*positive}, $\rho_*$, $\tau_*$, $t_*$ in~\eqref{bounds1} solve the equations above. 
$\qquad \square$

\begin{lemma}\label{Lie brackets}
Let  $w<u\le u_0$;
$Y\in {\cal O}^3_{u_0}$, $W\in {\cal O}^3_{u}$. Then
\begin{align}\VERT{\cal L}_Y[W]\VERT^{u_0-u+w}_{u-w}\le \VERT Y\VERT^{w}_{u-w}\VERT W\VERT ^{u_0-u+w}_{u}+\VERT W\VERT^{u_0-u+w}_{u-w}\VERT Y\VERT^{u_0-u+w}_{u_0}\, . \end{align}
\end{lemma}
{\bf Proof\ }
One has 
\begin{align}
\VERT{\cal L}_Y[W]\VERT^{u_0-u+w}_{u-w}&=\VERT J_W Y-J_Y W\VERT^{u_0-u+w}_{u-w}\nonumber\\
&\le  \VERT J_W Y\VERT^{u_0-u+w}_{u-w}+\VERT J_Y W\VERT^{u_0-u+w}_{u-w}
\end{align}
Now,  $(J_W Y)_i=\partial_{{\rm I}} W_i Y_1+\partial_{y} W_i Y_2+\partial_{\psi} W_i Y_3$, so, using Cauchy inequalities,
\begin{align}
 \|(J_W Y)_i\|_{u-w}&\le  \|\partial_{{\rm I}} W_i\|_{u-w}\| Y_1\|_{u-w}+\|\partial_{y} W_i\|_{u-w} \|Y_2\|_{u-w}+\|\partial_{\psi} W_i\|_{u-w} \|Y_3\|_{u-w}\nonumber\\
 &\le w_1^{-1}\| W_i\|_{u}\| Y_1\|_{u-w}+w_2^{-1}\|W_i\|_{u} \|Y_2\|_{u-w}+w_3^{-1}\|W_i\|_{u} \|Y_3\|_{u-w}\nonumber\\
 &=\VERT Y\VERT^w_{u-w}\|W_i\|_{u} 
\end{align}
Similarly,
\begin{align}\|(J_Y W)_i\|_{u-w}\le \VERT W\VERT_{u-w}^{u_0-u+w}\|Y_i\|_{u_0}\, .  \end{align}
Taking the $u_0-u+w$--weighted norms, the thesis follows. $\quad \square$
\begin{lemma}\label{iterateL}
Let $0<w<u\in {\mathbb R}^3$, $Y\in {\cal O}^3_{u+w}$, $W\in {\cal O}^3_{u}$. Then
\begin{align}\VERT{\cal L}^{{k}}_Y[W]\VERT^{w}_{u-w}\le 3^{{k}} {{k}}!\left(\VERT Y\VERT ^w_{u+w}\right)^{{k}}\VERT W\VERT^w _{u-w}\, . \end{align}
\end{lemma}

{\bf Proof\ } We apply Lemma~\ref{Lie brackets} with $W$ replaced by ${\cal L}^{i-1}_Y[W]$, 
$u$ replaced by $u-(i-1)w/{{k}}$, 
$w$ replaced by $w/{{k}}$ and, finally, $u_0=u+w$. With $\VERT\cdot\VERT_i^{w}=\VERT\cdot\VERT^{w}_{u-i\frac{w}{{{k}}}}$, $0\le i\le {{k}}$, so that $\VERT\cdot\VERT^w_0=\VERT\cdot\VERT^w_{u}$ and $\VERT\cdot\VERT^w_{{k}}=\VERT\cdot\VERT^w_{u-w}$,
\begin{align}
\VERT{\cal L}^i_Y[W]\VERT^{w+w/{{k}}}_i&=\left\VERT\left[Y, {\cal L}^{i-1}_Y[W]\right]\right\VERT^{w+w/{{k}}}_i\nonumber\\
&\le  \VERT Y\VERT^{w/{{k}}}_{i}\VERT{\cal L}^{i-1}_Y[W]\VERT^{w+w/{{k}}}_{i-1}+
\VERT Y\VERT^{w+w/{{k}}}_{u+w}\VERT{\cal L}^{i-1}_Y[W]\VERT^{w+w/{{k}}}_{i}\,.
\end{align}
Hence, de--homogenizating,
\begin{align}
\frac{{{k}}}{{{k}}+1}\VERT{\cal L}^i_Y[W]\VERT^{w}_i&\le  {{k}} \frac{{{k}}}{{{k}}+1}\VERT Y\VERT^{w}_{i}\VERT{\cal L}^{i-1}_Y[W]\VERT^{w}_{i-1}+
\frac{{{k}}^2}{({{k}}+1)^2}\VERT Y\VERT^{w}_{u+w}\VERT{\cal L}^{i-1}_Y[W]\VERT^{w}_{i}\nonumber\\
&\le \frac{{{k}}^2}{{{k}}+1}\left(1+\frac{1}{{{k}}+1}\right)\VERT Y\VERT^{w}_{u+w}\VERT{\cal L}^{i-1}_Y[W]\VERT^{w}_{i-1}
\end{align}
Eliminating the common factor $\frac{{{k}}}{{{k}}+1}$ and 
 iterating ${{k}}$ times from $i={{k}}$, by Stirling, we get
\begin{align}
\VERT {\cal L}^{{k}}_Y[W]\VERT^w _{u-w}
&\le   {{k}}^{{k}}\left(1+\frac{1}{{k}}\right)^{{k}}\left(\VERT Y\VERT ^w_{u+w}\right)^{{k}}\VERT W\VERT^w _{u-w}\nonumber\\
&\le  e^{{k}} {{k}}!\left(\VERT Y\VERT ^w_{u+w}\right)^{{k}}\VERT W\VERT^w _{u-w}\nonumber\\
&< 3^{{k}} {{k}}!\left(\VERT Y\VERT ^w_{u+w}\right)^{{k}}\VERT W\VERT^w _{u-w}
\end{align}
as claimed. $\quad \square$

\begin{proposition}\label{Lie Series}
Let $0<w<u$, $Y\in {\cal O}^3_{u+w}$, 
\begin{align}\label{q}q:=3\VERT Y\VERT^w_{u+w}<1\, . \end{align}
Then the Lie series $e^{{\cal L}_Y}$
defines an operator
\begin{align}e^{{\cal L}_Y}:\quad  {\cal O}^3_{u}\to {\cal O}^3_{u-w}\end{align}
and its tails
\begin{align}e^{{\cal L}_Y}_m=\sum_{{{k}}\ge m}\frac{{\cal L}^{{k}}_Y}{{{k}}!}\end{align}
verify
\begin{align}\left\VERT e^{{\cal L}_Y}_m W\right\VERT^w_{u-w}\le \frac{q^m}{1-q}\VERT W\VERT_{u}^w\qquad \forall\ W\in {\cal O}^3_{u}\, . \end{align}
\end{proposition}

\paragraph{\bf Proof of {Lemma}~\ref{iteration lemma}} We look for $Y$ such that $X_+:=e^{{\cal L}_Y} X$ has the desired properties. 
\begin{align}
e^{{\cal L}_Y} X&=e^{{\cal L}_Y}\left(N
+P\right)=N
+P+{\cal L}_Y N+e_2^{{\cal L}_Y} N
+e_1^{{\cal L}_Y}P\nonumber\\
&=N
+P-{\cal L}_N Y+P_+
\end{align}
with
$P_+=e_2^{{\cal L}_Y} N
+e_1^{{\cal L}_Y}P$.
We 
choose $Y$ so that the homological equation 
\begin{align}{\cal L}_N Y=P\end{align}
is satisfied. By Proposition~\ref{homeq1}, this equation has a solution $Y\in {\cal O}^3_{r, \sigma, s-3s_1}$ verifying 
\begin{align} q:=3\VERT Y\VERT^{w_*} _{r, \sigma, s-3s_1}&\le  
3{\rm diam}({\mathbb Y}_\sigma)\left\|\frac{1}{v}\right\|_{r, \sigma}\VERT P\VERT^w_u=Q \VERT P\VERT^w_u<1\, . \end{align}
By Proposition~\ref{Lie Series}, the Lie series $e^{{\cal L}_Y}$
defines an operator
\begin{align}e^{{\cal L}_Y}:\quad W\in {\cal O}_{u_*+w_*}\to {\cal O}_{u_*}\end{align}
and its tails $e^{{\cal L}_Y}_m$
verify
\begin{align}
\left\VERT e^{{\cal L}_Y}_m W\right\VERT^{w_*}_{u_*}&\le  \frac{q^m}{1-q}\VERT W\VERT^{w_*}_{u_*+w_*}\nonumber\\
&\le \frac{\left(Q\VERT P\VERT^w_u\right)^m}{1-Q\VERT P\VERT^w_u}\VERT W\VERT^{w_*}_{u_*+w_*}
\end{align}
for all $W\in {\cal O}^3_{u_*+w_*}$.
In particular, $e^{{\cal L}_Y}$ is well defined on ${\cal O}^3_{u}\subset {\cal O}^3_{u_*+w_*}$, hence $P_+\in {\cal O}^3_{u_*}$. The bounds on $P_+$ are obtained as follows. Using the homological equation, one finds
\begin{align}\label{pertbound1}
\VERT e_2^{{\cal L}_Y}N\VERT^{w_*}_{u_*}&= \left\VERT \sum_{{{k}}=1}^{\infty}\frac{{\cal L}^{{{k}}+1}_Y N}{({{k}}+1)!}\right\VERT^{w_*}_{u_*}\nonumber\\
&\le 
\sum_{{{k}}=1}^{\infty}\frac{1}{({{k}}+1)!}\left\VERT {\cal L}^{{{k}}+1}_Y N\right\VERT^{w_*}_{u_*}\nonumber\\
&=\sum_{{{k}}=1}^{\infty}\frac{1}{({{k}}+1)!}\left\VERT {\cal L}^{{{k}}}_Y P\right\VERT^{w_*}_{u_*}\nonumber\\
&\le  \sum_{{{k}}=1}^{\infty}\frac{1}{{{k}}!}\left\VERT {\cal L}^{{{k}}}_Y P\right\VERT^{w_*}_{u_*}\nonumber\\
&\le \frac{Q
\left(\VERT P\VERT^w_u\right)^2
}{1-Q
\left\VERT P\right\VERT^w_u}
\end{align}
The bound
\begin{align}\label{pertbound2}
 \VERT e_1^{{\cal L}_Y}P\VERT^{w_*}_{u_*}
\le \frac{Q
\left(\VERT P\VERT^w_u\right)^2
}{1-Q
\left\VERT P\right\VERT^w_u}\end{align}
is even more straightforward. $\quad \square$

\subsection{Proof of the Normal Form Theorem}\label{Proof of NFL}
The proof of {\sc nft} is obtained -- following~\cite{poschel93} -- via iterate applications of the Step Lemma. 
At the base step, we let\footnote{{With slight abuse of notations, here and during the proof of Theorem~\ref{Normal Form LemmaNEW}, the sub--fix $j$ will denote   the value of a given quantity at the $j^{\rm th}$ step  of the iteration.}}
\begin{align}X=X_0:=N+P_0\, ,\quad w=w_0:=(\rho,\tau, t)\, ,\quad u=u_0:=(r,\sigma, s)\end{align}
with $X_0=N+P_0\in {\cal O}^3_{u_0}$. We let
\begin{align}
Q_0:=3\,{\rm diam}({\mathbb Y}_\sigma)\left\|\frac{1}{v}\right\|_{r, \sigma}\end{align}
Conditions~\eqref{newcond1}--\eqref{newsmallness}   are implied by the assumptions~\eqref{theta1}--\eqref{NEWnewsmallness}. We then conjugate $X_0$ to $X_1=N+P_1\in {\cal O}^3_{u_1}$, where
\begin{align}u_1=(r-4\rho, \sigma-4\tau e^{s_2}, s-5t)=:(r_1, \sigma_1, s_1)\, . \end{align}
Then we have 
\begin{align}\label{base step}\VERT P_1\VERT_{u_1}^{w_0}\le 8 e^{s_2}Q_0  \left(\VERT P_0\VERT_{u_0}^{w_0}\right)^2\le \frac{1}{2} \VERT P_0\VERT_{u_0}^{w_0}\, .  \end{align}
We assume, inductively, that, for some $1\le j\le {{p}}$, 
we have  
\begin{align}\label{induction}X_{j}=N+P_{j}\in {\cal O}^3_{u_j}\, ,\qquad \VERT P_{j}\VERT_{u_j}^{w_0}<2^{-(j-1)}
\VERT P_{1}\VERT _{u_{1}}^{w_0}
\end{align} where \begin{align}\label{uj}u_j=(r_j, \sigma_j, s_j)\end{align} with
\begin{align}r_j:=r_1-4 (j-1)\frac{\rho}{ {{p}}}\, ,\quad \sigma_j:=\sigma_1-4 e^{s_2}(j-1)\frac{\tau}{ {{p}}}\, ,\quad s_j:=s_1-5 (j-1)\frac{t}{ {{p}}}\, . \end{align} The case $j=1$  trivially   reduces to the identity $\VERT P_{1}\VERT _{u_{1}}^{w_0}=\VERT P_{1}\VERT _{u_{1}}^{w_0}$. 
We  aim to apply {Lemma}~\ref{simplifiedsteplemma}  with $u=u_j$ as in~\eqref{uj} and \begin{align}w=w_1:=\frac{w_0}{{{p}}}\, ,\qquad \forall\ 1\le j\le {{p}}\, . \end{align}
Conditions~\eqref{newcond1},~\eqref{newcond2} and~\eqref{u+positive} are easily seen to be implied by~\eqref{theta1},~\eqref{NEWnewcond2},~\eqref{NEWu+positive} and the first condition in~\eqref{NEWnewsmallness} combined with the inequality ${{p}}\eta^2<1$, implied by the choice of ${{p}}$. We check condition~\eqref{newsmallness}.
By homogeneity,
\begin{align}\label{eq}\VERT P_j\VERT_{u_j}^{w_1}={{p}}\VERT P_j\VERT_{u_j}^{w_0}\le {{p}}\VERT P_1\VERT_{u_1}^{w_0}\le 8{{p}} e^{s_2}Q_0  \left(\VERT P_0\VERT_{u_0}^{w_0}\right)^2\end{align}
whence, using
\begin{align}
Q_j=3\,{\rm diam}({\mathbb Y}_{\sigma_j})\left\|\frac{1}{v}\right\|_{r_j, \sigma_j}\le Q_0\end{align}
we see that condition~\eqref{newsmallness} is met:
\begin{align}2 Q_j\VERT P_j\VERT_{u_j}^{w_1}\le 16{{p}} e^{s_2}Q^2_0  \left(\VERT P_0\VERT_{u_0}^{w_0}\right)^2<1\, . \end{align}
Then the Iterative Lemma can be applied and we get $X_{j+1}=N+P_{j+1}\in {\cal O}^3_{u_{j+1}}$, with \begin{align}\VERT P_{j+1}\VERT _{u_{j+1}}^{w_1}\le 8 e^{s_2}  Q_j \left(\VERT P_j\VERT _{u_j}^{w_1}\right)^2\le 8 e^{s_2}  Q_0 \left(\VERT P_j\VERT _{u_j}^{w_1}\right)^2\,.\end{align}
Using homogeneity again to the extreme sides of this inequality and combining it with~\eqref{induction},~\eqref{base step} and~\eqref{NEWnewsmallness}, we get
\begin{align}\label{iterate}\VERT P_{j+1}\VERT _{u_{j+1}}^{w_{0}}&\le  8 {{p}} e^{s_2}  Q_0 \left(\VERT P_j\VERT _{u_j}^{w_0}\right)^2
\le 8 {{p}} e^{s_2}  Q_0  \VERT P_1\VERT _{u_1}^{w_0}\VERT P_j\VERT _{u_j}^{w_0}\nonumber\\
&\le 
64{{p}}e^{2s_2}  Q_0^2\left(\VERT P_0\VERT _{u_0}^{w_0}\right)^2
\VERT P_j\VERT _{u_j}^{w_0}
\le
\frac{1}{2}\VERT P_j\VERT _{u_j}^{w_0}\nonumber\\
&<2^{-j}\VERT P_1\VERT _{u_1}^{w_0}\, . 
 \end{align}
After ${{p}}$ iterations, 
\begin{align}\VERT P_{{{p}}+1}\VERT _{u_{{{p}}+1}}^{w_{0}}<2^{-{{p}}}\VERT P_1\VERT _{u_1}^{w_0}<2^{-({{p}}+1)}\VERT P_0\VERT _{u_0}^{w_0}\end{align}
so we can take $X_\star=X_{{{p}}+1}$, $P_\star=P_{{{p}}+1}$,  $u_\star=u_{{{p}}+1}$. $\qquad \square$

\subsection{A generalisation when the dependence on $\psi$ is smooth}\label{A generalisation when the dependence}
\begin{definition}\rm We denote ${\cal C}^3_{u, \ell_*}$, with $u=(r, \sigma)$, the class of vector--fields 
\begin{align}  ({\rm I}, y, \psi):\ {\mathbb P}_{u}:={\mathbb I}_r \times{\mathbb Y}_\sigma \times {{\mathbb T}}\to X=(X_1, X_2, X_3)\in{\mathbb C}^3\qquad u=(r, \sigma)\end{align}
where each $X_i\in {\cal C}_{u, \ell_*}$, meaning that $X_i$ is $C^{\ell_*}$ in ${\mathbb P}:={\mathbb I} \times{\mathbb Y} \times {{\mathbb T}}$, $X_i(\cdot, \cdot, \psi)$ is holomorphic in ${\mathbb I}_r \times{\mathbb Y}_\sigma$ for each fixed $\psi$ in ${\mathbb T}$.
\end{definition}

\noindent
In this section we generalise {Theorem}~\ref{Normal Form Lemma} to the case that $X\in {\cal C}^3_{u, \ell_*}$. We use techniques going back to J. Nash and J. Moser~\cite{nash1956, moser1961, moser1962}.

\noindent
First of all, we need a different definition of norms\footnote{The series   in~\eqref{normfOLD} is in general diverging when $f\in {\cal C}_{u, \ell_*}$.} and, especially,  {\it smoothing} operators.

\paragraph{1. Generalised weighted norms}
We let
\begin{align}\label{normX}
\VERT X \VERT_{u, \ell}^{w}&:= \sum_i w^{-1}_i\|X_i\|_{u, \ell}\ ,\qquad 0\le \ell\le \ell_*\end{align}
 where
$w=(w_1$, $w_2$, $w_3)\in {\mathbb R}_+^3$
where, if
$f:\ {\mathbb P}_{r, \sigma}:={\mathbb I}_r \times{\mathbb Y}_\sigma \times {{\mathbb T}}\to {\mathbb C}$, 
then
\begin{align}\label{normf}\|f\|_{u}:=\sup_{{\mathbb I}_r\times{\mathbb Y}_\sigma\times{\mathbb T}}|f|
\ ,\quad \|f\|_{u, \ell}:=\max_{0\le j\le \ell}\{\|\partial_\varphi^j\, f\|_u\}\qquad u=(r, \sigma)\,.
\end{align}
Clearly, the class ${\cal O}^3_{r, \sigma, s}$ defined in Section~\ref{Weighted norms} is a proper subset of ${\cal C}^3_{u, \ell_*}$

\noindent
Observe that the norms~\eqref{normX} still verify monotonicity and 
homogeneity in
\eqref{monotonicity} and~\eqref{homogeneity}.

\paragraph{2. Smoothing}
We call {\it smoothing} a family of operators
\begin{align}T_K:\qquad f\in {\cal C}_{u, \ell_*}\to T_Kf\in {\cal C}_{u, \ell_*}\ ,\quad K\in {\mathbb N}\end{align}
verifying the following. Let $R_K:=I-T_K$. There exist $c_0>0$, $\delta\ge 0$ such that for all $f\in {\cal C}_{u, \ell_*}$, for all $K$, $0\le {j}\le \ell\le {\ell_*}$,

\begin{itemize}\item[\tiny\textbullet] $\|T_K\,f\|_{u, \ell}\le c_0\,K^{(\ell-{j}+\delta)}\|f\|_{u, {j}}\quad \forall\,0\le \ell\le \ell_*$\\
\item[\tiny\textbullet] $\|R_K\,f\|_{u, {j}}\le c_0\,K^{-(\ell-{j}-\delta)}\|f\|_{u, \ell}\quad \forall\,0\le \ell\le \ell_*$
\end{itemize}
As an example, as suggested in~\cite{arnold63}, one can  take 
\begin{align}T_K\,f({\rm I}, y, \psi):=\sum_{k\in {\mathbb Z}, |k|_1\le K} f_k({\rm I}, y)e^{{\rm i} k\psi}\end{align}
which, with the definitions~\eqref{normX}--\eqref{normf}, verifies the inequalities above with $\delta=2$.

\noindent
We name Generalised Normal Form Theorem ({\sc gnft}) the following

\begin{theorem}[{\sc gnft}]\label{Normal Form LemmaNEW}
Let $u=(r, \sigma)$; $X=N+P\in {\cal C}^3_{u, \ell_*}$, ${{p}}$, $\ell$, $K\in \natural$ and let $w_K=\left(\rho, \tau, \frac{1}{c_0\,K^{1+\delta}}\right)\in {\mathbb R}_+^3$  and assume that for some $s_1$, $s_2\in {\mathbb R}_+$, the following inequalities are satisfied.
Put
\begin{align}\label{delta1NEW}
 Q:=3\, e^{s_1}{\rm diam}({\mathbb Y}_\sigma)\left\|\frac{1}{v}\right\|_{r, \sigma}\end{align}
then assume:
 \begin{align}\label{NEWu+positiveNEW}
0<\rho<\frac{r}{8}\, ,\quad 0<\tau< e^{-s_2}\frac{\sigma}{8}\end{align}
and
 \begin{align}\label{NEWnewcond2NEW}\chi&:= \max\left\{\frac{{\rm diam}({\mathbb Y}_\sigma)}{s_1}\left\|\frac{\omega}{v}\right\|_{r, \sigma}\, ,\ \frac{{\rm diam}({\mathbb Y}_\sigma)}{s_2}\left\|\frac{\partial_y v}{v}\right\|_{r, \sigma}\right\}
\le 1	\\
\label{theta1NEW}\theta_1&:= 2\,e^{s_1+s_2}{\rm diam}({\mathbb Y}_\sigma)\left\|\frac{\partial_y\omega}{v}\right\|_{r, \sigma}c_0\,K^{1+\delta}\tau\le 1\\
\theta_2&:= \label{theta2NEW}4\,e^{s_1}{\rm diam}({\mathbb Y}_\sigma)\left\|\frac{\partial_{\rm I} v}{v}\right\|_{r, \sigma}\frac{\rho}{\tau}\le 1\\
\theta_3&:= \label{theta3NEW}8\,{e^{s_1}}{\rm diam}({\mathbb Y}_\sigma)
\left\|\frac{\partial_{\rm I}\omega}{v}\right\|_{r, \sigma}c_0\,K^{1+\delta}\rho\le 1\\\\
\label{NEWnewsmallnessNEW}\eta&:= 2^4\,e^{s_2}Q  \VERT P\VERT_{u}^{{w_K}}<\frac{1}{\sqrt{{{p}}}}\, . \end{align}
Then, with
\begin{align}u_*=(r_\star, \sigma_\star)\, ,\quad r_\star:=r-8\rho\, , \quad\sigma_\star=\sigma-8 e^{s_2}\tau\end{align}
 there exists  a real--analytic change of coordinates $\Phi_\star$
such that $X_\star:=\Phi_\star X\in {\cal C}^3_{u_\star, \ell_*}$ and
 $X_\star=N+P_\star$,
with
\begin{align}\label{P*NEW}\VERT P_\star\VERT^{w_K}_{u_\star}\le \max\left\{2^{-({{p}}+1)}\VERT P\VERT^{w_K}_{u}\, ,\ 2 c_0\,K^{-\ell+\delta}\VERT P\VERT^{w_K}_{u, \ell} \right\} \qquad \forall\ 0\le \ell\le \ell_*\,. \end{align}
\end{theorem}

\noindent
The result generalising {Lemma}~\ref{iteration lemma} is
\begin{lemma}\label{iteration lemmaNEW}
Let $X=N
+P\in {\cal C}^3_{u, \ell_*}$, with $u=(r, \sigma)$, $N$ as in~\eqref{NZ}, $\ell$, $K\in {\mathbb N}$. Assume~\eqref{existence}
and that 
$P$  is so small that

\begin{align}\label{smallcondNEW}
Q \VERT P\VERT^{w_K} _{u}<1\qquad Q:=3e^{s_1}{\rm diam}({\mathbb Y}_\sigma)\left\|\frac{1}{v}\right\|_{r, \sigma}\, ,\quad w_K=\left(\rho, \tau, \frac{1}{c_0\,K^{1+\delta}}\right)
\end{align}
Let $\rho_*$, $\tau_*$ be defined via 
\begin{align}\label{bounds1NEW}
\frac{1}{\rho_*}&=\frac{1}{\rho}-{\rm diam}({\mathbb Y}_\sigma)\left\|\frac{\partial_{\rm I} v}{v}\right\|_{r, \sigma}\left(\frac{e^{s_1}}{\tau}-e^{2s_1+s_2}{\rm diam}({\mathbb Y}_\sigma)
\left\|\frac{\partial_y\omega}{v}\right\|_{r, \sigma}c_0\,K^{1+\delta}
\right)\nonumber\\
&- {\rm diam}({\mathbb Y}_\sigma)
\left(e^{s_1}\left\|\frac{\partial_{\rm I}\omega}{v}\right\|_{r, \sigma}+e^{{2s_1+s_2}}{\rm diam}({\mathbb Y}_\sigma)\left\|\frac{\partial_{\rm I} v}{v}\right\|_{r, \sigma} \left\|\frac{\partial_y\omega}{v}\right\|_{r, \sigma}\right)c_0\,K^{1+\delta}\nonumber\\
\frac{1}{\tau_*}&=\frac{e^{-s_2}}{\tau}-e^{s_1}{\rm diam}({\mathbb Y}_\sigma)\left\|\frac{\partial_y\omega}{v}\right\|_{r, \sigma}c_0\,K^{1+\delta}
\end{align}
assume
\begin{align}\label{w*positiveNEW}\hat w_*=(\rho_*, \tau_*)\in {\mathbb R}_+^2\, ,\qquad u_*=(r-2\rho_*, \sigma-2\tau_*)\in {\mathbb R}_+^2\end{align} 
and put
\begin{align}w_{*, K}:=\left(\hat w_*, \frac{1}{c_0\,K^{1+\delta}}\right)\, . \end{align}
Then there exists $Y\in T_K{\cal C}^3_{u_*+\hat w_*, \ell_*}$  such that
$X_+:=e^{{\cal L}_Y}X\in {\cal C}^3_{u_*, \ell_*}$ and 
$X_+=N+
P_+$, with
\begin{align}\label{P+NEW} \VERT P_+\VERT^{w_{*, K}}_{u_*}
\le \frac{2Q
\left(\VERT P\VERT^{w_K}_u\right)^2
}{1-Q
\left\VERT P\right\VERT^{w_K}_u}+c\,K^{-\ell+\delta}\VERT P\VERT^ {w_{K}}_{u, \ell}\qquad \forall\ 0\le \ell\le \ell_*
\end{align}
\end{lemma}
The simplified form of
{Lemma}~\ref{iteration lemmaNEW}, corresponding to {Lemma}~\ref{simplifiedsteplemma}, is

\begin{lemma}[Generalised Step Lemma]\label{simplifiedsteplemmaNEW}
Assume~\eqref{existence} and replace~\eqref{smallcondNEW} and~\eqref{bounds1NEW}  with  
\begin{align}\label{newcond1NEW}&2\,e^{s_1+s_2}{\rm diam}({\mathbb Y}_\sigma)\left\|\frac{\partial_y\omega}{v}\right\|_{r, \sigma}c_0\,K^{1+\delta}\tau\le 1\nonumber\\
&4\,e^{s_1}{\rm diam}({\mathbb Y}_\sigma)\left\|\frac{\partial_{\rm I} v}{v}\right\|_{r, \sigma}\frac{\rho}{\tau}\le 1\nonumber\\
&8\,{e^{s_1}}{\rm diam}({\mathbb Y}_\sigma)
\left\|\frac{\partial_{\rm I}\omega}{v}\right\|_{r, \sigma}c_0\,K^{1+\delta}\rho\le 1\\
\label{u+positiveNEW}
&0<\rho<\frac{r}{4}\, ,\quad 0<\tau<\frac{\sigma}{4}e^{-s_2}\\
\label{newsmallnessNEW} 
&2Q\VERT P\VERT ^{w_K}_u<1
\end{align}
then $X_+=N+P_+\in {\cal C}^3_{u_+, \ell_*}$ and \begin{align}\label{finalineqNEW}\VERT P_+\VERT_{u_+}^{w_K}\le 8 e^{s_2}Q
(\VERT P \VERT_u^{w_K})^2+c\,K^{-\ell+\delta}\VERT P\VERT^ {w_{K}}_{u, \ell}
\end{align}
with
\begin{align}u_+:=(r-4\rho, \sigma-4\tau e^{s_2})\, . \end{align}
\end{lemma}
{\bf Proof\ }The  inequalities in~\eqref{newcond1NEW} guarantee 
\begin{align}\frac{1}{\rho_*}\ge \frac{1}{2\rho}\, ,\quad \frac{1}{\tau_*}\ge \frac{e^{-s_2}}{2\tau}\end{align}
whence
\begin{align}w_{*, K}<2 e^{s_2} w_K\, ,\qquad u_*\ge u_+>0\, . \end{align}
Then
~\eqref{finalineqNEW} is implied by~\eqref{P+NEW}, monotonicity and homogeneity and the inequality in~\eqref{newsmallnessNEW}. $\quad\square$

\vskip.2in
\noindent
 Let now
${\cal F}_{v, \omega}$ and ${\cal G}_{v, \omega}$ be as in~\eqref{FandG}. First of all, observe that
${\cal F}_{v, \omega}$, ${\cal G}_{v, \omega}$ take $ T_K{\cal C}_{u, \ell_*}$ to  itself. Moreover, generalising Lemma~\ref{estimates},

\begin{lemma}\label{estimatesNEW}
If inequalities~\eqref{existence} hold,
then
\begin{align}{\cal F}_{v, \omega}\, ,\ {\cal G}_{v, \omega}:\quad {\cal C}_{u, \ell_*}\to {\cal C}_{u, \ell_*}\end{align}
and
\begin{align}\|{\cal F}_{v, \omega}[g]\|_{r, \sigma}\le e^{s_1}{\rm diam}({\mathbb Y}_\sigma)\left\|\frac{g}{v}\right\|_{r, \sigma}\, ,\quad \|{\cal G}_{v, \omega}[g]\|_{r, \sigma}\le e^{s_1+s_2}\,{\rm diam}({\mathbb Y}_\sigma) \left\|\frac{g}{v}\right\|_{r, \sigma}\, . \end{align}
 \end{lemma}

 \begin{proposition}\label{homeq1NEW}
Let 
\begin{align}\label{NZNEW}N=(0, v({\rm I}, y), \omega({\rm I}, y))\, ,\qquad Z=(Z_1({\rm I}, y, \psi), Z_2({\rm I}, y, \psi), Z_3({\rm I}, y, \psi))\end{align}
belong to ${\cal C}^3_{u, \ell_*}$ and assume~\eqref{existence}.
Then  the ``homological equation''
\begin{align}\label{homeqNEW}{\cal L}_N[Y]=Z\end{align}
has a solution $Y\in {\cal C}_{u, \ell_*}$ verifying
\begin{align}\label{boundsNEW}
\VERT Y\VERT_{u}^{\rho_*, \tau_*, t_*}\le e^{s_1}{\rm diam}({\mathbb Y}_\sigma)\left\|\frac{1}{v}\right\|_{u}\VERT Z\VERT_{u}^{\rho, \tau, t}\quad u=(r, \sigma)
\end{align}
with $\rho_*$, $\tau_*$, $t_*$ defined via
\begin{align}\label{bounds1NEW*}
\frac{1}{\rho_*}&=\frac{1}{\rho}-{\rm diam}({\mathbb Y}_\sigma)\left\|\frac{\partial_{\rm I} v}{v}\right\|_{u}\left(\frac{e^{s_1}}{\tau}-e^{2s_1+s_2}{\rm diam}({\mathbb Y}_\sigma)
\left\|\frac{\partial_y\omega}{v}\right\|_{u}\frac{1}{t}
\right)\nonumber\\
&- {\rm diam}({\mathbb Y}_\sigma)
\left(e^{s_1}\left\|\frac{\partial_{\rm I}\omega}{v}\right\|_{u}+e^{{2s_1+s_2}}{\rm diam}({\mathbb Y}_\sigma)\left\|\frac{\partial_{\rm I} v}{v}\right\|_{u} \left\|\frac{\partial_y\omega}{v}\right\|_{u}\right)\frac{1}{t}\nonumber\\
\frac{1}{\tau_*}&=\frac{e^{-s_2}}{\tau}-e^{s_1}{\rm diam}({\mathbb Y}_\sigma)\left\|\frac{\partial_y\omega}{v}\right\|_{u}\frac{1}{t}\nonumber\\
t_*&=t
\end{align}
and provided that
\begin{align}\label{w*positiveNEW*}
(\rho_*, \tau_*)\in {\mathbb R}_+^2\, . 
\end{align}
 In particular, if $Z\in T_{K}{\cal C}^3_{u, \ell_*}$ for some $K\in {\mathbb N}$, then also $Y\in T_{K}{\cal C}^3_{u, \ell_*}$.
\end{proposition}
{\bf Proof\ }The solution~\eqref{Yi} 
satisfies
\begin{align}
\|Y_1\|_{u}&\le  e^{s_1}{\rm diam}({\mathbb Y}_\sigma)\left\|\frac{1}{v}\right\|_{u}\left\|Z_1\right\|_{u}\nonumber\\
\|Y_2\|_{u}&\le  e^{s_1+s_2}{\rm diam}({\mathbb Y}_\sigma)\left\|\frac{1}{v}\right\|_{u}\left\|Z_2\right\|_{u}+e^{2s_1+s_2}{\rm diam}({\mathbb Y}_\sigma)^2\left\|\frac{1}{v}\right\|_{u}\left\|\frac{\partial_{\rm I} v}{v}\right\|_{u}\left\|Z_1\right\|_{u}\nonumber\\
\|Y_3\|_{u}&\le  e^{s_1} {\rm diam}({\mathbb Y}_\sigma)\left\|\frac{1}{v}\right\|_{u}\left\|Z_3\right\|_{u}
+e^{2s_1+s_2}{\rm diam}({\mathbb Y}_\sigma)^2\left\|\frac{1}{v}\right\|_{u}\left\|\frac{\partial_y\omega}{v}\right\|_{u}
\left\|Z_2\right\|_{u}
\nonumber\\
&+ {\rm diam}({\mathbb Y}_\sigma)^2\left\|\frac{1}{v}\right\|_{u}
\left(e^{2s_1}\left\|\frac{\partial_{\rm I}\omega}{v}\right\|_{u}+e^{{3s_1+s_2}}{\rm diam}({\mathbb Y}_\sigma)\left\|\frac{\partial_{\rm I} v}{v}\right\|_{u} \left\|\frac{\partial_y\omega}{v}\right\|_{u}\right)\left\|Z_1\right\|_{u}
\end{align}
Multiplying the inequalities above by $\rho_*^{-1}$, $\tau_*^{-1}$, $t_*^{-1}$ respectively and taking the sum, we find~\eqref{boundsNEW}, with 
\begin{align}
\frac{1}{\rho}&=\frac{1}{\rho_*}+e^{s_1+s_2}{\rm diam}({\mathbb Y}_\sigma)\left\|\frac{\partial_{\rm I} v}{v}\right\|_{u}\frac{1}{\tau_*}+{\rm diam}({\mathbb Y}_\sigma)
\left(e^{s_1}\left\|\frac{\partial_{\rm I}\omega}{v}\right\|_{u}+e^{{2s_1+s_2}}{\rm diam}({\mathbb Y}_\sigma)\left\|\frac{\partial_{\rm I} v}{v}\right\|_{u} \left\|\frac{\partial_y\omega}{v}\right\|_{u}\right)\frac{1}{t_*}\nonumber\\
\frac{1}{\tau}&=\frac{e^{s_2}}{\tau_*}+e^{s_1+s_2}{\rm diam}({\mathbb Y}_\sigma)\left\|\frac{\partial_y\omega}{v}\right\|_{u}\frac{1}{t_*}\nonumber\\
\frac{1}{t}&=\frac{1}{t_*}\, . 
\end{align}
We recognise that, under conditions~\eqref{w*positiveNEW*},  $\rho_*$, $\tau_*$, $t_*$ in~\eqref{bounds1NEW*} solve the equations above. Observe that if $Z\in T_{K}{\cal C}^3_{u, \ell_*}$, then also $Y\in T_{K}{\cal C}^3_{u, \ell_*}$, as ${\cal F}_{v, \omega}$, ${\cal G}_{v, \omega}$ do so.
$\qquad \square$

\begin{lemma}\label{Lie bracketsNEW}
Let  $u_0\ge u>w\in {\mathbb R}_+^2\times \{0\}$;
$Y\in T_K{\cal C}^3_{u_0, \ell_*}$, $W\in T_K{\cal C}^3_{u, \ell_*}$. Put $w_K:=\left(w_1, w_2, \frac{1}{c_0\,K^{1+\delta}}\right)$. Then
\begin{align}\VERT{\cal L}_Y[W]\VERT^{u_0-u+w_K}_{u-w}\le \VERT Y\VERT^{w_K}_{u-w}\VERT W\VERT ^{u_0-u+w_K}_{u}+\VERT W\VERT^{u_0-u+w_K}_{u-w}\VERT Y\VERT ^{u_0-u+w_K}_{u_0}\, . \end{align}
\end{lemma}
{\bf Proof\ }By Cauchy inequalities, the definitions~\eqref{normX}--\eqref{normf} and the smoothing properties,
\begin{align}
 \|(J_W Y)_i\|_{u-w}&\le  \|\partial_{{\rm I}} W_i\|_{u-w}\| Y_1\|_{u-w}+\|\partial_{y} W_i\|_{u-w} \|Y_2\|_{u-w}+\|\partial_{\psi} W_i\|_{u-w} \|Y_3\|_{u-w}\nonumber\\
 &\le w_1^{-1}\| W_i\|_{u}\| Y_1\|_{u-w}+w_2^{-1}\|W_i\|_{u} \|Y_2\|_{u-w}+\|W_i\|_{u, 1} \|Y_3\|_{u-w}\nonumber\\
 &\le w_1^{-1}\| W_i\|_{u}\| Y_1\|_{u-w}+w_2^{-1}\|W_i\|_{u} \|Y_2\|_{u-w}+c_0\,K^{1+\delta}\|W_i\|_{u} \|Y_3\|_{u-w}\nonumber\\
 &=\VERT Y\VERT_{u-w}^{w_K}\|W_i\|_{u}
\end{align}
Similarly,
\begin{align}\|(J_Y W)_i\|_{u-w}\le \VERT W\VERT_{u-w}^{u_0-u+w_K}\|Y_i\|_{u_0}\, .  \end{align}
Taking the $u_0-u+w_K$--weighted norms, the thesis follows. $\quad \square$

\begin{lemma}\label{iterateLNEW}
Let $0<w<u\in {\mathbb R}^2_+\times\{0\}$, $w_K:=\left(w_1, w_2, \frac{1}{c_0\,K^{1+\delta}}\right)$;
$Y\in T_K{\cal C}^3_{u+w, \ell_*}$, $W\in T_K{\cal C}^3_{u, \ell_*}$. Then
\begin{align}\VERT{\cal L}^n_Y[W]\VERT^{w_K}_{u-w}\le 3^n n!\left(\VERT Y\VERT ^{w_K}_{u+w}\right)^n\VERT W\VERT^{w_K} _{u-w}\, . \end{align}
\end{lemma}

{\bf Proof\ } The proof copies the one of Lemma~\ref{iterateL}, up to invoke Lemma~\ref{Lie bracketsNEW} at the place of Lemma~\ref{Lie brackets} and hence replace the $w$'s ``up'' with $w_K$. $\qquad \square$

\begin{proposition}\label{Lie SeriesNEW}
Let $0<w<u\in {\mathbb R}_+^2\times\{0\}$, $w_K:=\left(w_1, w_2, \frac{1}{c_0\,K^{1+\delta}}\right)$, $Y\in T_K{\cal C}^3_{u+w, \ell_*}$, 
\begin{align}\label{q}q:=3\VERT Y\VERT^{w_K}_{u+w}<1\, . \end{align}
Then the Lie series $e^{{\cal L}_Y}$
defines an operator
\begin{align}e^{{\cal L}_Y}:\quad  T_K{\cal C}^3_{u, \ell_*}\to T_K{\cal C}^3_{u-w, \ell_*}\end{align}
and its tails
\begin{align}e^{{\cal L}_Y}_m=\sum_{{{k}}\ge m}\frac{{\cal L}^n_Y}{{{k}}!}\end{align}
verify
\begin{align}\left\VERT e^{{\cal L}_Y}_m W\right\VERT^{w_K}_{u-w}\le \frac{q^m}{1-q}\VERT W\VERT_{u}^{w_K}\qquad \forall\ W\in T_K{\cal C}^3_{u, \ell_*}\, . \end{align}
\end{proposition}
\paragraph{\bf Proof of {Lemma}~\ref{iteration lemmaNEW}} {All the remarks before Lemma}~\ref{estimates} 
{continue holding also in this case, except for the fact that,}  differently from {Lemma}~\ref{iteration lemma} here we need a ``ultraviolet cut--off'' of the perturbing term. Namely, we split
\begin{align}
e^{{\cal L}_Y} X&=e^{{\cal L}_Y}\left(N
+P\right)=N
+P+{\cal L}_Y N+e_2^{{\cal L}_Y} N
+e_1^{{\cal L}_Y}P\nonumber\\
&=N
+T_KP-{\cal L}_N Y+P_+
\end{align}
with
$P_+=e_2^{{\cal L}_Y} N
+e_1^{{\cal L}_Y}P+R_KP$.
We 
choose $Y$ so that the homological equation 
\begin{align}{\cal L}_N Y=T_KP\end{align}
is satisfied. By Proposition~\ref{homeq1NEW}, this equation has a solution $Y\in T_K{\cal C}^3_{u, \ell_*}$ verifying 
\begin{align} q:=3\VERT Y\VERT^{w_*} _{u}&\le  
3{e^{s_1}}{\rm diam}({\mathbb Y}_\sigma)\left\|\frac{1}{v}\right\|_{u}\VERT P\VERT^{w_K}_u=Q \VERT P\VERT^{w_K}_u<1\, . \end{align}
with $w_*=(\rho_*, \tau_*, t_*)$ as in~\eqref{bounds1NEW*}. As
$t_*=t=\frac{1}{c_0\,K^{1+\delta}}$,  
We let \begin{align}w_{*, K}:=w_*\ ,\quad \hat w_{*}:=(\rho_*, \tau_*)\end{align}
with
$(\rho_*, \tau_*)$ as in~\eqref{bounds1NEW}.
By Proposition~\ref{Lie SeriesNEW}, the Lie series $e^{{\cal L}_Y}$
defines an operator
\begin{align}e^{{\cal L}_Y}:\quad W\in T_K{\cal C}_{u_*+\hat w_*, \ell_*}\to T_K{\cal C}_{u_*, \ell_*}\end{align}
and its tails $e^{{\cal L}_Y}_m$
verify
\begin{align}
\left\VERT e^{{\cal L}_Y}_m W\right\VERT^ {w_{*, K}}_{u_*}\le\frac{\left(Q\VERT P\VERT^{w_K}_u\right)^m}{1-Q\VERT P\VERT^{w_K}_u}\VERT W\VERT^ {w_{*, K}}_{u_*+\hat w_*}
\end{align}
for all $W\in T_K{\cal C}^3_{u_*+\hat w_*, \ell_*}$.
In particular, $e^{{\cal L}_Y}$ is well defined on $T_K{\cal C}^3_{u, \ell_*}\subset T_K{\cal C}^3_{u_*+\hat w_*, \ell_*}$, hence $P_+\in {\cal C}^3_{u_*, \ell_*}$. The bounds on $P_+$ are obtained as follows. The terms $\VERT e_2^{{\cal L}_Y}N\VERT^ {w_{*, K}}_{u_*}$ and $\VERT e_1^{{\cal L}_Y}P\VERT^{w_{*, K}}_{u_*}$ are treated quite similarly as~\eqref{pertbound1} and~\eqref{pertbound2}:
\begin{align}
\VERT e_2^{{\cal L}_Y}N\VERT^ {w_{*, K}}_{u_*}\le \frac{Q
\left(\VERT P\VERT^{w_K}_u\right)^2
}{1-Q
\left\VERT P\right\VERT^{w_K}_u}\ ,\quad \VERT e_1^{{\cal L}_Y}P\VERT^{w_{*, K}}_{u_*}
\le \frac{Q
\left(\VERT P\VERT^{w_K}_u\right)^2
}{1-Q
\left\VERT P\right\VERT^{w_K}_u}
\end{align}
The moreover, here we have the term $R_KP$, which is obviously bounded as
\begin{align}
\VERT R_KP\VERT^ {w_{*, K}}_{u_*}\le c\,K^{-\ell+\delta}\VERT P\VERT^ {w_{*, K}}_{u_*, \ell}\le c\,K^{-\ell+\delta}\VERT P\VERT^ {w_{K}}_{u, \ell}\, . \quad \square
\end{align}

\noindent
We are finally ready for the

\paragraph{Proof of {Theorem~\ref{Normal Form LemmaNEW}}}

Analogously as in the proof of {\sc nft}, we proceed by  iterate applications of the Generalised Step Lemma. 
At the base step, we let
\begin{align}X=X_0:=N+P_0\, ,\quad w_0:=w_{0, K}:=\left(\rho,\tau, \frac{1}{c_0 K^{1+\delta}}\right)\, ,\quad u_0:=(r,\sigma)\end{align}
with $X_0=N+P_0\in {\cal C}^3_{u_0, \ell_*}$. We let
\begin{align}
Q_0:=3\,e^{s_1}{\rm diam}({\mathbb Y}_\sigma)\left\|\frac{1}{v}\right\|_{u_0}\end{align}
Conditions~\eqref{newcond1NEW}--\eqref{newsmallnessNEW}   are implied by the assumptions~\eqref{delta1NEW}--\eqref{NEWnewsmallnessNEW}. We then conjugate $X_0$ to $X_1=N+P_1\in {\cal C}^3_{u_1, \ell_*}$, where
\begin{align}u_1=(r-4\rho, \sigma-4\tau e^{s_2})=:(r_1, \sigma_1)\, . \end{align}
Then we have 
\begin{align}\VERT P_1\VERT_{u_1}^{w_0}\le 
8 e^{s_2}Q_0  \left(\VERT P_0\VERT_{u_0}^{w_0}\right)^2+ c_0\,K^{-\ell+\delta}
\VERT P_0\VERT_{u_0, \ell}^{w_0}
\, .  \end{align}
If $8 e^{s_2}Q_0  \left(\VERT P_0\VERT_{u_0}^{w_0}\right)^2\le
c_0\,K^{-\ell+\delta}
\VERT P_0\VERT_{u_0, \ell}^{w_0} $,
the proof finishes here. So, we assume the opposite inequality, which gives
\begin{align}\label{base stepNEW}\VERT P_1\VERT_{u_1}^{w_0}\le 
16 e^{s_2}Q_0  \left(\VERT P_0\VERT_{u_0}^{w_0}\right)^2\le \frac{1}{2} \VERT P_0\VERT_{u_0}^{w_0}
\, .  \end{align}
We assume, inductively, that, for some $1\le j\le {{p}}$, 
we have  
\begin{align}\label{inductionNEW}X_{j}=N+P_{j}\in {\cal C}^3_{u_j, \ell_*}\, ,\qquad \VERT P_{j}\VERT_{u_j}^{w_0}<2^{-(j-1)}
\VERT P_{1}\VERT _{u_{1}}^{w_0}
\end{align} where \begin{align}\label{ujNEW}u_j=(r_j, \sigma_j)\end{align} with
\begin{align}r_j:=r_1-4 (j-1)\frac{\rho}{ {{p}}}\, ,\quad \sigma_j:=\sigma_1-4 e^{s_2}(j-1)\frac{\tau}{ {{p}}}\, . \end{align} The case $j=1$  is trivially  true because it is the identity $\VERT P_{1}\VERT _{u_{1}}^{w_0}=\VERT P_{1}\VERT _{u_{1}}^{w_0}$. 
We  aim to apply {Lemma}~\ref{simplifiedsteplemmaNEW}  with $u=u_j$ as in~\eqref{ujNEW} and \begin{align}w=w_1:=\frac{w_0}{{{p}}}\, ,\qquad \forall\ 1\le j\le {{p}}\, . \end{align}
Conditions~\eqref{newcond1NEW} and~\eqref{u+positiveNEW} correspond to~\eqref{NEWnewcond2NEW}--\eqref{theta3NEW}, while~\eqref{newsmallnessNEW} is implied by~\eqref{NEWnewsmallnessNEW}. We check condition~\eqref{newsmallnessNEW}.
By homogeneity,
\begin{align}\label{eqNEW}\VERT P_j\VERT_{u_j}^{w_1}={{p}}\VERT P_j\VERT_{u_j}^{w_0}\le {{p}}\VERT P_1\VERT_{u_1}^{w_0}\le 16{{p}} e^{s_2}Q_0  \left(\VERT P_0\VERT_{u_0}^{w_0}\right)^2\end{align}
whence, using
\begin{align}
Q_j=3\,{\rm diam}({\mathbb Y}_{\sigma_j})\left\|\frac{1}{v}\right\|_{r_j, \sigma_j}\le Q_0\end{align}
we see that condition~\eqref{newsmallness} is met:
\begin{align}2 Q_j\VERT P_j\VERT_{u_j}^{w_1}\le 32\,{{p}} e^{s_2}Q^2_0  \left(\VERT P_0\VERT_{u_0}^{w_0}\right)^2<1\, . \end{align}
Then the Iterative Lemma can be applied and we get $X_{j+1}=N+P_{j+1}\in {\cal C}^3_{u_{j+1}, \ell_*}$, with \begin{align}\VERT P_{j+1}\VERT _{u_{j+1}}^{w_1}\le 8 e^{s_2}  Q_j \left(\VERT P_j\VERT _{u_j}^{w_1}\right)^2\le 8 e^{s_2}  Q_0 \left(\VERT P_j\VERT _{u_j}^{w_1}\right)^2\end{align}
Using homogeneity again to the extreme sides of this inequality and combining it with~\eqref{inductionNEW},~\eqref{base stepNEW} and~\eqref{NEWnewsmallnessNEW}, we get
\begin{align}\label{iterateNEW}\VERT P_{j+1}\VERT _{u_{j+1}}^{w_{0}}&\le  8 {{p}} e^{s_2}  Q_0 \left(\VERT P_j\VERT _{u_j}^{w_0}\right)^2
\le 8 {{p}} e^{s_2}  Q_0  \VERT P_1\VERT _{u_1}^{w_0}\VERT P_j\VERT _{u_j}^{w_0}\nonumber\\
&\le 
128\,{{p}}e^{2s_2}  Q_0^2\left(\VERT P_0\VERT _{u_0}^{w_0}\right)^2
\VERT P_j\VERT _{u_j}^{w_0}
\le
\frac{1}2\VERT P_j\VERT _{u_j}^{w_0}\nonumber\\
&<2^{-j}\VERT P_1\VERT _{u_1}^{w_0}\, . 
 \end{align}
After ${{p}}$ iterations, 
\begin{align}\VERT P_{{{p}}+1}\VERT _{u_{{{p}}+1}}^{w_{0}}<2^{-{{p}}}\VERT P_1\VERT _{u_1}^{w_0}<2^{-({{p}}+1)}\VERT P_0\VERT _{u_0}^{w_0}\end{align}
so we can take $X_\star=X_{{{p}}+1}$, $P_\star=P_{{{p}}+1}$,  $u_\star=u_{{{p}}+1}$. $\qquad \square$

\section{Symplectic tools}\label{Tools}
In this section we describe various sets of canonical coordinates that are needed to our application. We remark that  during the proof of Theorem B, we shall not use any of such sets completely, but rather a ``mix'' of action--angle and regularising coordinates, described below.

\subsection{Starting coordinates}

 We begin with the coordinates

 \begin{align}\label{coord}
\left\{\begin{array}{l}\displaystyle {\rm C}=\|{\mathbf x}\times {\mathbf y}+{\mathbf x}'\times {\mathbf y}'\|\\\\
\displaystyle {\rm G}=\|{\mathbf x}\times {\mathbf y}\|\\\\
\displaystyle {\rm R}=\frac{\mathbf y'\cdot \mathbf x'}{\|\mathbf x'\|}\\\\
\displaystyle \Lambda= \sqrt{ a}
\end{array}\right.\qquad\qquad \left\{\begin{array}{l}\displaystyle \gamma =\alpha_{\mathbf k}(\mathbf i, \mathbf x')+\frac{\pi}{2}\\\\
\displaystyle  {\rm g}=\alpha_{\mathbf k}({\mathbf x'},\mathbf P)+\pi\\\\
\displaystyle  {\rm r}=\|\mathbf x'\|\\\\
\displaystyle \ell={\rm mean\ anomaly\ of\ {\mathbf x}\ in\ \mathbb E}
\end{array}\right.
\end{align}
where:

 \begin{itemize}
   \item[{\tiny\textbullet}] ${\mathbf i}=\left(
\begin{array}{lll}
1\\
0\\
0
\end{array}
\right)$, $ {\mathbf j}=\left(
\begin{array}{lll}
0\\
1\\
0
\end{array}
\right)$ is a hortonormal frame in ${\mathbb R}^2\times\{\mathbf 0\}$ and ${\mathbf k}={\mathbf i}\times {\mathbf j}$ (``$\times$'' denoting, as usual, the ``skew--product''); 
  \item[{\tiny\textbullet}] after 
 fixing  a set of values of $({\mathbf y},  {\mathbf x})$ where the Kepler Hamiltonian~\eqref{Kep}
 takes negative values, ${\mathbb E}$ denotes the elliptic orbit with initial values $({\mathbf y}_0, {\mathbf x}_0)$ in such set;
 \item[{\tiny\textbullet}]
 $a$  is the semi--major axis of ${\mathbb E}$; 
  \item[{\tiny\textbullet}] ${\mathbf P}$, with $\|{\mathbf P}\|=1$, the direction of the perihelion  of ${\mathbb E}$, assuming  ${\mathbb E}$ is not a circle;
    \item[{\tiny\textbullet}] $\ell$ is the mean anomaly of ${\mathbf x}$  on ${\mathbb E}$, defined, mod $2\pi$, as the area of the elliptic sector spanned  from ${\mathbf P}$ to ${\mathbf x}$, normalized to $2\pi$;
     \item[{\tiny\textbullet}]    $\alpha_{\mathbf w}({\mathbf u}, {\mathbf v})$ is the oriented angle  from ${\mathbf u}$ to ${\mathbf v}$ relatively to the positive orientation established by ${\mathbf w}$, if ${\mathbf u}$, ${\mathbf v}$ and ${\mathbf w}\in {\mathbb R}^3\setminus\{\mathbf 0\}$, with ${\mathbf u}$, ${\mathbf v}\perp{\mathbf w}$.
    \end{itemize}
The canonical\footnote{Namely, the change of coordinate~\eqref{coord} satisfies $\sum_{i=1}^2(d{\mathbf y}_i\wedge d{\mathbf x}_i+d{\mathbf y}_i'\wedge d{\mathbf x}_i')=dC\wedge d\gamma +d{\rm G}\wedge d{\rm g}+d{\rm R}\wedge d{\rm r}+d\Lambda\wedge d\ell$.} character of the coordinates~\eqref{coord} has been discussed, in a more general setting, in~\cite{pinzari19}. The shifts $\frac{\pi}{2}$ and $\pi$ in~\eqref{coord} serve only to be consistent with the spatial coordinates of~\cite{pinzari19}.

\subsection{Energy--time coordinates}\label{Energy--time coordinates}

We now describe the ``energy--time'' change of coordinates
\begin{align}\label{transf}\phi_{\rm et}:\qquad ({\cal R}, {\cal E}, {\rm r}, \tau)\to ({\rm R}, {\rm G}, {\rm r}, {\rm g})=({\cal R}+\rho({\cal E}, {\rm r}, \tau),\ \widetilde{\rm G}({\cal E}, {\rm r}, \tau),\ {\rm r},\ \widetilde{\rm g}({\cal E}, {\rm r}, \tau))\end{align}
which integrates the function ${\rm E}({\rm r}, {\rm G}, {\rm g})$ in~\eqref{E}, where ${\cal E}$ (``energy'') denotes the generic level--set of ${\rm E}$, while $\tau$ is its conjugated (``time'') coordinate. The domain of the coordinates~\eqref{transf} is
\begin{align}\label{range}{\cal R}\in {\mathbb R}\, ,\quad 0\le {\rm r}<2\, ,\quad -{\rm r}<{\cal E}<1+\frac{{\rm r}^2}{4}\, ,\quad \tau\in {\mathbb R}\, ,\quad {\cal E}\notin\{{\rm r}, 1\}\, . \end{align}
The extremal values of ${\cal E}$ are taken to be the minimum and the maximum of the function ${\rm E}$ for $0\le {\rm r}<2$. The values ${\rm r}$ and $1$ have been excluded because they correspond, in the $({\rm g}, {\rm G})$--plane, to the curves ${\cal S}_0({\rm r})$ and ${\cal S}_1({\rm r})$ in Figure~\ref{figure1}, where periodic motions do not exist.

\noindent
The functions  $\widetilde{\rm G}({\cal E}, {\rm r}, \cdot)$, $\widetilde{\rm g}({\cal E}, {\rm r}, \cdot)$ and $\rho({\cal E}, {\rm r}, \cdot)$ appearing in~\eqref{transf} are, respectively, $2\tau_{\rm p}$ periodic, $2\tau_{\rm p}$ periodic, $2\tau_{\rm p}$ quasi--periodic,  meaning that they satisfy
\begin{align}\label{full real}
{\cal P}_{er}:\quad \left\{\begin{array}{l}\displaystyle \widetilde{\rm G}({\cal E}, {\rm r}, \tau+2j\tau_{\rm p})=\widetilde{\rm G}({\cal E}, {\rm r}, \tau)\\\\ \displaystyle \widetilde{\rm g}({\cal E}, {\rm r}, \tau+2j\tau_{\rm p})=\widetilde{\rm g}({\cal E}, {\rm r}, \tau)\\\\
\displaystyle  \rho({\cal E}, {\rm r}, \tau+2j\tau_{\rm p})=\rho({\cal E}, {\rm r}, \tau)+2j\rho({\cal E}, {\rm r}, \tau_{\rm p})
\end{array}\right.\qquad \forall\ \tau\in {\mathbb R}\, ,\ \forall\ j\in {\mathbb Z}
\end{align}
with $\tau_{\rm p}=\tau_{\rm p}({\cal E}, {\rm r})$ the period, defined below. Note that one can find a unique splitting
\begin{align}\label{split}\rho({\cal E}, {\rm r}, \tau)={\cal B}({\cal E}, {\rm r}) \tau+\widetilde\rho({\cal E}, {\rm r}, \tau)\end{align}
such that $\widetilde\rho({\cal E}, {\rm r}, \cdot)$ is $2\tau_{\rm p}$--periodic. It is obtained taking
\begin{align}\label{B}{\cal B}({\cal E}, {\rm r})=\frac{\rho({\cal E}, {\rm r}, \tau_{\rm p}({\cal E}, {\rm r}))}{\tau_{\rm p}({\cal E}, {\rm r})}\, ,\quad \widetilde\rho({\cal E}, {\rm r}, \tau)=\rho({\cal E}, {\rm r}, \tau)-\frac{\rho({\cal E}, {\rm r}, \tau_{\rm p}({\cal E}, {\rm r}))}{\tau_{\rm p}({\cal E}, {\rm r})}\tau\,.\end{align}

\noindent
The transformation~\eqref{transf} turns to satisfy also the following ``half--parity'' symmetry:
 \begin{align}\label{full period}{\cal P}_{1/2}:\quad \left\{\begin{array}{l}\displaystyle \widetilde{\rm G}({\cal E}, {\rm r}, \tau)=\widetilde{\rm G}({\cal E}, {\rm r}, -\tau)\\\\
\displaystyle  \widetilde{\rm g}({\cal E}, {\rm r}, \tau)=
2\pi-\widetilde{\rm g}({\cal E}, {\rm r}, -\tau)\\\\
\displaystyle  \rho({\cal E}, {\rm r}, \tau)=-\rho({\cal E}, {\rm r}, -\tau)\end{array}\right.\qquad \forall\ -\tau_{\rm p}<\tau<\tau_{\rm p}\, . \end{align}

\noindent
 In addition, when  $-{\rm r}<{\cal E}<{\rm r}$, one has the following  ``quarter--parity''
 \begin{align}\label{half domain}{\cal P}_{1/4}:\quad \left\{\begin{array}{l}\displaystyle  \widetilde{\rm G}({\cal E}, {\rm r}, \tau)=-{\rm G}\left({\cal E}, {\rm r}, \tau_{\rm p}-\tau\right)\\\\
\displaystyle \widetilde{\rm g}({\cal E}, {\rm r}, \tau)=\widetilde{\rm g}\left({\cal E}, {\rm r}, \tau_{\rm p}-\tau\right)\\\\
\displaystyle  \rho({\cal E}, {\rm r}, \tau)=\rho\left({\cal E}, {\rm r}, \tau_{\rm p}\right)-\rho\left({\cal E}, {\rm r}, \tau_{\rm p}-\tau\right)\end{array}\right.\quad \forall\ 0\le \tau\le \tau_{\rm p}\, . \end{align}

\noindent
The change~\eqref{transf} will be constructed using, as  generating function, a solution of the Hamilton--Jacobi equation
\begin{align}\label{HJ}{\rm E}({\rm r}, {\rm G}, \partial_{\rm G} S_{\rm et})={\rm G}^2+{\rm r}\sqrt{1-{\rm G}^2}\cos\left(\partial_{\rm G} S_{\rm et}\right)={\cal E}\,.\end{align}
We choose the solution
\begin{align}S_{\rm et}^+({\cal R}, {\cal E},{\rm r}, {\rm G})=\left\{\begin{array}{l}\displaystyle \pi\sqrt{\alpha_+({\cal E}, {\rm r})}-\int_{{\rm G}}^{\sqrt{\alpha_+({\cal E}, {\rm r})}} \cos^{-1}\frac{{\cal E}-\Gamma ^2}{{\rm r}\sqrt{1-\Gamma ^2}}d\Gamma +{\cal R}{\rm r}\quad -{\rm r}\le{\cal E}<1\\\\
\displaystyle \pi-\int_{{\rm G}}^{\sqrt{\alpha_+({\cal E}, {\rm r})}} \cos^{-1}\frac{{\cal E}-\Gamma ^2}{{\rm r}\sqrt{1-\Gamma ^2}}d\Gamma +{\cal R}{\rm r}\quad\quad 1\le {\cal E}\le  1+\frac{{\rm r}^2}{4}
\displaystyle 
\end{array}\right.
\end{align}
where we denote as 
\begin{align}\label{alphapm}
\alpha_\pm({\cal E}, {\rm r})={\cal E}-\frac{{\rm r}^2}{2}\pm{\rm r}\sqrt{1+\frac{{\rm r}^2}{4}-{\cal E}}
\end{align}
 the real roots 
of
\begin{align}\label{equation}
x^2-2\left({\cal E}-\frac{{\rm r}^2}{2}\right)x+{\cal E}^2-{\rm r}^2=0\end{align} 
Note that  the equation  in~\eqref{equation} has always a positive real root all ${\rm r}$, ${\cal E}$ as in~\eqref{range}, so $\alpha_+({\cal E}, {\rm r})$ is positive.
$S_{\rm et}^+$ generates the following equations  

\begin{align}\label{ftau}\left\{\begin{array}{l}
\displaystyle  {\rm g}=-\cos^{-1}\frac{{\cal E}-{\rm G}^2}{{\rm r}\sqrt{1-{\rm G}^2}}\\\\
\displaystyle  \tau=+\int_{\widetilde{\rm G}({\cal E}, {\rm r}, \tau)}^{\sqrt{\alpha_+({\cal E}, {\rm r})}}\frac{d\Gamma }{\sqrt{(\Gamma ^2-\alpha_-({\cal E}, {\rm r}))(\alpha_+({\cal E}, {\rm r})-\Gamma ^2)}}\\\\
\displaystyle  {\rm R}={\cal R}-\frac{1}{{\rm r}}\int_{\widetilde{\rm G}({\cal E}, {\rm r}, \tau)}^{\sqrt{\alpha_+({\cal E}, {\rm r})}}\frac{({\cal E}-\Gamma ^2)d\Gamma }{\sqrt{(\Gamma ^2-\alpha_-({\cal E}, {\rm r}))(\alpha_+({\cal E}, {\rm r})-\Gamma ^2)}}=:{\cal R}+\rho({\cal E}, {\rm r}, \tau)\\\\
\displaystyle  {\rm r}={\rm r}\end{array}\right.\end{align}
The equations for ${\rm g}$ and ${\rm r}$ are immediate. We check the equation for $\tau$. Letting, for short,  $\sigma({\cal E}, {\rm r}):=\sqrt{\alpha_+({\cal E}, {\rm r})}$, we have

\begin{align}\label{Ederivative}
\tau&=\partial_{\cal E} S^+_{\rm et}({\cal R}, {\cal E}, {\rm r}, {\rm G})\nonumber\\
&=\left\{
\begin{array}{llll}
\displaystyle \pi\partial_{\cal E}\sigma({\cal E}, {\rm r})-{\partial_{\cal E}\sigma({\cal E}, {\rm r})}{\rm g}_+({\cal E}, {\rm r}) -\int_{\rm G}^{{\sigma({\cal E}, {\rm r})} }\partial_{\cal E}\cos^{-1}\frac{{\cal E}-\Gamma ^2}{{\rm r}\sqrt{1-\Gamma ^2}}d\Gamma \ &-{\rm r}\le {\cal E}< 1\\\\
\displaystyle  -{\partial_{\cal E}\sigma({\cal E}, {\rm r})}{\rm g}_+({\cal E}, {\rm r})-\int_{\rm G}^{{\sigma({\cal E}, {\rm r})} }\partial_{\cal E}\cos^{-1}\frac{{\cal E}-\Gamma ^2}{{\rm r}\sqrt{1-\Gamma ^2}}d\Gamma  &1\le {\cal E}\le 1+\frac{{\rm r}^2}{4}
\end{array}
\right.\nonumber\\
&=-\int_{\rm G}^{{\sigma({\cal E}, {\rm r})} }\partial_{\cal E}\cos^{-1}\frac{{\cal E}-\Gamma ^2}{{\rm r}\sqrt{1-\Gamma ^2}}d\Gamma \nonumber\\
&=
\int_{\widetilde{\rm G}({\cal E}, {\rm r}, \tau)}^{\sqrt{\alpha_+({\cal E}, {\rm r})}}\frac{d\Gamma }{\sqrt{(\Gamma ^2-\alpha_-({\cal E}, {\rm r}))(\alpha_+({\cal E}, {\rm r})-\Gamma ^2)}}
\end{align}
having let ${\rm g}_+({\cal E}, {\rm r}):=\cos^{-1}\frac{{\cal E}-\sigma({\cal E}, {\rm r})^2}{{\rm r}\sqrt{1-\sigma({\cal E}, {\rm r})^2}}$ and used, by~\eqref{alphapm},
\begin{align}\label{gplus}
{\rm g}_+({\cal E}, {\rm r})&=\cos^{-1}{\, \rm sign\, }\left(\frac{{\rm r}}{2}-\sqrt{1+\frac{{\rm r}^2}{4}-{\cal E}}\right)=\left\{\begin{array}{l}\displaystyle  \pi\quad -{\rm r}\le {\cal E}< 1\\\\
\displaystyle  0\quad 1\le {\cal E}\le 1+\frac{{\rm r}^2}{4}\end{array}\right.
\end{align}
Observe that $({\rm g}_+$, $\sigma)$ are the coordinates of the point where ${\rm E}$ reaches its maximum on each level set (Figure~\ref{figure1}).
The equation for ${\rm R}$ is analogous.

\noindent
Equations~\eqref{ftau} define the segment of the transformation~\eqref{transf} with $0\le  \tau\le \tau_{\rm p}$, where
\begin{align}\label{T}\tau_{\rm p}({\cal E}, {\rm r}):=\int_{\beta({\cal E}, {\rm r})}^{\sqrt{\alpha_+({\cal E}, {\rm r})}}\frac{d\Gamma }{\sqrt{(\Gamma ^2-\alpha_-({\cal E}, {\rm r}))(\alpha_+({\cal E}, {\rm r})-\Gamma ^2)}}\end{align}
is the half--period,
with
\begin{align}\label{beta}\beta({\cal E}, {\rm r})=\left\{\begin{array}{l}-\sqrt{\alpha_+({\cal E}, {\rm r})}\quad {\rm if}\quad \alpha_-({\cal E}, {\rm r})<0\\\\
\phantom{-}\sqrt{\alpha_-({\cal E}, {\rm r})}\quad {\rm if}\quad \alpha_-({\cal E}, {\rm r})>0\, . 
\end{array}\right.\end{align}
 The transformation is  prolonged to $-\tau_{\rm p}<\tau<0$  choosing the solution \begin{align}S_{\rm et}^-:=-2\pi{\rm G}-S_{\rm et}^+\end{align} of~\eqref{HJ}. It can be checked that this choice provides  the symmetry relation described in~\eqref{full period}.
Considering next the functions $S_k^\pm=S_{\rm et}^\pm+2 k \Sigma({\cal E}, {\rm r})$, where $\Sigma$ solves\footnote{The existence of the function $\Sigma({\cal E}, {\rm r})$ follows from the arguments of the next section: compare the formula in~\eqref{A1A2}.}
\begin{align}\partial_{\cal E}\Sigma=
 \tau_{\rm p}({\cal E}, {\rm r})\, ,\quad \partial_{\rm r}\Sigma=
 \rho({\cal E}, {\rm r}, \tau_{\rm p}({\cal E}, {\rm r}))\end{align}
 one obtains the extension of the transformation to $\tau\in {\mathbb R}$ verifying~\eqref{full real}.
 
 \noindent
 Observe that quarter period symmetry~\eqref{transf}, holding in the case $-{\rm r}<{\cal E}<{\rm r}$,  is an immediate consequence of the definitions~\eqref{ftau}.

 \noindent
The coordinates $({\cal R}, {\cal E}, {\rm r}, \tau)$ are referred to as {\it energy--time coordinates}.

\noindent  
The regularity of the functions $\widetilde{\rm G}({\cal E}, {\rm r}, \tau)$, $\widetilde\rho({\cal E}, {\rm r}, \tau)$, ${\cal B}({\cal E}, {\rm r})$ and $\tau_{\rm p}({\cal E}, {\rm r})$, which are relevant for the paper,
are studied in detail in Section~\ref{Regularity of the energy--time coordinates}. Their holomorphy  is not discussed.
 
 \subsection{Action--angle coordinates}
We look at the   transformation 
\begin{align}\phi_{\rm aa}:\qquad ({\cal R}_*, A_*, {\rm r}_*, \varphi_*)\to ({\cal R}, {\cal E}, {\rm r}, \tau)\end{align}
defined by equations
\begin{align}\label{actionangle}\left\{\begin{array}{l}\displaystyle  A_*={\cal A}({\cal E}, {\rm r})\\\\
\displaystyle  \varphi_*=\pi\frac{\tau}{\tau_{\rm p}({\cal E}, {\rm r})}\\\\
\displaystyle  {\rm r}_*={\rm r}\\\\
\displaystyle  {\cal R}_*={\cal R}+{\cal B}({\cal E}, {\rm r})\tau
\end{array}\right.
\end{align}
with ${\cal B}({\cal E}, {\rm r})$  as in~\eqref{B}, $\tau_{\rm p}({\cal E}, {\rm r})$ as in~\eqref{T} and ${\cal A}({\cal E}, {\rm r})$ the ``action function'', defined as

\begin{align}{\cal A}({\cal E}, {\rm r}):=\left\{\begin{array}{l}\displaystyle  \sqrt{\alpha_+({\cal E}, {\rm r})}-\frac{1}{\pi}\int_{\beta({\cal E}, {\rm r})}^{\sqrt{\alpha_+({\cal E}, {\rm r})}}\cos^{-1}\frac{{\cal E}-\Gamma ^2}{{\rm r}\sqrt{1-\Gamma ^2}}d\Gamma \qquad -{\rm r}\le {\cal E}\le1\\\\
\displaystyle  1-\frac{1}{\pi}\int_{\beta({\cal E}, {\rm r})}^{\sqrt{\alpha_+({\cal E}, {\rm r})}}\cos^{-1}\frac{{\cal E}-\Gamma ^2}{{\rm r}\sqrt{1-\Gamma ^2}}d\Gamma \qquad 1< {\cal E}\le1+\frac{{\rm r}^2}{4}
\end{array}\right.\end{align}
with $\alpha_+({\cal E}, {\rm r})$ and $\beta({\cal E}, {\rm r})$ being defined in~\eqref{alphapm},~\eqref{beta}. \\
Geometrically, ${\cal A}({\cal E}, {\rm r})$ represents the area of the region encircled by the level curves of ${\rm E}$ in Figure~\ref{figure1} in the former case, the area of its complement in the second case, divided by $2\pi$.
\\
The canonical character of the transformation~\eqref{actionangle} is recognised looking at the generating function
\begin{align}\label{S}S_{\rm aa}({\cal R}, {\cal E}, {\rm r}_*, \varphi_*)= \varphi_*{\cal A}({\cal E}, {\rm r}_*)+{\cal R}{\rm r}_*\end{align}
and using the following relations (compare the formulae in~\eqref{ftau} and~\eqref{T})
\begin{align}\label{A1A2}{\cal A}_{{\rm r}}({\cal E}, {\rm r})&=-\frac{1}{\pi{\rm r}}\int_{\beta({\cal E}, {\rm r})}^{\sqrt{\alpha_+({\cal E}, {\rm r})}}\frac{({\cal E}-\Gamma ^2)d\Gamma }{\sqrt{(\Gamma ^2-\alpha_-({\cal E}, {\rm r}))(\alpha_+({\cal E}, {\rm r})-\Gamma ^2)}}\nonumber\\
&=\frac{1}{\pi}\rho({\cal E}, {\rm r}, \tau_{\rm p})\nonumber\\
{\cal A}_{{\cal E}}({\cal E}, {\rm r})&=\frac{1}{\pi}\int_{\beta({\cal E}, {\rm r})}^{\sqrt{\alpha_+({\cal E}, {\rm r})}}\frac{d\Gamma }{\sqrt{(\Gamma ^2-\alpha_-({\cal E}, {\rm r}))(\alpha_+({\cal E}, {\rm r})-\Gamma ^2)}}\nonumber\\
&=\frac{1}{\pi}\tau_{\rm p}({\cal E}, {\rm r})
\end{align}
which allow us to rewrite~\eqref{actionangle}
as the transformation generated by~\eqref{S}:
 \begin{align}\label{generating equations}\left\{\begin{array}{l}\displaystyle  A_*={\cal A}({\cal E}, {\rm r})\\\\
\displaystyle  \varphi_*=\frac{\tau}{{\cal A}_{{\cal E}}({\cal E}, {\rm r})}\\\\
\displaystyle  {\rm r}_*={\rm r}\\\\
\displaystyle  {\cal R}_*={\cal R}+\frac{{\cal A}_{{\rm r}}({\cal E}, {\rm r})}{{\cal A}_{{\cal E}}({\cal E}, {\rm r})}\tau\,.
\end{array}\right.
 \end{align}
 The coordinates $({\cal R}_*, A_*, {\rm r}_*, \varphi_*)$ are referred to as {\it action--angle coordinates}.
 
\begin{remark}\rm
We conclude this section observing a non--negligible advantage while using {\it action--angle coordinates} 
compared to {\it energy--time} -- besides the obvious one of dealing with a constant period.
It is the law that relates ${\rm R}$ to ${\cal R}_*$, which is (see~\eqref{transf},~\eqref{split} and~\eqref{actionangle})
\begin{align}\label{good relation}{\rm R}={\cal R}_*+\rho_*(A_*, {\rm r}_*, \varphi_*)\, ,\quad {\rm with}\quad \rho_*(A_*, {\rm r}_*, \varphi_*):=\widetilde \rho\circ \phi_{\rm aa}(A_*, {\rm r}_*, \varphi_*)\end{align} where $\widetilde\rho$ is as in~\eqref{split}. Here $\rho_*(A_*, {\rm r}_*, \varphi_*)$ is a {\it periodic function} because so is the function  $\widetilde\rho$. This benefit is evident comparing with the corresponding formula with {\it energy--time} coordinates:
\begin{align}{\rm R}={\cal R}+{\cal B}({\cal E}, {\rm r})\tau+\widetilde\rho({\cal E}, {\rm r}, \tau)\end{align}
which would include the uncomfortable linear term ${\cal B}({\cal E}, {\rm r})\tau$. Incidentally, such term  would unnecessarily complicate the computations we are going to present in the next Section~\ref{Proof of Theorem}. 
\end{remark}

\subsection{Regularising coordinates}\label{Regularising coordinates}
In this section we define the the {\it regularising coordinates}. First of all we rewrite ${\cal S}_0({\rm r})$ in~\eqref{S0}  in terms of $(A_*, \varphi_*)$:
\begin{align}{\cal S}_0({\rm r}_*)=\Big\{(A_*,  \varphi_*):\quad A_*={\cal A}_{\rm s}({\rm r}_*)\, ,\ \varphi_*\in {\mathbb R}\Big\}\qquad 0<{\rm r}_*<2\end{align}
with ${\cal A}_{\rm s}({\rm r}_*)$ being the limiting value of ${\cal A}({\cal E}, {\rm r}_*)$ when ${\cal E}={\rm r}_*$:
\begin{align}{\cal A}_{\rm s}({\rm r}_*)=\left\{\begin{array}{lll}\displaystyle \sqrt{{\rm r}_*(2-{\rm r}_*)}-\frac{1}{\pi}\int_{0}^{\sqrt{{\rm r}_*(2-{\rm r}_*)}}\cos^{-1}\frac{{\rm r}_*-\Gamma ^2}{{\rm r}_*\sqrt{1-\Gamma ^2}}d\Gamma \quad &0<{\rm r}_*<1\\\\
\displaystyle  1-\frac{1}{\pi}\int_{0}^{\sqrt{{\rm r}_*(2-{\rm r}_*)}}\cos^{-1}\frac{{\rm r}_*-\Gamma ^2}{{\rm r}_*\sqrt{1-\Gamma ^2}}d\Gamma  &1<{\rm r}_*<2
\end{array}
\right.\end{align}

\noindent
We observe that
the function ${\cal A}_{\rm s}({\rm r}_*)$ is continuous in $[0, 2]$ (in particular, ${\cal A}_{\rm s}(1^-)={\cal A}_{\rm s}(1^+)$), with
\begin{align}{\cal A}_{\rm s}(0)=0\, ,\quad  {\cal A}_{\rm s}(2)=1\end{align}
and increases smoothly between those two values, as it results from the analysis of its derivative. Indeed, letting, for short, $\sigma_0({\rm r}_*):=\sqrt{{\rm r}_*(2-{\rm r}_*)}$ and proceeding  analogously as~\eqref{Ederivative}, we get
\begin{align}\label{As}
{\cal A}_{\rm s}'({\rm r}_*)
%
&=-\frac{1}{\pi}\int_{0}^{{\sigma_0({\rm r}_*)} }\partial_{{\rm r}_*}\cos^{-1}\frac{{\rm r}_*-\Gamma ^2}{{\rm r}_*\sqrt{1-\Gamma ^2}}d\Gamma \nonumber\\
&=\frac{1}{\pi{\rm r}_*}\int_{0}^{{\sigma_0({\rm r}_*)} }\frac{\Gamma  d\Gamma }{\sqrt{\sigma_0({\rm r}_*)^2-\Gamma ^2}}\nonumber\\
&=\frac{1}{\pi}\sqrt{\frac{2-{\rm r}_*}{{\rm r}_*}}\qquad \forall\ 0< {\rm r}_*< 2
\end{align}
We denote as $A_*\to{\rm r}_{\rm s}(A_*)$  the inverse function \begin{align}\label{inverseA}{\rm r}_{\rm s}:={\cal A}_{\rm s}^{-1}\end{align} and
we define two different changes of coordinates
\begin{align}\phi^k_{\rm rg}:\quad (Y_k, A_k, y_k, \varphi_k)\to ({\cal R}_*, A_*, {\rm r}_*, \varphi_*)\qquad k=\pm 1\end{align}
 via the formulae
\begin{align}\label{reg}\left\{\begin{array}{l}\displaystyle {\cal R}_*=Y_k e^{ky_k }\\\\
\displaystyle  A_*=A_k\\\\
\displaystyle {\rm r}_*=-ke^{-ky_k }+{\rm r}_{\rm s}(A_k)\\\\
\displaystyle \varphi_*=\varphi_k+Y_k e^{ky}{\rm r}_{\rm s}'(A_k)\end{array}\right.\end{align}
The transformations~\eqref{reg} are canonical, being generated by
\begin{align}
S^k_{\rm rg}(Y_k, A_k, {\rm r}_*, \varphi_*):=-\frac{Y_k}{k}\log\left|\frac{{\rm r}_{\rm s}(A_k)-{\rm r}_*}{k}\right|+A_k\varphi_*\, . 
\end{align}
The coordinates $(Y_k, A_k, y_k, \varphi_k)$ with $k=\pm 1$ are called {\it regularising coordinates}.

\section{A deeper insight into energy--time coordinates}
\label{Regularity of the energy--time coordinates}
In this section we study the functions $\widetilde{\rm G}({\cal E}, {\rm r}, \tau)$, $\widetilde\rho({\cal E}, {\rm r}, \tau)$, ${\cal B}({\cal E}, {\rm r})$ and $\tau_{\rm p}({\cal E}, {\rm r})$, described in Section~\ref{Energy--time coordinates}. We prove that $\widetilde{\rm G}({\cal E}, {\rm r}, \tau)$, $\widetilde\rho({\cal E}, {\rm r}, \tau)$ are $C^\infty$ provided that $({\cal E}, {\rm r})$ vary in a compact subset set of
\eqref{range} and we study the behaviour of ${\cal B}({\cal E}, {\rm r})$ and $\tau_{\rm p}({\cal E}, {\rm r})$ closely to ${\cal S}_0({\rm r})$.

\noindent
It reveals to be useful to perform this study via suitable other functions $\breve{\rm G}({\kappa}, \theta)$, $\breve\rho({\kappa}, \theta)$, ${\cal A}({\kappa})$ and $T_0({\kappa})$, which we now define.
We rewrite 
\begin{align}\label{changeofcoord}\widetilde{\rm G}({\cal E}, {\rm r}, \tau)=\sigma({\cal E}, {\rm r})\breve{\rm G}\big({\kappa}({\cal E}, {\rm r}), \theta({\cal E}, {\rm r}, \tau)\big)\ ,\quad \tau_{\rm p}({\cal E}, {\rm r})=\frac{T_{\rm p}\big({\kappa}({\cal E}, {\rm r})\big)}{\sigma({\cal E}, {\rm r})}\end{align}
and
\begin{align}\label{changeofcoord1}
\rho({\cal E}, {\rm r}, \tau)=-\frac{{\cal E}\tau}{{\rm r}}+\frac{\sigma({\cal E}, {\rm r})}{{\rm r}}\widehat\rho({\kappa}({\cal E}, {\rm r}), \theta({\cal E}, {\rm r}, \tau))\qquad 0\le \theta\le T_{\rm p}({\kappa})
\end{align}
where (changing, in the integrals in~\eqref{ftau}, the integration variable   $\Gamma ={\sigma}\xi$)
$\breve{\rm G}({\kappa}, \theta)$ is the unique solution of
\begin{align}\label{f(k, theta)}\int^{1}_{\breve{\rm G}({\kappa}, \theta)}\frac{d\xi}{{\sqrt{(1-\xi^2)(\xi^2-{\kappa})}}}=\theta\ ,\qquad 0\le \theta\le T_{\rm p}({\kappa})\end{align}
\begin{align}\label{hatrho}
&\widehat\rho({\kappa}, \theta)=\int^{1}_{\breve{\rm G}({\kappa}, \theta)}\frac{\xi^2 d\xi}{{\sqrt{(1-\xi^2)(\xi^2-{\kappa})}}}\qquad 0\le \theta\le T_{\rm p}({\kappa})\end{align}
and
\begin{align}\label{tmax} T_{\rm p}({\kappa})=\left\{\begin{array}{l}\displaystyle  T_{0}({\kappa})\quad 0<{\kappa}<1\\\\
\displaystyle  2{T_{0}({\kappa})}\quad {\kappa}<0 
\end{array}\right.\end{align}
with
\begin{align}\label{f0T0}T_0({\kappa}):=\int_{{\rm G}_0({\kappa})}^1\frac{d\xi}{{\sqrt{(1-\xi^2)(\xi^2-{\kappa})}}}\ ,\quad{\rm where}\quad {\rm G}_0({\kappa}):=\left\{\begin{array}{l}\displaystyle  {\sqrt{\kappa}}\quad 0<{\kappa}<1\\\\
\displaystyle  0\quad {\kappa}<0
\end{array}\right.\end{align}
The function $\widehat\rho({\kappa}, \theta)$ in~\eqref{hatrho} is further split as 
\begin{align}\label{hatrhosplit}\widehat\rho({\kappa}, \theta)={\cal A}({\kappa})\theta+\breve\rho({\kappa}, \theta)\end{align}
where
\begin{align}\label{A(k)NEW}{\cal A}({\kappa})=\frac{\widehat\rho({\kappa}, T_{\rm p}({\kappa}))}{T_{\rm p}({\kappa})}\ ,\quad \breve\rho({\kappa}, \theta)=\widehat\rho(\kappa, \theta)-{\cal A}(\kappa)\theta\,.\end{align}
 Finally, $\sigma({\cal E}, {\rm r})$, $\kappa({\cal E}, {\rm r})$ and $\theta({\cal E}, {\rm r}, \tau)$ are given by
\begin{align}\label{sigmakappa}
&\sigma({\cal E}, {\rm r}):=\sqrt{\alpha_+({\cal E}, {\rm r})}=\sqrt{{\cal E}-\frac{{\rm r}^2}{2}+{\rm r}\sqrt{1+\frac{{\rm r}^2}{4}-{\cal E}}}\nonumber\\
&  \kappa({\cal E}, {\rm r}):=\frac{\alpha_-({\cal E}, {\rm r})}{\alpha_+({\cal E}, {\rm r})}=\frac{{\cal E}^2-{\rm r}^2}{\left({\cal E}-\frac{{\rm r}^2}{2}+{\rm r}\sqrt{1+\frac{{\rm r}^2}{4}-{\cal E}}\right)^2}\nonumber\\
&\theta({\cal E}, {\rm r}, \tau):=\tau\sqrt{{\cal E}-\frac{{\rm r}^2}{2}+{\rm r}\sqrt{1+\frac{{\rm r}^2}{4}-{\cal E}}}\,.
\end{align}
The periodicity of $\breve\rho(\kappa, \cdot)$ (see equation~\eqref{extendhatrho} below), the uniqueness of the splitting~\eqref{split} and the formulae in~\eqref{changeofcoord1} and~\eqref{hatrhosplit} imply that ${\cal A}(\kappa)$ and $\breve\rho(\kappa, \theta)$ are related to ${\cal B}({\cal E}, {\rm r})$ and   $\widetilde\rho({\cal E}, {\rm r}, \tau)$  in~\eqref{split} via
\begin{align}\label{tilderhoB}{\cal B}({\cal E}, {\rm r})=-\frac{{\cal E}}{{\rm r}}+\frac{\sigma({\cal E}, {\rm r})^2}{{\rm r}}{\cal A}(\kappa)\ ,\quad \widetilde\rho({\cal E}, {\rm r}, \tau)=\frac{\sigma({\cal E}, {\rm r})}{{\rm r}}\breve\rho(\kappa({\cal E}, {\rm r}), \theta({\cal E}, {\rm r}, \tau))\,.\end{align}

\vskip.2in
\noindent
In view of relations~\eqref{changeofcoord},~\eqref{tmax} and~\eqref{tilderhoB}, we focus on the functions $\breve{\rm G}(\kappa, \theta)$, $\breve\rho(\kappa, \theta)$, ${\cal A}(\kappa)$ and $T_0(\kappa)$.
The  proofs of the following statements are postponed at the end of the section.

\noindent
Let us denote
$\breve{\rm G}_{ij}(\kappa, \theta):=\partial^{i+j}_{\kappa^i\theta^j}\breve{\rm G}(\kappa, \theta)$, $\breve\rho_{ij}(\kappa, \theta):=\partial^{i+j}_{\kappa^i\theta^j}\breve\rho(\kappa, \theta)$.

\begin{proposition}\label{Gderiv}
 Let $0\ne\kappa<1$ fixed. The functions  $\breve{\rm G}_{ij}(\kappa, \cdot)$ and $\breve\rho_{ij}(\kappa, \cdot)$  
are continuous for all $\theta\in {\mathbb R}$. 
\end{proposition}
This immediately implies
\begin{corollary}
Let ${\mathbb K}\subset{\mathbb R}$ a compact set, with $0$, $1\notin{\mathbb K}$. Then $\breve{\rm G}$, $\breve\rho$ are $C^\infty({\mathbb K}\times{\mathbb T})$.
\end{corollary}

\noindent
Concerning $T_0(\kappa)$, we have

\begin{proposition}\label{real period}
Let $0\ne \kappa<1$, and let  $T_0(\kappa)$ be as in~\eqref{f0T0}.
Then one can find two real numbers $C^*$, ${\cal R}^*$, ${\cal S}^*$ and two  functions ${\cal R}(\kappa)$, ${\cal S}(\kappa)$ verifying
\begin{align}{\cal R}(0)=1={\cal S}(0)\, ,\quad 0\le {\cal R}(\kappa)\le {\cal R}^*\, ,\quad 0\le {\cal S}(\kappa)\le {\cal S}^*\qquad \forall\ \kappa\in (-1, 1)\end{align}
such that  
\begin{align}\label{T0preal}T_0'(\kappa)={-\frac{{\cal R}(\kappa)}{2\kappa}}\, ,\quad  T_0''(\kappa)={\frac{{\cal S}(\kappa)}{4\kappa^2}}\, , \qquad \forall\ 0\ne \kappa<1\end{align}
In particular,
\begin{align}\label{T0}|T_0(\kappa)|\le 
\frac{{\cal R}^*}{2}\Big|\log|\kappa|\Big|+C^*\,,\quad|T'_0(\kappa)|\le 
\frac{{\cal R}^*}{2}\Big|\kappa\Big|^{-1}\,,\quad|T''_0(\kappa)|\le 
\frac{{\cal S}^*}{4}\Big|\kappa\Big|^{-2}\, .  \end{align} 
\end{proposition}

\noindent
Finally, as for ${\cal A}(\kappa)$, we have

\begin{proposition}\label{calA} Let $0\ne \kappa<1$, and let  ${\cal A}(\kappa)$ be as in 
\eqref{A(k)NEW}.
Then one can find  $C^*>0$
such that  
\begin{align}|{\cal A}(\kappa)|\le 
C^*\Big|\log|\kappa|\Big|^{-1}\,,\quad|{\cal A}'(\kappa)|\le 
C^*\Big|\kappa\Big|^{-1}\,,\quad|{\cal A}''(\kappa)|\le 
C^*\Big|\kappa\Big|^{-2} \end{align} \, . 
\end{proposition}


\paragraph{Proofs of Propositions~\ref{Gderiv},~\ref{real period} and~\ref{calA}}
 Relations 
~\eqref{full real},~\eqref{full period} and~\eqref{half domain}
provide

 \begin{align}\label{fsym}
& \left\{\begin{array}{llll}
\displaystyle  \breve{\rm G}(\kappa, \theta+2j T_{\rm p})=\breve{\rm G}(\kappa, \theta)\qquad &\forall\ \theta\in {\mathbb R}\, ,\ j\in {\mathbb Z}&\forall\ 0\ne \kappa<1\\\\
\displaystyle  \breve{\rm G}(\kappa, -\theta)=\breve{\rm G}(\kappa, \theta)&\forall\ 0\le \theta\le T_{\rm p}(\kappa)&\forall\ 0\ne \kappa<1\\\\
\displaystyle \breve{\rm G}(\kappa, T_{\rm p}-\theta)=-\breve{\rm G}(\kappa, \theta)&\forall\  0\le \theta\le T_0(\kappa)\quad &\forall\  \kappa<0.
\end{array}
\right.
 \\\nonumber\\\nonumber\\\label{extendhatrho}
&\left\{\begin{array}{llll}
\displaystyle \breve\rho(\kappa, \theta+2j T_{\rm p})=\breve\rho(\kappa, \theta)
\, ,\qquad\   &\forall\ \theta\in{\mathbb R}\, ,\ j\in {\mathbb Z}&\forall\ 0\ne \kappa<1\\\\
\displaystyle \breve\rho(\kappa, -\theta)=-\breve\rho(\kappa, \theta)\, &\forall\ 0\le \theta\le T_{\rm p}(\kappa)&\forall\ 0\ne \kappa<1\\\\
\displaystyle  \breve\rho(\kappa, T_{\rm p}-\theta)=
- \breve\rho(\kappa, \theta) & \forall\ 0\le \theta\le T_0(\kappa)  & \forall\  \kappa<0.
\end{array}
\right.
\end{align}
The following lemmata are obvious
\begin{lemma}\label{dercont}
Let $g(\kappa, \cdot)$ verify~\eqref{fsym} with $T_{\rm p}(\kappa)=\pi$ for all $\kappa$ and $T_0$ as in~\eqref{tmax}. Then the functions $g_{ij}(\kappa, \theta):=\partial^{i+j}_{\kappa^i,\theta^j} g(\kappa, \theta)$  are continuous on ${\mathbb R}$ if and only if they are continuous in $[0, T_0]$ and verify
\begin{align}\label{continuity conditions}
\left\{ \begin{array}{lllll}{\rm no\ further\ condition}\quad &{\rm if}\quad &j\in 2{\mathbb N}\,,\quad &0<\kappa<1\\\\
g_{ij}(\kappa, \frac{\pi}{2})=0&{\rm if}\quad &j\in 2{\mathbb N}\,,\quad &\kappa<0\\\\
g_{ij}(\kappa, 0)=0=g_{ij}(\kappa, \pi)&{\rm if}\quad &j\in 2{\mathbb N}+1\,,\quad &0<\kappa<1\\\\
g_{ij}(\kappa, 0)=0&{\rm if}\quad &j\in 2{\mathbb N}+1\,,\quad &\kappa<0
\end{array}
\right.
\end{align}
\end{lemma}

\begin{lemma}\label{periodicity}
Let $g(\kappa, \cdot)$ verify~\eqref{extendhatrho} with $T_{\rm p}(\kappa)=\pi$ for all $\kappa$ and $T_0$ as in~\eqref{tmax}. Then $g_{ij}(\kappa, \cdot)$, where $g_{ij}(\kappa, \theta):=\partial^{i+j}_{\kappa^i,\theta^j} g(\kappa, \theta)$,  are continuous on ${\mathbb R}$ if and only if they are continuous in $[0, T_0(\kappa)]$ and verify
\begin{align}\label{continuity condition2}
\left\{
\begin{array}{llllll}g_{ij}(\kappa, 0)= g_{ij}(\kappa, \pi)=0\quad&{\rm if}\ &j\in 2{\mathbb N}\,,\ &0<\kappa<1\\\\
g_{ij}(\kappa, 0)= g_{ij}(\kappa, \frac{\pi}{2})=0&{\rm if}&j\in 2{\mathbb N}&\kappa<0\\\\
\displaystyle  {\rm no\ further\ condition}&{\rm if}&j\in 2{\mathbb N}+1&
\end{array}
\right.
\end{align}
\end{lemma}

\subparagraph{\it Proof of Proposition~\ref{Gderiv}}
(i) The function $\breve{\rm G}(\kappa, \cdot)$ is $C^\infty({\mathbb R})$ for all $0\ne \kappa<1$~\cite{freitagB2005}. Then so is the function
$g(\kappa, \cdot)$, where
$g(\kappa, \theta):=\breve{\rm G}(\kappa, \frac{T_{\rm p}(\kappa)}{\pi}\theta )$. Then~\eqref{continuity conditions} hold true for $g(\kappa, \theta)$ with $i=0$. Hence, the derivatives $g_{ij}(\kappa, \theta)$, which exist for all $0\ne \kappa<1$, also verify~\eqref{continuity conditions}. Then $g_{ij}(\kappa, \cdot)$ are continuous for all $0\ne \kappa<1$ and so are the $\breve{\rm G}_{ij}(\kappa, \cdot)$. 

\noindent
(ii) We  check conditions~\eqref{continuity condition2} for the function $g(\kappa, \theta):=\breve\rho(\kappa, \frac{T_{\rm p}(\kappa)}{\pi}\theta)$, in the case
 $j=0$.
 Using~\eqref{hatrho},~\eqref{f(k, theta)} and~\eqref{A(k)NEW}, we get, for $0<\kappa<1$,
\begin{align}\label{id1}
g\left(\kappa, 0\right)=\breve\rho(\kappa, 0)=0\ ,\quad g\left(\kappa, \pi\right)=	{\breve\rho(\kappa, T_{\rm p}(\kappa))=\widehat\rho(\kappa, T_{\rm p})-\frac{\widehat\rho(\kappa, T_{{\rm p}})}{T_{\rm p}} T_{\rm p}=0.}
\end{align}
while, for $\kappa<0$,
\begin{align}\label{id2}
g\left(\kappa, 0\right)=\breve\rho(\kappa, 0)=0\ ,\quad g\left(\kappa, \frac{\pi}{2}\right)=	{\breve\rho(\kappa, T_0(\kappa))=\widehat\rho(\kappa, T_0)-\frac{\widehat\rho(\kappa, T_{0})}{T_0} T_0=0.}
\end{align}
The identities~\eqref{id1} and~\eqref{id2} still hold replacing $g$ with any $g_{i0}(\kappa, \theta)$,  with $i\in{\mathbb N}$, therefore,
 any $g_{i0}(\kappa, \theta)$ satisfies~\eqref{continuity condition2}.
Let us now consider the case $j\ne 0$. Again by~\eqref{hatrho},~\eqref{f(k, theta)} and~\eqref{A(k)NEW},
\begin{align}\label{breverhotheta}
\breve\rho_\theta(\kappa, \theta)&=\breve{\rm G}(\kappa, \theta)^2-{\cal A}(\kappa)
\end{align}
so, for any $j\ne 0$,
\begin{align}
\breve\rho_{ij}(\kappa, \theta)=\partial^{i+j-1}_{\kappa^i\theta^{j-1}}\Big(\breve{\rm G}(\kappa, \theta)^2\Big)
\end{align}
Then the $\breve\rho_{ij}(\kappa, \cdot)$ with $j\ne 0$ are continuous because so is $\breve{\rm G}_{ij}(\kappa, \cdot)$. $\qquad \square$

\subparagraph{\it Proof of Proposition~\ref{real period}}
The function $T_0(\kappa)$ in~\eqref{f0T0} is studied in detail in Appendix~\ref{elliptic integrals}. Combining Lemma~\ref{splitperiod} and Proposition~\ref{period} and taking the $\kappa$--primitive of such relations, one obtains Proposition~\ref{real period}.

\subparagraph{\it Proof of Proposition~\ref{calA}}
\begin{align}\label{hatrho1}
{\cal A}(\kappa)&=\frac{1}{T_0(\kappa)}\int^1_{{\rm G}_{0}(\kappa)}\frac{\sqrt{\xi^2-\kappa}}{{\sqrt{1-\xi^2}}}d\xi+\kappa 
\end{align}
\begin{align}\label{Ap}
{{\cal A}'(\kappa)}&=\frac{1}{2}+(\kappa-{\cal A}(\kappa))\frac{T'_0(\kappa)}{T_0(\kappa)}=\frac{1}{2}-(\kappa-{\cal A}(\kappa))\frac{{\cal R}(\kappa)}{2\kappa T_0(\kappa)}\nonumber\\
&=\frac{1}{2}-\frac{{\cal R}(\kappa)}{2T_0(\kappa)}+\frac{{\cal A}(\kappa){\cal R}(\kappa)}{2\kappa T_0(\kappa)}
\end{align}
and 
 \begin{align}\label{cApprime}
 	{{\cal A}''(\kappa)}&=(1-{\cal A}'(\kappa))\frac{T'_0(\kappa)}{T_0(\kappa)}+(\kappa-{\cal A}(\kappa))\left(\frac{T''_0(\kappa)}{T_0(\kappa)}-\frac{\big(T'_0(\kappa)\big)^2}{\big(T_0(\kappa))^2}\right)\nonumber\\
 	&=\frac{T'_0(\kappa)}{2T_0(\kappa)}-2(\kappa-{\cal A}(\kappa))\frac{\big(T'_0(\kappa)\big)^2}{\big(T_0(\kappa))^2}+(\kappa-{\cal A}(\kappa))\frac{T''_0(\kappa)}{T_0(\kappa)}\nonumber\\
	&=-\frac{{\cal R}(\kappa)}{4 \kappa T_0(\kappa)}-2(\kappa-{\cal A}(\kappa))\frac{{\cal R}(\kappa)^2}{4\kappa^2 T_0(\kappa)^2}+(\kappa-{\cal A}(\kappa))\frac{{\cal S}(\kappa)}{4\kappa^2 T_0(\kappa)}\qquad \square
 \end{align}

\section{The function ${\rm F}({\cal E}, {\rm r})$}\label{appB}

In this section we study the function ${\rm F}({\cal E}, {\rm r})$  in~\eqref{relation***}. Specifically, we aim to prove the following

\begin{proposition}\label{Ffunct}
${\rm F}({\cal E}, {\rm r})$  is well defined and smooth for all $({\cal E}, {\rm r})$ with $0\le {\rm r}<2$ and $-{\rm r}\le {\cal E}<1+\frac{{\rm r}^2}{4}$, ${\cal E}\ne {\rm r}$. Moreover, there exists a number $C>0$ and a neighbourhood ${\cal O}$  of $0\in {\mathbb R}$ such that, for all $0\le {\rm r} <2$ and all $-{\rm r}\le {\cal E}<1+\frac{{\rm r}^2}{4}$ such that ${\cal E}-{\rm r}\in {\cal O}$,
\begin{align}\label{FFineq}|{\rm F}({\cal E}, {\rm r})|\le C\log|{\cal E}-{\rm r}|^{-1}\,,\quad |\partial_{{\cal E}, {\rm r}}{\rm F}({\cal E}, {\rm r})|\le C|{\cal E}-{\rm r}|^{-1}\,,\quad |\partial^2_{{\cal E}, {\rm r}}{\rm F}({\cal E}, {\rm r})|\le C|{\cal E}-{\rm r}|^{-2}\, . \end{align}

\end{proposition}

\noindent
To prove Proposition~\ref{Ffunct} we need an analytic representation  of the function ${\rm F}$, which we proceed to provide.
In terms of the coordinates~\eqref{coord}, the function ${\rm U}$ in~\eqref{3bpav} is given by (recall we have fixed $\Lambda=1$)
\begin{align}\label{U}{\rm U}({\rm r}, {\rm G}, {\rm g})&=\frac{1}{2\pi}\int_0^{2\pi}\nonumber\\
&\frac{\big(1-\sqrt{1-{\rm G}^2}\cos\xi\big)d\xi}{\sqrt{
(1-\sqrt{1-{\rm G}^2}\cos\xi)^2+2{\rm r}\Big(
(\cos\xi-\sqrt{1-{\rm G}^2})\cos{\rm g}-{\rm G}\sin\xi\sin{\rm g}
\Big)+{\rm r}^2}
} \nonumber\\\end{align}
where $\xi$ is the eccentric anomaly.
By~\cite{pinzari19},  ${\rm U}$  remains constant along the level curves,  at ${\rm r}$ fixed, of the function ${\rm E}({\rm r}, \cdot, \cdot)$ in~\eqref{E}. Therefore, the function ${\rm F}({\cal E}, {\rm r})$ which realises~\eqref{relation***} is nothing else than the value that ${\rm U}({\rm r}, \cdot, \cdot)$ takes at a chosen fixed point$({\rm G}_0({\cal E}, {\rm r}), {\rm g}_0({\cal E}, {\rm r}))$  of the level set ${\cal E}$ in Figure~\ref{figure1}. For the purposes\footnote{Compare~\eqref{F(E,r)} with the simpler formula proposed in~\cite{pinzari20a}, however valid only for values of ${\cal E}$ in the interval $[-{\rm r}, {\rm r})$.} of the paper, we choose such point to be  the point where the ${\cal E}$--level curve attains its maximum.  It follows from the discussion in Section~\ref{Energy--time coordinates} that the coordinates of such point are
\begin{align}\label{G+g+}\left\{\begin{array}{l}\displaystyle {\rm G}_+({\cal E}, {\rm r})=\sqrt{\alpha_+({\cal E}, {\rm r})}\\		\\
\displaystyle {\rm g}_+({\cal E}, {\rm r})=\left\{\begin{array}{l}\displaystyle  \pi\quad -{\rm r}\le {\cal E}< 1\\\\
\displaystyle  0\quad 1\le {\cal E}\le 1+\frac{{\rm r}^2}{4}\end{array}\right.\end{array}\right.\end{align}
where $\alpha_+({\cal E}, {\rm r})$ is as in~\eqref{alphapm}. Replacing~\eqref{G+g+} into~\eqref{U}, we obtain

\begin{align}\label{F(E,r)}{\rm F}({\cal E}, {\rm r})=
\frac{1}{2\pi}\int_0^{2\pi}\frac{(1-|e({\cal E}, {\rm r})|\cos\xi)d\xi}{\sqrt{(1-|e({\cal E}, {\rm r})|\cos\xi)^2+2s({\cal E}, {\rm r}){\rm r}(\cos\xi-|e({\cal E}, {\rm r})|)+{\rm r}^2}}
\end{align}
with
\begin{align}e({\cal E}, {\rm r})=\frac{{\rm r}}{2}-\sqrt{1+\frac{{\rm r}^2}{4}-{\cal E}}\,,\quad s({\cal E}, {\rm r}):={\, \rm sign\, } \big(e({\cal E}, {\rm r})\big)=\left\{\begin{array}{l}\displaystyle  -1\quad -{\rm r}\le {\cal E}< 1\\\\
\displaystyle  +1\quad 1< {\cal E}\le 1+\frac{{\rm r}^2}{4}\end{array}\right.\end{align}
To study the regularity of ${\rm F}$, it turns to be useful
to rewrite the integral~\eqref{F(E,r)} as twice the integral on the half period $[0, \pi]$ and next to make two subsequent changes of variable. The first time, with $z= s({\cal E}, {\rm r})\cos x$. It gives the following formula, which will be used below. 
 \begin{align}\label{first formula}{\rm F}({\cal E}, {\rm r})=
\frac{1}{\pi}\int_{-1}^{1}\frac{1}{\sqrt{1-z^2}}\frac{(1-e({\cal E}, {\rm r})z)d z}{\sqrt{(1-e({\cal E}, {\rm r})z)^2+2{\rm r}(z-e({\cal E}, {\rm r}))+{\rm r}^2}}
\end{align}
We denote as
\begin{align}\label{zetapm}z_{\pm}({\cal E}, {\rm r}):=\frac{e({\cal E}, {\rm r})-{\rm r}}{e({\cal E}, {\rm r})^2}\pm\frac{\sqrt{
{\rm r}({\rm r}-2e({\cal E}, {\rm r}))(1-e({\cal E}, {\rm r})^2)
}}{e({\cal E}, {\rm r})^2}\end{align}
the roots of the polynomial under the square root, which, as we shall see below, are real under conditions~\eqref{range}. As a second change, we let $z=\frac{1-\beta^2 t^2}{1+\beta^2 t^2}$. This leads to  write ${\rm F}({\cal E}, {\rm r})$ as
\begin{align}\label{reprF}{\rm F}({\cal E}, {\rm r})&=\frac{2(1-e({\cal E}, {\rm r})) }{\pi|e({\cal E}, {\rm r})|\sqrt{(z_-({\cal E}, {\rm r})+1)(z_+({\cal E}, {\rm r})-1)}}\left(\frac{1+e({\cal E}, {\rm r})}{1-e({\cal E}, {\rm r})}j_0(\kappa({\cal E}, {\rm r}))\right.\nonumber\\
&\left.-\frac{2e({\cal E}, {\rm r})}{1-e({\cal E}, {\rm r})}j_{\beta({\cal E}, {\rm r})}(\kappa({\cal E}, {\rm r}))\right)\end{align}
where 
 $j_\beta eta(\kappa)$ is the elliptic integral
\begin{align}\label{h} j_\beta (\kappa):=\int_{0}^{+\infty}\frac{1}{1+\beta t^2}\frac{dt}{\sqrt{\big(1+ t^2\big)\big(1+\kappa t^2\big)}}
\end{align}
and $\beta$,  $\kappa$ are taken to be
\begin{align}\label{betagammakappa}\beta({\cal E}, {\rm r}):=\frac{z_-({\cal E}, {\rm r})-1}{1+z_-({\cal E}, {\rm r})}
\,,\quad  \kappa({\cal E}, {\rm r}):=\frac{(1+z_+({\cal E}, {\rm r}))(
z_-({\cal E}, {\rm r})-1)}{(1+z_-({\cal E}, {\rm r}))(z_+({\cal E}, {\rm r})-1)}\, . \end{align}
The elliptic integrals $j_\beta (\kappa)$ in~\eqref{h} are studied in Appendix~\ref{elliptic integrals}: compare Proposition~\ref{period}.

\noindent In terms of $(e, {\rm r})$, the inequalities in
~\eqref{range} become
\begin{align}\label{rhate}{\rm r}\in [0,\ 2]\ ,\quad e\in \left[-1, \ \frac{{\rm r}}{2}\right]\setminus\{0,\ {\rm r}-1\}\subset [-1, 1] \end{align}
where
$\{e=-1\}$ corresponds to the minimum level $\{{\cal E}=-{\rm r}\}$;
$\{e={\rm r}-1\}$
corresponds to the separatrix level ${\cal S}_0({\rm r})$; $\{e=0\}$
corresponds to the separatrix level ${\cal S}_1({\rm r})$ and, finally, 
$\{e=\frac{{\rm r}}{2}\}$ corresponds to maximum level $\{{\cal E}=1+\frac{{\rm r}^2}{4}\}$. It is so evident that the discriminant in~\eqref{zetapm} is not negative under conditions~\eqref{rhate}, so $z_\pm ({\cal E}, {\rm r})$ are real under~\eqref{range}, as claimed.
In addition, one can easily verify that,for any $({\rm r}, e)$ as~\eqref{rhate},
it is $e^2+e-{\rm r}\le0$. This implies
\begin{align}
\displaystyle  z_++1=\frac{e({\cal E}, {\rm r})^2+e({\cal E}, {\rm r})-{\rm r}}{e({\cal E}, {\rm r})^2}+\frac{\sqrt{
{\rm r}({\rm r}-2e({\cal E}, {\rm r}))(1-e({\cal E}, {\rm r})^2)
}}{e({\cal E}, {\rm r})^2}<0\quad\forall\ e\ne {\rm r}-1\,. 
\end{align}
Moreover, since
 \begin{align}\label{zpm0}z_-({\cal E}, {\rm r})<z_+({\cal E}, {\rm r})\quad \forall\ {\rm r}\ne 0\ ,\ {\cal E}\ne 1+\frac{{\rm r}^2}{2}\,,\ {\cal E}\ne -{\rm r}\,,\ ({\cal E}, {\rm r})\ne (2,2)\end{align}
we  have
\begin{align}\beta({\cal E},{\rm r})>0\quad \forall\ ({\cal E}, {\rm r})\ {\rm as\ in\ }\eqref{zpm0}\end{align}
and
\begin{align}0<\kappa({\cal E},{\rm r})<1\quad \forall\ ({\cal E}, {\rm r})\ {\rm as\ in\ }\eqref{zpm0}\quad {\rm and}\quad {\cal E}\ne {\rm r}-1\, . \end{align}

\noindent Combining these informations with the formula in~\eqref{reprF} and with Proposition~\ref{period}, we  conclude that ${\rm F}({\cal E}, {\rm r})$ is smooth for all ${\rm r}\ne 0\ ,\ {\cal E}\ne 1\ ,\ {\cal E}\ne 1+\frac{{\rm r}^2}{2}\,,\ {\cal E}\ne \pm{\rm r}\,,\ ({\cal E}, {\rm r})\ne (2,2)$ and that~\eqref{FFineq} holds.
However, the representation in~\eqref{first formula} allows to extend regularity for ${\rm F}({\cal E}, {\rm r})$ to the domain $0\le {\rm r}<2$,  $-{\rm r}\le {\cal E}<1+\frac{{\rm r}^2}{4}$, ${\cal E}\ne {\rm r}$, as claimed. $\qquad \square$

\section{Proof of Theorem B}\label{Proof of Theorem}
In this section we state and prove a more precise statement of Theorem B, which is Theorem~\ref{main*} below. 

\noindent
The framework is  as follows:

\begin{itemize}
\item[\tiny\textbullet]  fix a energy level $c$;
\item[\tiny\textbullet]  change the time via \begin{align}\label{tau}\frac{dt}{d{t'}}=e^{- 2k y}\qquad k=\pm 1\end{align} where $t'$ is the new time and $t$ the old one. The new time $t'$ is soon renamed $t$;
\item[\tiny\textbullet]  look at the ODE
\begin{align}\label{Xk}\partial_{t} q_k=X^{{ (k) }}(q_k; c)\end{align}
for the triple $q_k=(A_k, y_k, \psi)$ where $A_k$, $y_k$ are as in~\eqref{reg}, while $\psi=\varphi_*$, with $\varphi_*$ as in~\eqref{actionangle}  in ${\mathbb P}_k$, where
\begin{align}{\mathbb P}_k(\varepsilon_-, \varepsilon_+, L_-, L_+, \xi)&:= \Big\{(A_k, y_k, \psi):\ 1-2\varepsilon_+<A_k\le 1-2\varepsilon_-\ ,\ L_-+2\xi\le ky_k\le L_+-2\xi,\nonumber\\
&\ \psi\in {\mathbb T}\Big\}
\end{align}
with $\xi<(L_+-L_-)/4$.
Observe that
\begin{itemize}
\item[{\tiny\textbullet}] the projection of ${\mathbb P}_+$ in the plane $({\rm g}, {\rm G})$ in Figure~\ref{figure1} is an  
inner region of  ${\cal S}_0({\rm r})$ and ${\rm r}$ varies in a $\varepsilon$--left neighburhood of $2$;
\item[{\tiny\textbullet}]  the projection of ${\mathbb P}_-$ in the plane $({\rm g}, {\rm G})$ in Figure~\ref{figure1} is an 
outer region of  ${\cal S}_0({\rm r})$ and ${\rm r}$ varies in a $\varepsilon$--left neighburhood of $2$; 
\item[{\tiny\textbullet}] the boundary of ${\mathbb P}_\kappa$ includes ${\cal S}_0$  
if $L_+=\infty$; it has a positive distance from it if $L_+<+\infty$.  
\end{itemize}
\end{itemize}
We shall prove
\begin{theorem}\label{main*}
There exist a graph ${{\cal G}}_k\subset {\mathbb P}_k(\varepsilon_-, \varepsilon_+, L_-, L_+, \xi)$ and a number $L_\star>1$ such that for any $L_->L_\star$  there exist $\varepsilon_-$, $\varepsilon_+$, $L_+$, $\xi$,  an open neighbourhood $W_k\supset{{\cal G}}_k$
such that along  any orbit $q_k(t)$ such that $q_k(0)\in W_k$, 
\begin{align}|A(q_k(t))-A(q_k(0))|\le C_0 \epsilon  e^{-L_-^3}\, t\qquad \forall\ t:\ |t|< t_{\rm ex}\end{align}
where $t_{\rm ex}$ is the first $t$ such that  $q(t)\notin W_k$ and $\epsilon$ is an upper bound for $\|P_1\|_{W_k}$ (with $P_1$ being the first component of $P$).\end{theorem}

{\bf Proof\ }
For definiteness, from now on we discuss the case $k=+1$ (outer orbits). The case $k=-1$ (inner orbits) is pretty similar. We neglect to write the sub--fix ``$+1$'' everywhere. As the proof is long and technical, we divide it in paragraphs. We shall take
\begin{align}{{\cal G}}=\Big\{(A, y, \psi_\circ(A, y)),\ 1-2\varepsilon_+\le A\le 1-2\varepsilon_-\, ,\ L_-+2\xi\le  y\le L_+-2\xi\Big\}\subset {\mathbb P}\end{align}
with  $\varepsilon_-$, $\varepsilon_+$, $L_-$, $L_+$, $\psi_\circ$ to be chosen below.

\paragraph{\it Step 1. The vector--field $X$}
As $\psi$ is one of the {\it action--angle coordinates}, while $A$, $y$ are two among the {\it regularising coordinates},  we need the expressions of the Hamiltonian~\eqref{3bpav} written in terms of those two sets.
The Hamiltonian~\eqref{3bpav} written in {\it action--angle coordinates} is
\begin{align}\label{H*}{\rm H}_{\rm aa}({\cal R}_*, A_*, {\rm r}_*,  \varphi_*)=\frac{({\cal R}_*+\rho_*(A_*, {\rm r}_*,  \varphi_*))^2}{2}+\alpha{\rm F}_*(A_*, {\rm r}_*)+\frac{({\rm C}-{\rm G}_*(A_*, {\rm r}_*,  \varphi_*))^2}{2{\rm r}_*^2}-\frac{\beta}{{\rm r}_*}\end{align}
where
\begin{align}\label{G*r*}{\rm G}_*(A_*, {\rm r}_*, \varphi_*):={\rm G}\circ \phi_{\rm aa}(A_*, {\rm r}_*, \varphi_*)\, ,\qquad {\rm F}_*(A_*, {\rm r}_*):={\rm F}\circ \phi_{\rm aa}(A_*, {\rm r}_*)
\end{align}  with $\phi_{\rm aa}$ as in~\eqref{actionangle}, while $\widetilde{\rm G}({\cal E}, {\rm r}, \tau)$,  ${\rm F}({\cal E}, {\rm r})$ as in~\eqref{transf},~\eqref{relation***}, respectively, $\rho_*$ is as in~\eqref{good relation}. The Hamiltonian~\eqref{3bpav} written in {\it regularising coordinates} is
\begin{align}{\rm H}_{\rm rg}(Y, A, y,  \varphi)&=\frac{(Y e^{y}+\rho_*(A, {\rm r}_\circ(A, y),  \varphi_\circ(Y, A, y,  \varphi)))^2}{2}+\alpha{\rm F}_*(A, {\rm r}_\circ(A, y))\nonumber\\
&+ \frac{({\rm C}-{\rm G}_*(A, {\rm r}_\circ(A, y),  \varphi_\circ(Y, A, y,  \varphi)))^2}{2{\rm r}_\circ(A, y)^2}-\frac{\beta}{{\rm r}_\circ(A, y)}\end{align}
where
${\rm r}_\circ(A, y)$, $\varphi_\circ(Y, A, y,  \varphi)$ are the right hand sides of the equations for ${\rm r}_*$, $\varphi_*$ in~\eqref{reg}, with $k=+1$.\\
Taking the $\varphi_*$--projection of Hamilton equation of ${\rm H}_{\rm aa}$, and
the $(A, y)$--projection of Hamilton equation of ${\rm H}_{\rm rg}$,
changing the time as prescribed in~\eqref{tau} and reducing the energy via
 \begin{align}
{\cal R}_*+\rho_{*}(A, {\rm r}_\circ(A, y), \psi)=Ye^{y }+\rho_{*}(A, {\rm r}_\circ(A, y), \psi)={\cal Y}(A, y, \psi; c)\nonumber\\
\end{align}
with
\begin{align}\label{Ycanc}{\cal Y}(A, y, \psi; c):=\pm\sqrt{2\left(
c-\alpha{\rm F}_*(A, {\rm r}_\circ(A, y))
-\frac{({\rm C}-{\rm G}_*(A, {\rm r}_\circ(A, y), \psi))^2}{2{\rm r}_\circ(A, y)^2}+\frac{\beta }{{\rm r}_\circ(A, y)}
\right)}\end{align}
 we find that the evolution for the triple $q=(A, y, \psi)$ during the time $t$ is governed by the vector--field
\begin{align}\label{X}\left\{\begin{array}{l}\displaystyle  X_1(A, y, \psi; c)=e^{-2y}\frac{{\rm C}-{\rm G}_*(A, {\rm r}_\circ(A, y), \psi)}{{\rm r}_\circ(A, y)^2}{\rm G}_{*, 3}(A, {\rm r}_\circ(A, y), \psi)-e^{-2y}\rho_{*, 3}(A, {\rm r}_\circ(A, y), \psi){\cal Y}(A, y, \psi; c)\\\\
\displaystyle  X_2(A, y, \psi; c)=-e^{-y}\frac{{\rm C}-{\rm G}_*(A, {\rm r}_\circ(A, y), \psi)}{{\rm r}_\circ(A, y)^2}{\rm G}_{*, 3}(A, {\rm r}_\circ(A, y), \psi){\rm r}_{\rm s}'(A)\\\ \ \qquad+
\displaystyle  e^{-y}\big(1+\rho_{*, 3}(A, {\rm r}_\circ(A, y), \psi){\rm r}_{\rm s}'(A)\big){\cal Y}(A, y, \psi; c)\\\\
\displaystyle   X_3(A, y, \psi; c)=\alpha\,e^{-2y}{\rm F}_{*, 1}(A, {\rm r}_\circ(A, y))-e^{-2y}\frac{{\rm C}-{\rm G}_*(A, {\rm r}_\circ(A, y), \psi)}{{\rm r}_\circ(A, y)^2}{\rm G}_{*, 1}(A, {\rm r}_\circ(A, y), \psi)\\\ \ \displaystyle \qquad+e^{-2y}\rho_{*, 1}(A, {\rm r}_\circ(A, y), \psi){\cal Y}(A, y, \psi; c)
\end{array}\right.\nonumber\\
\end{align}
where we have used the notation, for $f=\rho_*$, ${\rm G}_*$, ${\rm F}_*$,
\begin{align}f_{1}(A, {\rm r}_*, \psi):=\partial_{A} f(A, {\rm r}_*, \psi)\, ,\quad f_{3}(A, {\rm r}_*, \psi):=\partial_{\psi} f(A, {\rm r}_*, \psi)\, . \end{align}

\paragraph{\it Step 2. Splitting the vector--field} 
We write
\begin{align}\label{splitX}X(A, y, \psi; c)=N(A, y; c)+P(A, y, \psi; c)\end{align}
 with
\begin{align}\label{N}\left\{\begin{array}{l}\displaystyle  N_1(A, y; c)=0\\\\
\displaystyle  N_2(A, y; c)=v(A, y; c):=e^{-y }\sqrt{2\big(c-\alpha{\rm F}_*(A, {\rm r}_\circ(A, y))\big)}\\\\
N_3(A, y; c)= \omega(A, y; c):=\alpha\,e^{-2y }{\rm F}_{*, 1}(A, {\rm r}_\circ(A, y))\end{array}\right.\end{align}
hence,
\begin{align}\label{P}
\left\{\begin{array}{l}\displaystyle  P_1=e^{-2y}\frac{{\rm C}-{\rm G}_*(A, {\rm r}_\circ(A, y), \psi)}{{\rm r}_\circ(A, y)^2}{\rm G}_{*, 3}(A, {\rm r}_\circ(A, y), \psi)-e^{-2y}\rho_{*, 3}(A, {\rm r}_\circ(A, y), \psi){\cal Y}(A, y, \psi; c)\\\\
\displaystyle  P_2=-e^{-y}\frac{{\rm C}-{\rm G}_*(A, {\rm r}_\circ(A, y), \psi)}{{\rm r}_\circ(A, y)^2}{\rm G}_{*, 3}(A, {\rm r}_\circ(A, y), \psi){\rm r}_{\rm s}'(A)+
e^{-y}\rho_{*, 3}(A, {\rm r}_\circ(A, y), \psi){\rm r}_{\rm s}'(A)\nonumber\\
\qquad \displaystyle \cdot{\cal Y}(A, y, \psi; c)
+
\displaystyle  e^{-y}\left(
{\cal Y}(A, y, \psi; c)
-\sqrt{2\big(c-\alpha{\rm F}_*(A, {\rm r}_\circ(A, y))\big)}\right)\\\\
\displaystyle  P_3=-e^{-2y}\frac{{\rm C}-{\rm G}_*(A, {\rm r}_\circ(A, y), \psi)}{{\rm r}_\circ(A, y)^2}{\rm G}_{*, 1}(A, {\rm r}_\circ(A, y), \psi)+e^{-2y}\rho_{*, 1}(A, {\rm r}_\circ(A, y), \psi){\cal Y}(A, y, \psi; c)
\end{array}\right.\nonumber\\
\end{align}

\noindent
The application of {\sc nft} relies on the smallness of the perturbing term $P$. In the case in point, the ``greatest'' term of $P$ is the  component $P_2$, and precisely $\rho_{*, 3}$. This function is not uniformly small. For this reason, we need to look at its zeroes and localise around them. The localisation (described in detail below) carries the  holomorphic perturbation $P$ to a perturbation $\widetilde P$, which is smaller, but {\it no longer holomorphic}. We shall apply {\sc gnft} to the new vector--field $\widetilde X=N+\widetilde P$.

\paragraph{\it Step 3. Localisation about  non--trivial zeroes of $\rho_{*,3}$}

The following lemma gives an insight on the term $\rho_{*, 3}$, appearing in~\eqref{P}. It will be proved in Appendix~\ref{Technicalities}.

\begin{lemma}\label{zeroes}
For any 
${\cal A}_{\rm s}({\rm r}_*)<A<1$ $(0<A<{\cal A}_{\rm s}({\rm r}_*))$
 there exists  $0<\psi_*(A, {\rm r}_*))<\pi$ {\rm (}$0<\psi_*(A, {\rm r}_*))<\pi/2${\rm)} such that $\rho_{*,3}(A, {\rm r}_*, \psi_*(A, {\rm r}_*)))\equiv0$ {\rm(and} $\rho_{*,3}(A, {\rm r}_*, \pi-\psi_*(A, {\rm r}_*)))\equiv0${\rm)}. Moreover, there exists $C>0$ such that, for any  $\d>0$ one can find a neighbourhood $V_*(A, {\rm r}_*; \delta)$ of $\psi_*(A, {\rm r}_*))$ {\rm(}and a neighbourhood $V'(A, {\rm r}_*; \delta)$ of $\pi-\psi_*(A, {\rm r}_*))${\rm)} such that
\begin{align}\label{boundrho3}
&|\rho_{*,3}(A, {\rm r}_*, \psi)|\le C\frac{\sigma_*(A, {\rm r}_*)}{{\rm r}_*}{\delta}\qquad \forall\ \psi\in V_*(A, {\rm r}_*; \delta)\nonumber\\
&\left(|\rho_{*,3}(A, {\rm r}_*, \psi)|\le C\frac{\sigma_*(A, {\rm r}_*)}{{\rm r}_*}{\delta}\qquad\forall\ \psi\in V_*(A, {\rm r}_*; \delta)\cup V'(A, {\rm r}_*; \delta)\, . \right)\end{align}
\end{lemma}

\noindent
We now let
\begin{align}\psi_\circ(A, y):=\psi_*(A, {\rm r}(A, y))\ ,\quad V_\circ(A, y; \delta):=V_*(A, {\rm r}(A, y); \delta)\, . \end{align}

\noindent
For definiteness, from now on,  we focus on orbits with initial datum $(A_0, y_0, \psi_0)$ such that $\psi_0$ is close to $\psi_\circ(A_0, y_0)$. The symmetrical cases  can be similarly treated.

\noindent
Let $W_\circ(A, y; \delta)\subset V_\circ(A, y; \delta)$ an open set and let $g(A, y, \cdot)$ be a $C^\infty$, $2\pi$--periodic function such that, in each period $[\psi_\circ(A, y)-\pi, \psi_\circ(A, y)+\pi)$ satisfies
\begin{align}\label{g}g(A, y, \psi; \delta)\left\{\begin{array}{l}\displaystyle ~\equiv 1\quad \forall\ \psi\in  W_\circ(A, y; \delta)\\\\ \displaystyle \equiv   0\quad \forall\ \psi\in [\psi_\circ(A, y)-\pi, \psi_\circ(A, y)+\pi)\setminus V_\circ(A, y; \delta)\\\\
\displaystyle  \in (0, 1) \quad \forall\ \psi\in V_\circ(A, y; \delta)\setminus W_\circ(A, y; \delta)
\end{array}\right.
\end{align}
The function $g$ is chosen so that
\begin{align}\label{g0}
\sup_{0\le \ell<\ell_*}\|g\|_{u, \ell}\le 1\,.
 \end{align}
 As an example, one can take $g(A, y, \psi; \delta)=\chi(\psi-\psi_\circ(A, y))$, with
 \begin{align}
 \chi(\theta)=\left\{
 \begin{array}{lll}
 1\quad & |\theta|\le a\\
 1-\frac{\int_a^{\theta}e^{-\frac{\zeta}{(\theta-a)(b-\theta)}}d\zeta}{\int_a^{b}e^{-\frac{\zeta}{(\theta-a)(b-\theta)}}d\zeta}& a<\theta\le b\\
 0& \theta>b\\
 \chi(-\theta)& \theta<-a
 \end{array}
 \right.
 \end{align}
 with $0<a<$ $b$ so small that $B_a(\psi_\circ(A, y))\subset W_\circ(A, y; \delta)$, $B_b(\psi_\circ(A, y))\subset V_\circ(A, y; \delta)$. If $\zeta\in (0, 1)$ is sufficiently small (depending on $\ell_*$), then 
 \eqref{g0} is met.\\
Let
\begin{align}\label{tildeP}\widetilde P(A, y, \psi; \delta):= g(A, y, \psi;\delta) P(A, y, \psi)\, . \end{align}
We let
\begin{align}\widetilde X:=N+\widetilde P\end{align}
and
\begin{align}\label{domain}{\mathbb P}_{\varepsilon_-, \xi}={\mathbb A}_{\varepsilon_-}\times {\mathbb Y}_{\xi}\times {\mathbb T}\, ,\end{align}
where ${\mathbb A}=[1-2\varepsilon_+, 1-2\varepsilon_-]$, ${\mathbb Y}=[L_-+2\xi, L_+-2\xi]$ and $\varepsilon_-<\varepsilon_+$, $\xi$ are sufficiently small,
and  $u=(\varepsilon_-,  \xi)$.
By construction, $\widetilde X$ and $\widetilde P\in {\cal C}^3_{u, \infty}$. In particular, $\widetilde P\in {\cal C}^3_{u, \ell_*}$, for all $\ell_*\in{\mathbb N}$. Below, we shall fix a suitably large $\ell_*$.

\paragraph{\it Step 4. Bounds}
The following uniform bounds follow rather directly from the definitions. Their proof is deferred to Appendix~\ref{Technicalities}, in order not to interrupt the flow.
\begin{align}\label{bounds2}&
\left\|\frac{1}{v}\right\|_u \le C \frac{e^{L_+}}{\alpha L^{\frac{1}{2}}_-}\ ,\quad 
\left\|\frac{\partial_A v}{v}\right\|_u \le C \frac{e^{L_+}}{L_-\sqrt{\varepsilon_-}}\ ,\quad \left\|\frac{\partial_y v}{v}\right\|_u \le1+ C \frac{e^{L_+-L_-}}{L^2_-}\nonumber\\\nonumber\\
&\left\|\frac{\omega}{v}\right\|_u\le C\frac{e^{L_+-L_-}}{L_-^{3/2}}\ ,\quad \left\|\frac{\partial_A\omega}{v}\right\|_u\le C\frac{e^{2L_+-L_-}}{L_-^{3/2}\varepsilon_-^{\frac{1}2}}\ ,\quad \left\|\frac{\partial_y\omega}{v}\right\|_u\le C\frac{e^{2L_+-2L_-}}{L_-^{3/2}}\\\nonumber\\
\label{bounds3}
&\|\widetilde P_1\|_u\le C e^{-2 L_-}\max\Big\{
|{\rm C}|L_+\sqrt{\varepsilon_+},\ L_+\varepsilon_+,\ \delta\sqrt{\varepsilon_+}\sqrt{\alpha\,L_+}
\Big\}\nonumber\\\nonumber\\
&\|\widetilde P_2\|_u\le C e^{-L_-}\max\Big\{
|{\rm C}|L_+\sqrt{\frac{\varepsilon_+}{\varepsilon_-}},\ L_+\frac{\varepsilon_+}{\sqrt{\varepsilon_-}},\ \sqrt{\frac{\varepsilon_+}{\varepsilon_-}}\delta\sqrt{\alpha\,L_+}\, ,\  (\alpha{L_-})^{-\frac{1}{2}}\max\{|{\rm C}|^2,\ \varepsilon_+^2,\ \beta\}
\Big\}\nonumber\\\nonumber\\
&\|\widetilde P_3\|_u\le C e^{-2L_-}\max\Big\{
|{\rm C}|\frac{\sqrt{\varepsilon_+}}{\varepsilon_-},\ \frac{\varepsilon_+}{\varepsilon_-},\ \frac{\sqrt{\varepsilon_+}}{\varepsilon_-}\sqrt{\alpha\,L_+}
\Big\}
\end{align}
Here $C$  is a number not depending on $L_-$, $L_+$, $\xi$, $\varepsilon_-$, $\varepsilon_+$, $c$, $|{\rm C}|$, $\beta$, $\alpha$ and the norms are meant as in Section~\ref{A generalisation when the dependence}, in the domain~\eqref{domain}.
Remark that the validity of~\eqref{bounds3} is subject to condition 
\begin{align}\label{Lm}
L_-\ge C\alpha^{-1}\max\{|c|,\ |{\rm C}|^2,{\varepsilon_+}, \beta \}\, . 
\end{align}
which will be verified below.

\paragraph{\it Step 5. Application of {\sc gnft} and conclusion}
Fix $s_1$, $s_2>0$. Define
\begin{align}\rho:=\frac{\varepsilon_-}{16}\ ,\quad \tau:=e^{-s_2}\frac{\xi}{16}\ ,\quad w_K:=\left(\frac{\varepsilon_-}{16}\,, \frac{e^{-s_2}\xi }{16}\,, \frac{1}{c_0 K^{1+\delta}}\right)\end{align}
so that~\eqref{NEWu+positiveNEW} are satisfied.
With these choices, as a consequence of the bounds in~\eqref{bounds2}--\eqref{bounds3}, one has
\begin{align}\label{bounds10}
\chi&\le C(L_+-L_-) \max\left\{\frac{e^{L_+-L_-}}{s_1L_-^{3/2}},\ \frac{1}{s_2}\left(1+ C \frac{e^{L_+-L_-}}{L^2_-}\right)\right\}\nonumber\\
\theta_1&\le C e^{s_1}(L_+-L_-)\xi K^{1+\delta}\frac{e^{2L_+-2L_-}}{L_-^{3/2}}\nonumber\\
\theta_2&\le C e^{s_1+s_2}(L_+-L_-)\frac{\sqrt{\varepsilon_-}}{\xi}\frac{e^{L_+}}{L_-}
\nonumber\\
\theta_3&\le C{e^{s_1}}(L_+-L_-)K^{1+\delta}\sqrt{\varepsilon_-}\frac{e^{2L_+-L_-}}{L_-^{3/2}}\nonumber\\
\eta&\le C e^{s_1+s_2}(L_+-L_-)\frac{e^{L_+-L_-}}{\alpha L_-^{\frac{1}{2}}}\max\left\{
 e^{-L_-}\varepsilon_-^{-1}\max\Big\{
|{\rm C}|L_+\sqrt{\varepsilon_+},\ L_+\varepsilon_+,\ \delta\sqrt{\varepsilon_+}\sqrt{\alpha\,L_+}
\Big\}\right.\, ,\nonumber\\
&{e^{s_2}}\xi^{-1}\max\Big\{
|{\rm C}|L_+\sqrt{\frac{\varepsilon_+}{\varepsilon_-}},\ L_+\frac{\varepsilon_+}{\sqrt{\varepsilon_-}},\ \sqrt{\frac{\varepsilon_+}{\varepsilon_-}}\delta\sqrt{\alpha\,L_+}\, ,(\alpha{L_-})^{-\frac{1}{2}}\max\{|{\rm C}|^2,\ \varepsilon_+^2,\ \beta\}
\Big\}\, ,\nonumber\\
&\left. e^{-L_-}K^{1+\delta}\max\Big\{
|{\rm C}|\frac{\sqrt{\varepsilon_+}}{\varepsilon_-},\ \frac{\varepsilon_+}{\varepsilon_-},\ \frac{\sqrt{\varepsilon_+}}{\varepsilon_-}\sqrt{\alpha\,L_+}
\Big\}\right\}
\end{align}

\noindent
We now discuss  inequalities~\eqref{NEWu+positiveNEW}--\eqref{NEWnewsmallnessNEW} and~\eqref{Lm}. We choose  $s_i$, $L_\pm$,  $\varepsilon_\pm$ and $K$ to be the following functions of $L$ and $\xi$, with $0<\xi<1<L$:
\begin{align}
&L_-=L\,,\quad \varepsilon_\pm=c_\pm L^2e^{-2L}\,,\qquad L_+=L+10\xi\,,\quad s_1=C_1 \xi L^{-\frac{3}{2}}\,,\quad s_2=C_1\xi\,,\quad  K=\left[\left(\frac{c_1}{\xi\sqrt L}\right)^{\frac{1}{1+\delta}}\right]
\end{align}
with $0<c_1<1<C_1$ and $0<c_-<c_+<1$ suitably fixed, so as to have $K> 0$. A more stringent relation between $\xi$ and $L$ will be specified below.
We take
\begin{align}|{\rm C}|<c_1 L^2 e^{-2L}\,,\quad \b<c_1 L^4 e^{-4L}\,,\quad \d<c_1 L^{3/2} e^{-L}\end{align}
In view of~\eqref{bounds10}, it is immediate to check that there exist  suitable numbers $0<c_1<1<C_1$ depending only on $c$, $c_+$, $c_-$ and $\alpha$ such that
inequalities~\eqref{NEWu+positiveNEW}--\eqref{theta3NEW} and~\eqref{Lm} are satisfied and
\begin{align}\eta<C_2 L^{-\frac{3}{2}}\end{align}

\noindent
An application of {\sc gnft} conjugates $\widetilde X=N+\widetilde P$ to  a new vector--field $\widetilde X_\star=N+\widetilde P_\star$, with the first component of the vector $\widetilde P_*$ being bounded as

\begin{align}
\|\widetilde P_{\star, 1}\|_{u_\star}\le &\varepsilon_{-}
\VERT \widetilde P_\star\VERT^{w_K}_{u_\star}\le \varepsilon_{-}\max\left\{2^{-c_2L^3}\VERT \widetilde P\VERT^{w_K}_{u}\, ,\ 2 c_0\,K^{-\ell+\delta}\VERT \widetilde P\VERT^{w_K}_{u, \ell} \right\}\end{align}
Using~\eqref{g0},~\eqref{tildeP}, that $\widetilde P$ vanishes outside $V_\circ$, the chain rule and the holomorphy of $P(A, y, \cdot)$,
\begin{align}
\VERT\widetilde P\VERT^{w_K}_{u, \ell}\le 2^\ell\VERT P_{V_\circ}\VERT^{w_K}_{u, \ell}\le 2^\ell\, \frac{\ell!}{s^\ell}\,\VERT P_{(V_\circ)_s}\VERT^{w_K}_{u} \qquad \forall\ 0\le \ell\le \ell_*
\end{align}
where $P_{(V_\circ)_s}(A, y, \psi)$ denotes the restriction of $P(A, y, \cdot)$ on $(V_{\circ})_s$, while  $s$ is the analyticity radius of $P(A, y, \cdot)$.
We take $s$ so small that 
\begin{align}
\VERT  P_{(V_\circ)_s}\VERT^{w_K}_{u}\le 2 \VERT  P_{V_\circ}\VERT^{w_K}_{u}
\end{align}
Then we have
\begin{align}
\|\widetilde P_{\star, 1}\|_{u_\star}
\le &2 \varepsilon_{-}\max\left\{2^{-c_2L^3}\, ,\ c_0\,2^{\ell+1}\ell!s^{-\ell} K^{-\ell+\delta}\right\}\VERT  P_{V_\circ}\VERT^{w_K}_{u}\le 
2\varepsilon_{-}2^{-c_2L^3}\VERT  P_{(V_\circ)_s}\VERT^{w_K}_{u}\le 
2\varepsilon_{-}2^{-c_2L^3}Q^{-1}
\end{align}
where we have used the inequality
\begin{align}\label{bound4}
c_0\,2^{\ell+1}\ell!s^{-\ell} K^{-\ell+\delta}\le 2^{-c_2L^3}
\end{align}
which will be discussed below.
On the other hand, analogous techniques as the ones used to obtain~\eqref{bounds3} provide
 \begin{align}c \epsilon\le \|{P_1}_{V}\|_u\le \epsilon
 \,,\quad c L^{\frac{1}{2}} e^{-L}\le Q^{-1}\le C L^{\frac{1}{2}} e^{-L}\,.\end{align}
 with $\epsilon:=C L^3 e^{-4L}$ and $0<c<1$.
So,  
\begin{align}\|\widetilde P_{\star, 1}\|_{u_\star}
\le  C_32^{-c_3 L^3}\epsilon\end{align}
which is what we wanted to prove.
It remains to discuss~\eqref{bound4}. By Stirling and provided that $\ell>2\delta$,~\eqref{bound4} is implied by
\begin{align}
K>1\,,\quad \left(\frac{4c_0\sqrt{2\pi}\ell^{\frac{3}{2}}}{es \sqrt K}\right)^\ell\le 2^{-c_2L^3}
\end{align}
These inequalities are satisfied by choosing $\ell$, $\ell_*$ and $\xi$ to be  related to $L$ such in a way that
\begin{align}\ell&=
\max\left\{[c_2 L^3]+1\,, [2\d]+1\,,\left[\left(\frac{1}{2\pi}\frac{e^2\sigma^2}{64c_0^2}\right)^{\frac{1}{3}}\right]+1\right\}\,,\quad \ell_*>\ell\\ K&=\left[\left(\frac{c_1}{\xi\sqrt L}\right)^{\frac{1}{1+\delta}}\right]>2\pi\frac{64 c_0^2}{e^2\sigma^2}\ell^3>1\,.\quad \square\end{align}

\appendix
\section{The  elliptic integrals $T_0(\kappa)$ and $j_\beta (\kappa)$}\label{elliptic integrals}

The  functions $T_0(\kappa)$ in~\eqref{f0T0} and $j_\beta (\kappa)$ in~\eqref{h}    are complete elliptic integrals.  We use this appendix to store some useful material concerning such functions.

\noindent
First of all, in the definition of $T_0(\kappa)$, we change the integration variable, letting $\xi\to \frac{1}{\xi}$, so as to rewrite
\begin{align}\label{T0changed}
T_0(\kappa)=\int_1^{\frac{1}{{\rm G}_0(\kappa)}}\frac{d\xi}{\sqrt{(\xi^2-1)(1-\kappa\xi^2)}}\quad 0\ne \kappa<1\end{align}
with ${\rm G}_0(\kappa)$ as in~\eqref{f0T0}.
Next,
 we look at the complex--valued function

\begin{align}\label{ovlT}g(\kappa):=\int_1^{+\infty}\frac{d\xi}{\sqrt{(\xi^2-1)(1-\kappa\xi^2)}}\quad \kappa\in{\mathbb R}\setminus\{0, 1\}\end{align} 
which
is easily related to $T_0(\kappa)$ and $j_0(\kappa)$:

\begin{lemma}\label{splitperiod}
Let  $0\ne\kappa<1$. Then
\begin{align}\label{ovlTnew}
T_0(\kappa)=\left\{\begin{array}{l}\displaystyle 
g(\kappa)\quad {\rm if}\quad \kappa<0\\\\
j_0(\kappa)=\Re g(\kappa)\quad {\rm if}\quad 0<\kappa<1
\end{array}\right.
\end{align}
\end{lemma}
{\bf Proof\ }
We have only to prove that $T_0(\kappa)=j_0(\kappa)$ when $0<\kappa<1$, as the other relations are immediate, from~\eqref{T0changed} and~\eqref{ovlT}. 
 We write \begin{align}\label{split2}T_0(\kappa)=\left(\int_0^{+\infty}-\int_{0}^1-\int_{\frac{1}{\sqrt\kappa}}^{+\infty}
\right)\frac{d\xi}{\sqrt{(\xi^2-1)(1-\kappa\xi^2)}}\, . \end{align}
We deform the integration path of the first integral at right hand side stretching the real path $\xi\in[0, +\infty)$ to the purely imaginary line $z={\rm i} y$, with $y\in[0, +\infty)$, so that
\begin{align}\label{integral equality}\int_0^{\infty}\frac{d\xi}{\sqrt{(\xi^2-1)(1-\kappa\xi^2)}}=\int_0^{\infty}\frac{dy}{\sqrt{(y^2+1)(1+\kappa y^2)}}=j_0(\kappa)\end{align}
Combining this with the observation that, for $0<\kappa<1$,  $T_0(\kappa)$ and $j_0(\kappa)$ are real while the two latter integrals in~\eqref{split2} are purely imaginary, we have
$T_0(\kappa)=j_0(\kappa)$, as claimed.  $\qquad \square$

\begin{remark}\rm
It follows from the proof of Lemma~\ref{splitperiod}  (compare~\eqref{split2}--\eqref{integral equality}) that, in the sense of complex integrals,
\begin{align}\label{identity}\left(\int_{0}^1+\int_{\frac{1}{\sqrt\kappa}}^{+\infty}
\right)\frac{d\xi}{\sqrt{(\xi^2-1)(1-\kappa\xi^2)}}\equiv 0\,,\qquad \forall\ 0<\kappa<1\, . \end{align}
This identity can  be  also directly checked, using proper changes of coordinate combined with cuts of the complex plane, in order to make the square roots single--valued in a neighbourhood of the real axis.\end{remark}

\noindent The advantage of looking at $g(\kappa)$ instead of $T_0(\kappa)$ is that the integration path  in~\eqref{ovlT} is $\kappa$--independent, and this turns to be useful when taking $\kappa$--derivatives. The main result at this respect in this section is the following

\begin{proposition}\label{period} \item[\tiny\textbullet] Let
 $\kappa\in {\mathbb R}\setminus \{0, 1\}$ and let $g(\kappa)$ be as in~\eqref{ovlT}.  There  exists two positive real numbers $\overline{\cal R}^*$, $\overline{\cal S}^*$ and 
two complex numbers \begin{align} \overline{\cal R}(\kappa),\ \overline{\cal S}(\kappa)\in \left\{\begin{array}{l}\displaystyle  {\mathbb R}_+\quad {\rm if}\quad \kappa<0\\
{\mathbb C}\quad {\rm if}\quad 0<\kappa<1\\
\displaystyle  {\rm i}{\mathbb R}_+\quad {\rm if}\quad \kappa>1
\end{array}\right.\end{align} with
\begin{align}\Re\overline{\cal R}(0)=\Re\overline{\cal S}(0)=1\, ,\quad 0\le \Re\overline{\cal R}(\kappa)\le \overline{\cal R}^*\, ,\quad 0\le \Re\overline{\cal S}(\kappa)\le \overline{\cal S}^*\quad \forall\ \kappa\in (-1, 1)\end{align}
 such that
\begin{align}\label{T0p-} g'(\kappa)=-\frac{\overline{\cal R}(\kappa)}{2\kappa}\qquad g''(\kappa)={+\frac{\overline{\cal S}(\kappa)}{4\kappa^2}}\qquad \forall\ \kappa\in {\mathbb R}\setminus\{0, 1\}\, . 
\end{align} 
\item[\tiny\textbullet] Let $\beta\ge 0$; $0<\kappa<1$, $j_\beta (\kappa)$ as in~\eqref{h}. There exist two positive numbers ${\cal R}_\beta^*$ and ${\cal S}_\beta^*\in {\mathbb R}$  and two real functions ${\cal R}_\beta(\kappa)$, ${\cal S}_\beta(\kappa)$ satisfying
\begin{align}\label{RS}
&{\cal R}_\beta(0)={\cal S}_\beta(0)=\left\{\begin{array}{l}\displaystyle  1\ {\rm if}\ \b=0\\\\
\displaystyle  0\ {\rm if}\ \b>0
\end{array}\right.\nonumber\\\nonumber\\
& 0\le {\cal R}_\beta(\kappa)\le {\cal R}_0^*\, ,\quad 0\le {\cal S}_\beta(\kappa)\le {\cal S}_0^*\quad\forall\ \beta\ge 0\,\quad \forall\ \kappa\in (0, 1)
\end{align}
such that
\begin{align} j_\beta '(\kappa)=-\frac{{\cal R}_\beta(\kappa)}{2\kappa}\qquad j_\beta ''(\kappa)={+\frac{{\cal S}_\beta(\kappa)}{4\kappa^2}}\qquad \forall\ 0<\kappa<1\, . 
\end{align}

\end{proposition}
{\bf Proof\ }
We prove the first statement. We distinguish two  cases.

\noindent Case 1:  $\kappa<0$ or $\kappa>1$. The integral takes real values when $\kappa<0$; purely imaginary ones when $\kappa>1$:
\begin{align}g(\kappa)=\left\{\begin{array}{l}\displaystyle  \int_{1}^{+\infty}\frac{d\xi}{\sqrt{(\xi^2-1)(1-\kappa\xi^2)}}\quad \kappa<0\\\\
\displaystyle  -{\rm i}\int_{1}^{+\infty}\frac{d\xi}{\sqrt{(\xi^2-1)(\kappa\xi^2-1)}}\quad \kappa>1
\end{array}\right.\end{align}
The function under the integral
 is bounded above by $\frac{1}{\min\{1, \sqrt{|\kappa|}\}\sqrt{\xi^4-1}}$ when $\kappa<0$;  by $\frac{1}{{\xi^2-1}}$ when $\kappa>1$. Both such bounds are integrable.   Then it is possible to derive under the integral, and we obtain
\begin{align}\label{T0p-}
g'(\kappa)=
\left\{\begin{array}{l}\displaystyle  \frac{1}{2}\int_{1}^{+\infty}\frac{\xi^2d\xi}{\sqrt{(\xi^2-1)(1-\kappa\xi^2)^3}}\quad \kappa<0\\\\
\displaystyle  \frac{\rm i}{2}\int_{1}^{+\infty}\frac{\xi^2 d\xi}{\sqrt{(\xi^2-1)(\kappa\xi^2-1)^3}}\quad \kappa>1
\end{array}\right.
\end{align}
and
\begin{align}
g''(\kappa)=
\left\{\begin{array}{l}\displaystyle  \frac{3}{4}\int_{1}^{+\infty}\frac{\xi^4d\xi}{\sqrt{(\xi^2-1)(1-\kappa\xi^2)^5}}\quad \kappa<0\\\\
\displaystyle  -\frac{3}{4}{\rm i}\int_{1}^{+\infty}\frac{\xi^4 d\xi}{\sqrt{(\xi^2-1)(\kappa\xi^2-1)^5}}\quad \kappa>1
\end{array}\right.
\end{align}
We  change  variable 
$1-\kappa\xi^2=\eta$ when $\kappa<0$, $
\kappa\xi^2-1=\eta$ when $\kappa>1$ and
 rewrite 
\begin{align}
g'(\kappa)\left\{\begin{array}{l}\displaystyle  \frac{1}{4|\kappa|}\int_{1-\kappa}^{+\infty}\sqrt{\frac{{\eta-1}}{{\left({\eta-1+\kappa}\right)\eta^3}}} \,d \eta\quad\ \kappa<0\\
\displaystyle  \frac{{\rm i}}{4|\kappa|}\int_{\kappa-1}^{+\infty}\sqrt{\frac{{\eta+1}}{{\left({\eta+1-\kappa}\right)\eta^3}}} \,d \eta\quad\ \kappa>1
\end{array}\right.\end{align}
and
\begin{align}
g''(\kappa)=\left\{\begin{array}{l}\displaystyle  \frac{3}{8|\kappa|^2}\int_{1-\kappa}^{+\infty}(\eta-1)\sqrt{\frac{{\eta-1}}{{\left({\eta-1+\kappa}\right)\eta^5}}} \,d \eta\quad\ \kappa<0\\
\displaystyle  -\frac{3}{8|\kappa|^2}{\rm i}\int_{\kappa-1}^{+\infty}(\eta+1)\sqrt{\frac{{\eta+1}}{{\left({\eta+1-\kappa}\right)\eta^5}}} \,d \eta\quad\ \kappa>1
\end{array}\right.\end{align}
so we take
\begin{align}
\overline{\cal R}(\kappa)=\left\{\begin{array}{l}\displaystyle  \frac{1}{2}\int_{1-\kappa}^{+\infty}\sqrt{\frac{{\eta-1}}{{\left({\eta-1+\kappa}\right)\eta^3}}} \,d \eta\quad\ \kappa<0\\
\displaystyle  \frac{\rm i}{2}\int_{\kappa-1}^{+\infty}\sqrt{\frac{{\eta+1}}{{\left({\eta+1-\kappa}\right)\eta^3}}} \,d \eta\quad\ \kappa>1
\end{array}\right.\end{align}
and
\begin{align}
\overline{\cal S}(\kappa)=\left\{\begin{array}{l}\displaystyle  \frac{3}{2}\int_{1-\kappa}^{+\infty}(\eta-1)\sqrt{\frac{{\eta-1}}{{\left({\eta-1+\kappa}\right)\eta^5}}} \,d \eta\quad\ \kappa<0\\
\displaystyle  -\frac{3}{2}{\rm i}\int_{\kappa-1}^{+\infty}(\eta+1)\sqrt{\frac{{\eta+1}}{{\left({\eta+1-\kappa}\right)\eta^5}}} \,d \eta\quad\ \kappa>1
\end{array}\right.\end{align}
Observe that, if $-1<\kappa<0$, \begin{align}\Re\overline{\cal R}(0^-)=1=\Re\overline{\cal S}(0^-)\end{align} and \begin{align}0\le \Re\overline{\cal R}(\kappa)=\frac{1}{2}\int_{1-\kappa}^{+\infty}\sqrt{\frac{{\eta-1}}{{\left({\eta-1+\kappa}\right)\eta^3}}} \,d \eta\le \frac{1}{2} \int_1^{+\infty}\sqrt{\frac{\eta-1}{(\eta-2)\eta^3}}d\eta\end{align}
\begin{align}0\le \Re\overline{\cal S}(\kappa)\le \frac{3}{2} \int_1^{+\infty}(\eta-1)\sqrt{\frac{\eta-1}{(\eta-2)\eta^5}}d\eta\, . \end{align}

\noindent Case 2:  $0<\kappa<1$. We split $g(\kappa)$ into its real and imaginary part. Using~\eqref{ovlTnew} and~\eqref{identity}, we obtain

\begin{align*}
g(\kappa)=&+\int_1^{\frac{1}{\sqrt\kappa}}\frac{d\xi}{\sqrt{(\xi^2-1)(1-\kappa\xi^2)}}+
\int_{\frac{1}{\sqrt\kappa}}^{+\infty}\frac{d\xi}{\sqrt{(\xi^2-1)(1-\kappa\xi^2)}}\nonumber\\
=&+\int_0^{\infty}\frac{dy}{\sqrt{(y^2+1)(1+\kappa y^2)}}+{\rm i}\int_0^1\frac{d\xi}{\sqrt{(1-\xi^2)(1-\kappa\xi^2)}}
\end{align*}
Notice that also in this case, the  functions under the integrals may be bounded by integrable functions: $\frac{1}{\sqrt{\kappa}(y^2+1)}$ for the former; $\frac{1}{\sqrt{1-\xi^2}}\frac{1}{\sqrt{1-\kappa}}$ in the latter. Again, we can derive under the integral, and obtain
\begin{align}
g'(\kappa)=-\frac{1}{2}\int_0^{\infty}\frac{y^2dy}{\sqrt{(y^2+1)(1+\kappa y^2)^3}}+\frac{\rm i}{2}\int_0^1\frac{\xi^2d\xi}{\sqrt{(1-\xi^2)(1-\kappa\xi^2)^3}}
\end{align}
and
\begin{align}
g''(\kappa)=+\frac{3}{4}\int_0^{\infty}\frac{y^4dy}{\sqrt{(y^2+1)(1+\kappa y^2)^5}}+\frac{3}{4}{\rm i}\int_0^1\frac{\xi^4d\xi}{\sqrt{(1-\xi^2)(1-\kappa\xi^2)^5}}
\end{align}
Then, letting $1+\kappa y^2=\eta$ in the first respective integrals, and $1-\kappa\xi^2=\eta$ in the second ones, 
\begin{align}
g'(\kappa)=-\frac{1}{4\kappa}\int_1^{+\infty}\sqrt{\frac{\eta-1}{(\eta-1+\kappa)\eta^3}}+\frac{\rm i}{4\kappa}\int_{1-\kappa}^{1}\sqrt{\frac{{1-\eta}}{{\left({\eta-1+\kappa}\right)\eta^3}}} \,d \eta\end{align}
and
\begin{align}
g''(\kappa)=+\frac{3}{8\kappa^2}\int_1^{+\infty}(\eta-1)\sqrt{\frac{\eta-1}{(\eta-1+\kappa)\eta^5}}+\frac{3}{8\kappa^2}{\rm i}\int_{1-\kappa}^{1}(1-\eta)\sqrt{\frac{{1-\eta}}{{\left({\eta-1+\kappa}\right)\eta^5}}} \,d \eta\end{align}
and we can take 
\begin{align}\overline{\cal R}(\kappa):=\frac{1}{2}\int_1^{+\infty}\sqrt{\frac{\eta-1}{(\eta-1+\kappa)\eta^3}}-\frac{\rm i}{2}\int_{1-\kappa}^{1}\sqrt{\frac{{1-\eta}}{{\left({\eta-1+\kappa}\right)\eta^3}}} \,d \eta\end{align} 
and
\begin{align}
\overline{\cal S}(\kappa)=\frac{3}{2}\int_1^{+\infty}(\eta-1)\sqrt{\frac{\eta-1}{(\eta-1+\kappa)\eta^5}}+\frac{3}{2}{\rm i}\int_{1-\kappa}^{1}(1-\eta)\sqrt{\frac{{1-\eta}}{{\left({\eta-1+\kappa}\right)\eta^5}}} \,d \eta\end{align}
Notice now that \begin{align}\Re\overline{\cal R}(0^+)=1=\Re\overline{\cal S}(0^+)\end{align} and \begin{align}0\le \Re \overline{\cal R}(\kappa)=\frac{1}{2}\int_1^{+\infty}\sqrt{\frac{\eta-1}{(\eta-1+\kappa)\eta^3}}\le \frac{1}{2}\int_1^{+\infty}\eta^{-\frac{3}{2}}=1\end{align}
and
\begin{align}0\le \Re \overline{\cal S}(\kappa)=\frac{3}{2}\int_1^{+\infty}(\eta-1)\sqrt{\frac{\eta-1}{(\eta-1+\kappa)\eta^5}}\le \frac{3}{2}\int_1^{+\infty}\eta^{-\frac{5}{2}}=1\end{align}
 for all $0<\kappa<1$. \\
The proof for  $j_\beta (\kappa)$ is completely analogous to the case 2 above (with the difference that we do not have the imaginary part in that case). One finds
\begin{align}{\cal R}_\beta(\kappa)=\frac{1}{2}\int_1^{+\infty}\frac{1}{1+\frac{\beta}{\kappa}(\eta-1)}\sqrt{\frac{\eta-1}{(\eta-1+\kappa)\eta^3}}\end{align} 
and
\begin{align}
{\cal S}_\beta(\kappa)=\frac{3}{2}\int_1^{+\infty}\frac{\eta-1}{1+\frac{\beta}{\kappa}(\eta-1)}\sqrt{\frac{\eta-1}{(\eta-1+\kappa)\eta^5}}
\end{align}
which verify~\eqref{RS}. $\qquad \square$

\section{Technicalities}\label{Technicalities}

In this section of the appendix we prove the bounds in~\eqref{bounds2},~\eqref{bounds3} and Lemma~\ref{zeroes}.

\paragraph{Proof of~\eqref{bounds2}} 
   
  We let
\begin{align}\label{ders}
&{\cal E}_*(A_*, {\rm r}_*):={\cal E}\circ\phi_{\rm aa}(A_*, {\rm r}_*)\ ,\quad {\cal E}_\circ(A, y):={\cal E}_*\circ\phi_{\rm rg}(A, y)={\cal E}_*(A, {\rm r}_\circ(A, y))\nonumber\\\nonumber\\
& {\cal B}_*(A_*, {\rm r}_*):={\cal B}\circ\phi_{\rm aa}(A_*, {\rm r}_*)\ ,\quad {\cal B}_\circ(A, y):={\cal B}_*\circ\phi_{\rm rg}(A, y)={\cal B}_*(A, {\rm r}_\circ(A, y))\nonumber\\\nonumber\\
&  T_{\rm p, *}(A_*, {\rm r}_*):=T_{\rm p}\circ\phi_{\rm aa}(A_*, {\rm r}_*)\ ,\quad T_{\rm p, \circ}(A, y):=T_{\rm p, *}\circ\phi_{\rm rg}(A, y)=T_{\rm p, *}(A, {\rm r}_\circ(A, y))\nonumber\\\nonumber\\
&  {\rm F}_\circ(A, y):={\rm F}_*\circ\phi_{\rm rg}(A, y)={\rm F}_*(A, {\rm r}_\circ(A, y))={\rm F}({\cal E}_\circ(A, y), {\rm r}_\circ(A, y))\nonumber\\\nonumber\\
&  {\rm F}_{*, 1, \circ}(A, y):={\rm F}_{*, 1}\circ\phi_{\rm rg}(A, y)={\rm F}_{*, 1}(A, {\rm r}_\circ(A, y))
\end{align}
(with ${\rm F}$, $T_{\rm p}$, ${\cal B}$ as in~\eqref{G*r*},~\eqref{tmax}--\eqref{f0T0},~\eqref{actionangle}) so as to write, more rapidly,
\begin{align}v(A, y; c)=e^{-y}\sqrt{2(c-\alpha F_\circ(A, y))}\ ,\quad \omega(A, y; c)=\alpha e^{-2y}{\rm F}_{*, 1, \circ}(A, y)\end{align}
and
\begin{align}\label{derF}
&\frac{1}{v}=\frac{e^y}{\sqrt{2(c-\alpha F_\circ(A, y))}}\ ,\quad \frac{\partial_A v}{v}=-\frac{\alpha}{2}\frac{\partial_AF_\circ(A, y) }{c-\alpha F_\circ(A, y) }\ ,\quad \frac{\partial_y v}{v}=-1-\frac{\alpha}{2}\frac{\partial_yF_\circ(A, y) }{c-\alpha F_\circ(A, y) }\ ,\nonumber\\
& \frac{\omega}{v}=\alpha\,e^{-y}\frac{{\rm F}_{*, 1, \circ}(A, y)}{\sqrt{2(c-\alpha F_\circ(A, y))}}\ ,\quad  \frac{\partial_A\omega}{v}=\alpha\,e^{-y}\frac{\partial_A{\rm F}_{*, 1, \circ}(A, y)}{\sqrt{2(c-\alpha F_\circ(A, y))}}\nonumber\\
& \frac{\partial_y\omega}{v}=-2\alpha\,e^{-y}\frac{{\rm F}_{*, 1, \circ}(A, y)}{\sqrt{2(c-\alpha F_\circ(A, y))}}+\alpha\,e^{-y}\frac{\partial_y{\rm F}_{*, 1, \circ}(A, y)}{\sqrt{2(c-\alpha F_\circ(A, y))}}
\end{align}
We evaluate the right hand sides of~\eqref{derF}, by means of  the chain  rule:
\begin{align}\label{derfcirc}
&{\rm F}_{*, 1, \circ}=\frac{{\rm F}_{\cal E}({\cal E}_\circ, {\rm r}_\circ)}{\hat T_{\rm p, \circ}}\ ,\quad \partial_A{\rm F}_{*, 1, \circ}=\frac{\partial^2_{\cal E}{\rm F}({\cal E}_\circ, {\rm r}_\circ)\partial_A{\cal E}_\circ+\partial^2_{{\cal E}{\rm r}}{\rm F}{\rm r}_{\rm s}'(A)
}{\hat T_{\rm p, \circ}}-\frac{
\partial_{\cal E} \hat T_{\rm p}\partial_A{\cal E}_\circ+
\partial_{\rm r} \hat T_{\rm p}{\rm r}_{\rm s}'(A)
}{\hat T_{\rm p, \circ}^2}{\rm F}_{\cal E}\nonumber\\\nonumber\\
&\partial_y{\rm F}_{*, 1, \circ}=\frac{\partial^2_{\cal E}{\rm F}({\cal E}_\circ, {\rm r}_\circ)\partial_y{\cal E}_\circ-e^{-y}\partial^2_{{\cal E}{\rm r}}{\rm F}
}{\hat T_{\rm p, \circ}}-\frac{
\partial_{\cal E} \hat T_{\rm p}\partial_y{\cal E}_\circ-e^{-y}
\partial_{\rm r} \hat T_{\rm p}
}{\hat T_{\rm p, \circ}^2}{\rm F}_{\cal E}\nonumber\\\nonumber\\
&\partial_A {\rm F}_\circ={\rm F}_{\cal E}({\cal E}_\circ, {\rm r}_\circ)\partial_A{\cal E}_\circ+{\rm F}_{\rm r}({\cal E}_\circ, {\rm r}_\circ){\rm r}_{\rm s}'(A)\ ,\quad \partial_y {\rm F}_\circ={\rm F}_{\cal E}({\cal E}_\circ, {\rm r}_\circ)\partial_y{\cal E}_\circ-e^{-y}{\rm F}_{\rm r}({\cal E}_\circ, {\rm r}_\circ)
\end{align}
where we have neglected to write the arguments (e.g., ${\rm F}_{\cal E}({\cal E}_\circ(A, y), {\rm r}_\circ(A, y))$, etc) and where, again by the chain\footnote{Use  $\partial_{A_*}{\cal E}_*=\frac{1}{\partial{\cal A}_{\cal E}}\circ\phi_{\rm aa}=\frac{1}{\hat T_{\rm p, *}(A_*, {\rm r}_*)}$ and $\partial_{{\rm r}_*}{\cal E}_*=-\frac{\partial{\cal A}_{\rm r}}{\partial{\cal A}_{\cal E}}\circ\phi_{\rm aa}=-{\cal B}_*(A_*, {\rm r}_*)$, implied by~\eqref{generating equations}.} rule,

\begin{align}
\partial_A{\cal E}_\circ=\frac{1}{\hat T_{\rm p, \circ}(A, y)}-
{\rm r}_{\rm s}'(A){\cal B}_\circ(A, y)\ ,\quad 
 \partial_y{\cal E}_\circ=
 e^{-y}{\cal B}_\circ(A, y)
\end{align}
As a result of the discussions in Sections~\ref{Regularity of the energy--time coordinates},~\ref{appB} and Appendix~\ref{elliptic integrals}, 
the functions  ${\rm F}$, $T_{\rm p}$ and ${\cal B}$ in~\eqref{ders} verify
\begin{align}
&
C'\log|\kappa|^{-1} \le |{\rm F}|,\  |T_{\rm p}|,\ |1/{\cal B}|\le C\log|\kappa|^{-1}\ ,\quad C'|\kappa|^{-1} \le |\partial_{{\cal E}, {\rm r}} {\rm F}|,\  |\partial_{{\cal E}, {\rm r}} T_{\rm p}|,\ |\partial_{{\cal E}, {\rm r}} {\cal B}|\le C|\kappa|^{-1}\nonumber\\
&C'|\kappa|^{-2} \le |\partial^2_{{\cal E}, {\rm r}} {\rm F}|,\  |\partial^2_{{\cal E}, {\rm r}} T_{\rm p}|,\ |\partial^2_{{\cal E}, {\rm r}} {\cal B}|\le C|\kappa|^{-2}
\end{align}
with $\kappa={\rm O}({\cal E}-{\rm r})={\rm O}(e^{-y})$ so that 
\begin{align}\label{FF*}
&C' L_-\le |{\rm F}_\circ|,\ |T_{\rm p, \circ}|,\ |{\cal B}_\circ|\le C L_+\nonumber\\
&C' e^{L_-}
\le | \partial_{{\cal E}, {\rm r}}{\rm F}({\cal E}_\circ, {\rm r}_\circ)|\, ,\ | \partial_{{\cal E}, {\rm r}}T_{\rm p}({\cal E}_\circ, {\rm r}_\circ)|\, ,\ | \partial_{{\cal E}, {\rm r}}{\cal B}({\cal E}_\circ, {\rm r}_\circ)|\le  C e^{L_+}\nonumber\\
&C' e^{2L_-}
\le | \partial^2_{{\cal E}, {\rm r}}{\rm F}({\cal E}_\circ, {\rm r}_\circ)|\, ,\ | \partial^2_{{\cal E}, {\rm r}}T_{\rm p}({\cal E}_\circ, {\rm r}_\circ)|\, ,\ | \partial^2_{{\cal E}, {\rm r}}{\cal B}({\cal E}_\circ, {\rm r}_\circ)|\le  C e^{2L_+}
\end{align}
  Finally,  using~\eqref{As}--\eqref{inverseA}, one has
\begin{align}{\rm r}_{\rm s}'(A)=\left.\frac{1}{{\cal A}'_{\rm s}({\rm r})}\right|_{{\rm r}={\rm r}_{\rm s}(A)}=\pi\,\sqrt{\frac{{\rm r}_{\rm s}(A)}{2-{\rm r}_{\rm s}(A)}}\end{align}
whence
\begin{align}\label{rspNEW}|{\rm r}'_{\rm s}(A)|\le \frac{C}{\sqrt{\varepsilon_-}}\end{align}
and collecting the bounds above into~\eqref{derF},
we find~\eqref{bounds2}.

\paragraph{Proof of~\eqref{bounds3}} 

We  use some results from Section~\ref{Regularity of the energy--time coordinates}. Taking in count~\eqref{changeofcoord},~\eqref{changeofcoord1},~\eqref{A(k)NEW} and~\eqref{sigmakappa} and
letting
\begin{align}
\sigma_*(A, {\rm r}_*)&:= \sigma\circ\phi_{\rm aa}(A, {\rm r}_*)\, ,\quad \kappa_*(A, {\rm r}_*):=\kappa\circ\phi_{\rm aa}(A, {\rm r}_*)\nonumber\\ 
\hat T_{\rm p, *}(A, {\rm r}_*)&:= \sigma_*(A, {\rm r}_*)\hat \tau_{\rm p, *}(A, {\rm r}_*):=\frac{T_{\rm p}(\kappa_*(A, {\rm r}_*))}{\pi}\, ,
\end{align}
we have that 

\begin{align}\label{G*}
{\rm G}_{*}(A, {\rm r}_*, \psi)&=\sigma_*(A, {\rm r}_*)\breve{\rm G}\big(\kappa_*(A, {\rm r}_*), \hat T_{\rm p, *}(A, {\rm r}_*) \psi\big)\\
\label{rho*}
\rho_{*}(A, {\rm r}_*, \psi)&=\frac{\sigma_*(A, {\rm r}_*)}{{\rm r}_*}\breve\rho\big(\kappa_*(A, {\rm r}_*), \hat T_{\rm p, *}(A, {\rm r}_*) \psi\big)
\end{align}
By the chain rule
\begin{align}\label{G*3}
{\rm G}_{*, 3}(A, {\rm r}_*, \psi)&=\partial_\psi{\rm G}_{*}(A, {\rm r}_*, \psi)\nonumber\\
&=\sigma_*(A, {\rm r}_*)\partial_{\psi}\breve{\rm G}(\kappa_*(A, {\rm r}_*), \hat T_{\rm p, *}(A, {\rm r}_*) \psi)\nonumber\\
&=\sigma_*(A, {\rm r}_*)\hat T_{\rm p, *}(A, {\rm r}_*)\breve{\rm G}_{3}(\kappa_*(A, {\rm r}_*), \hat T_{\rm p, *}(A, {\rm r}_*) \psi)
\end{align}
Similarly, 
\begin{align}\label{breverho3}
\rho_{*, 3}(A, {\rm r}_*, \psi)&=\frac{\sigma_*(A, {\rm r}_*)}{{\rm r}_*}\hat T_{\rm p, *}(A, {\rm r}_*)\breve\rho_{3}(\kappa_*(A, {\rm r}_*), \hat T_{\rm p, *}(A, {\rm r}_*) \psi)
\end{align}

\noindent
By the definitions in~\eqref{g}--\eqref{tildeP}, 
if
\begin{align}{\widetilde{\mathbb P}_{\varepsilon, \xi}}:=\bigcup_{(A, y)\in {\mathbb A}_{\varepsilon_-}\times {\mathbb Y}_\xi}\{A\}\times \{y\}\times V_\circ(A, y; \delta)\end{align}
then
\begin{align}\|\widetilde P_i\|_{{{\mathbb P}_{\varepsilon, \xi}}}\le \|P_i(A, y, \psi)\|_{{\widetilde{\mathbb P}_{\varepsilon, \xi}}}\end{align}
so we proceed to uniformly upper bound the $|P_i|$ in ${\widetilde{\mathbb P}_{\varepsilon, \xi}}$.  

\noindent
{\tiny\textbullet} By Proposition~\ref{Gderiv}, 
\begin{align}|\breve{\rm G}(\kappa, \theta)|,\ |\breve{\rm G}_3(\kappa, \theta)|\le C\end{align}

\noindent
{\tiny\textbullet} By
\eqref{G*},~\eqref{G*3} and~\eqref{FF*},
\begin{align}\label{G*NEW}|{\rm G}_*(A, {\rm r}_\circ(A, y), \psi)|\le C\sqrt{\varepsilon_+}\ ,\quad |{\rm G}_{*, 3}(A, {\rm r}_\circ(A, y), \psi)|\le C L_+\sqrt{\varepsilon_+}\end{align}

 \noindent 
{\tiny\textbullet}  Both the inequalities in~\eqref{G*NEW} hold (with the same proof) if ${\rm r}_\circ(A, y)$ is replaced by a generic ${\rm r}\in \Im{\rm r}_*(A, \cdot)$. Then,
\begin{align}|{\rm G}_{*, 1}(A, {\rm r}_\circ(A, y), \psi)|\le C\frac{\sqrt{\varepsilon_+}}{\varepsilon_-}\end{align}

\noindent
{\tiny\textbullet} Similarly, by~\eqref{rho*},
$|\rho_*(A, {\rm r}_*, \psi)|\le \sqrt{\varepsilon_+}$, hence
\begin{align}|\rho_{*, 1}(A, {\rm r}_\circ(A, y), \psi)|\le C\frac{\sqrt{\varepsilon_+}}{\varepsilon_-}\end{align}

\noindent
{\tiny\textbullet} The function ${\cal Y}(A, y, \psi; c)$  defined in~\eqref{Ycanc} verifies
 \begin{align}|{\cal Y}|\le C \sqrt{\alpha\,L_+}\end{align}
 having used the simplifying assumption~\eqref{Lm}.

\noindent
{\tiny\textbullet} By Lemma~\ref{zeroes},
\begin{align}|\rho_{*, 3}(A, {\rm r}_\circ(A, y), \psi)|\le C\,\sqrt{\varepsilon_+}\,\delta\end{align}

\noindent
{\tiny\textbullet} Recall~\eqref{rspNEW}.

\noindent
{\tiny\textbullet} Using the previous bounds into~\eqref{P} and writing the last term in the definition of $P_2$ as
\begin{align}e^{-y}
\frac{
\frac{({\rm C}-{\rm G}_*(A, {\rm r}_\circ(A, y), \psi))^2}{2{\rm r}_\circ(A, y)^2}-\frac{\beta }{{\rm r}_\circ(A, y)}
}{{\cal Y}(A, y, \psi; c)
+\sqrt{2\big(c-\alpha{\rm F}_*(A, {\rm r}_\circ(A, y))\big)}}\end{align}
we obtain, for $\|P_i\|_{\widetilde{\mathbb P}_{\varepsilon, \xi}}$, the  bounds at the right hand sides of~\eqref{bounds3}.

\paragraph{\bf Proof of Lemma~\ref{zeroes}}
Recall~\eqref{breverho3} and the expression of $\breve\rho_\theta(\kappa, \theta)$ in equation~\eqref{breverhotheta}. Equation
\begin{align}\label{rho3}\breve\rho_\theta(\kappa, \theta)={\breve{\rm G}(\kappa, \theta)^2}-{\cal A}(\kappa)=0\end{align}
has a unique solution
\begin{align}0<\theta_*(\kappa)<T_0(\kappa)\end{align}
if and only if
\begin{align}{\rm G}_0(\kappa)^2<{\cal A}(\kappa)<1\,.\end{align}
On the other hand, it is immediate to check that such inequality holds for all $0\ne \kappa<1$. Indeed, if
$0<\kappa<1$, then ${\rm G}_0(\kappa)^2=\kappa$ and we have
\begin{align}\kappa<{\cal A}(\kappa)=\frac{\int^{1}_{{\sqrt\kappa}}\frac{\xi^2d\xi}{{\sqrt{(1-\xi^2)(\xi^2-\kappa)}}}}{\int^{1}_{{\sqrt\kappa}}\frac{d\xi}{\sqrt{(1-\xi^2)(\xi^2-\kappa)}}}<1\, . \end{align}
If
 $\kappa<0$, then ${\rm G}_0(\kappa)^2=0$ and we have
\begin{align}0<{\cal A}(\kappa)=\frac{\int^{1}_{0}\frac{\xi^2d\xi}{{\sqrt{(1-\xi^2)(\xi^2-\kappa)}}}}{\int^{1}_{0}\frac{d\xi}{\sqrt{(1-\xi^2)(\xi^2-\kappa)}}}<1\end{align}
As a consequence of the formula~\eqref{rho3},  combined with the continuity of $\breve{\rm G}(\kappa, \cdot)$, we  find $V(\kappa; \delta)\subset (0, T_0(\kappa))$ (and $V'(\kappa; \delta)\subset (0, T_0(\kappa))$ when $\kappa<0$) such that
\begin{align}|\breve\rho_{3}(\kappa, \theta)|\le \frac{C \delta}{T_{\rm p}(\kappa)}\qquad \forall\ \theta\in V(\kappa; \delta)\quad \Big(\forall\ \theta\in V(\kappa; \delta)\cup V'(\kappa; \delta)\Big)\end{align}
which implies~\eqref{boundrho3}, after using~\eqref{breverho3}.
 $\qquad \square$

\addcontentsline{toc}{section}{References}
  \bibliographystyle{plain}
\footnotesize{\def\cprime{$'$} \def\cprime{$'$}

}
\vskip.1in
 \noindent
 {\it Declaration of interest: none.}

\end{document}